# Two lectures on the enumeration of curves by means of floor diagrams

**Abstract.** We discuss, following Mikhalkin, Brugallé, and many others, the counting of curves on toric surfaces with prescribed genus, Newton polygon, and intersection pattern with the toric boundary divisor, both at assigned and unassigned points.

The first lecture is dedicated to the proof of a correspondence theorem (for plane curves) with the counting of floor diagrams, using a degeneration of the projective plane to a chain of rational ruled surfaces. This is due to Brugallé and does not involve any tropical geometry.

The second lecture explores the relations with tropical geometry, and contains an introduction to toric varieties and tropical geometry. We discuss the correspondence theorem of Mikhalkin, and show how the corresponding tropical enumerative problem can be formulated in terms of the combinatorial problem of counting floor diagrams. We give many examples throughout, inspired by the study of the enumerative geometry of $K3$ surfaces, by degeneration to unions of rational surfaces with dual complex a tiling of the $\mathbf{S}^2$ sphere.

**Note.** These two lectures are part of a series[1] from the seminar *Degenerations and enumeration of curves on surfaces*, held at the University of Roma "Tor Vergata" in the years 2015–2017.

[1] see https://www.math.univ-toulouse.fr/~tdedieu/#enumTV

# Contents

# Lecture IX
# Enumeration of curves by floor diagrams, via degenerations of the projective plane

## by Thomas Dedieu

The goal of this text is to present a combinatorial translation of the same enumerative problem considered in [VI], namely that of counting plane curves of given degree and genus, with a prescribed intersection pattern with a fixed line, both at assigned and unassigned point, and passing through the appropriate number of fixed general points in the plane. The result is the Correspondence Theorem (1.9), which says that it is equivalent to count curves as above, and to count certain combinatorial gadgets called marked floor diagrams. Although, as is well-known, being a problem combinatorial doesn't make it easy, it is possible to teach a computer to count floor diagrams, whereas teaching it to directly solve the initial algebro-geometric problem is arguably hopeless.

The Correspondence Theorem has first been obtained by Brugallé and Mikhalkin [3], using tropical geometry. Here we present a proof which is purely algebro-geometric, due to Brugallé, see [2]. This proof is by degeneration of the projective plane. The idea in Brugallé and Mikhalkin's first proof, see [X], was to strech the general base points of the enumerative problem along a preferred direction in the plane; in Brugallé's proof presented here, this idea takes the following form. Let $r$ be the number of fixed general points in the enumerative problem. One lets the plane degenerate to a chain of $r$ rational ruled surfaces $\mathbf{F}_1$, in such a way that the fixed points of the problem end up one on each of the surfaces of the chain. This has the effect that the curves one has to count in the limit in order to solve the enumerative problem are so to speak maximally degenerate, in the sense that all their irreducible components are rational curves, which are moreover either sections or fibres of rational normal scrolls. One is thus left with objects from which all irrationality has vanished, so that it is possible to encode them in combinatorial terms.

The recursive formula obtained by Caporaso and Harris about the same enumerative problem, presented in the previous lecture [VI], may be interpreted in terms of a degeneration of the projective plane to a transverse union of two rational surfaces, see [VII]. The base cases of this recursion are the enumerative problems for degree one plane curves, i.e., lines. In Brugallé's proof, one somehow lets all the degenerations corresponding to the various iterations of the Caporaso–Harris' recursive procedure, all the way to the base cases, happen at once.

In order to directly follow Brugallé's proof, it is necessary to consider spaces of maps instead of the genuine Severi varieties, which parametrize curves embedded in surfaces. Indeed, multiple curves (or rather non-birational maps) inevitably show up. It follows that the results on degenerations of nodal curves along degenerations of surfaces usually considered in this volume, see [I] say, are not directly applicable. Following Brugallé, we use instead Jun Li's theory of relative stable maps [7, 8]. This theory will regrettably not be presented in this volume, and here I shall directly state its output without trying to give any explanation. Based on my personal experience, I believe that they should nevertheless look reasonable enough to anyone familiar with the results presented in due form in this volume. Finally, note that it should be possible (albeit certainly less convenient) to recast this proof using the theory expound in [I]; indications to that effect are given in Section 5.

In Section 1 we define floor diagrams and their markings, and state the Correspondence Theorem. In Section 2 we introduce the moduli spaces of maps that are needed for the proof of the Correspondence Theorem. Section 3 is the heart of the text, it contains the description of the degeneration (see Proposition (3.3)) that will finally give the proof of the Correspondence Theorem in Section 4, when we interpret it in combinatorial terms. Lastly, Section 5 is devoted to the interpretation of the considerations in this text in terms of the theory expound in [I].

**Thanks.** Many thanks to Erwan Brugallé for explaining all this to me and for answering unflustered my endless questions. Thanks also for your observations on this text. I know the secret meaning of the "French adage" you quote in [2].

**(0.1) Notation.** We use freely the notation introduced in [III, $**$]. In particular, $\underline{\mathbf{N}}$ denotes the set of sequences $\alpha = [\alpha_i]_{i\geqslant 1}$ of non-negative integers with finitely many $\alpha_i$ non-zero, and for all $\alpha \in \underline{\mathbf{N}}$, $|\alpha| = \sum_{i\geqslant 1} \alpha_i$ and $I\alpha = \sum_{i\geqslant 1} i\alpha_i$.

For any two integers $a \leqslant b$, we let $[\![a, b]\!] = [a, b] \cap \mathbf{Z}$.

In regard to the projective plane $\mathbf{P}^2$, we let $H$ denote the line class. In regard to the Hirzebruch ruled surface $\mathbf{F}_1$ (equivalently, $\mathbf{P}^2$ blown-up at one point), we let $E$ be the class of the $(-1)$-curve, $F$ be the class of the fibers of the ruling $\mathbf{F}_1 \to \mathbf{P}^1$, and $H = E + F$.

# 1 – Definitions and statement of the result

## 1.1 – Floor diagrams of plane curves relative to a line

In this section we define the floor diagrams corresponding to plane algebraic curves relative to a line. They are a particular instance of the floor diagrams considered in [X], and they are the only floor diagrams that will be considered in this text. Thus we simply call them floor diagrams, with no additional epithet.

**(1.1) A special class of weighted oriented graphs.** We shall consider connected, acyclic, oriented graphs $\mathcal{D}$, with finitely many vertices and edges. We denote by $\mathrm{Vert}(\mathcal{D})$ and $\mathrm{Edges}(\mathcal{D})$ the sets of vertices and edges of $\mathcal{D}$, respectively (or simply Vert and Edges if no confusion is likely). We see Edges as a subset of $\mathrm{Vert} \times \mathrm{Vert}$, with the first (resp. second) factor corresponding to the origin (resp. target) of an edge. We let $\mathcal{D}_{\mathrm{top}}$ be the topological incarnation of $\mathcal{D}$ as a 1-dimensional CW-complex.

We let the *vertices at infinity* of $\mathcal{D}$ be those 1-valent vertices of $\mathcal{D}$ for which the unique adjacent edge is outgoing. The other vertices of $\mathcal{D}$ will be called *floors* (this name will be explained in [X]). We denote by $\mathrm{Vert}^{-\infty}$ and $\mathrm{Vert}^{\mathrm{fin}}$ the sets of vertices at infinity and floors, respectively. By definition, there cannot be any edge joining two vertices at infinity. Thus, there is a partition $\mathrm{Edges} = \mathrm{Edges}^\infty \coprod \mathrm{Edges}^{\mathrm{fin}}$, where $\mathrm{Edges}^\infty$ is the set of edges originating from a vertex at infinity, and $\mathrm{Edges}^{\mathrm{fin}}$ is the set of those edges joining two floors. The former edges are called *infinite*, and the latter *finite*.

We will in addition assume that our graphs are weighted, which means that they come equipped with a function $w : \mathrm{Edges} \to \mathbf{N}^*$, which we call *weight function* One then defines the *divergence* $\mathrm{div}(v)$ of a vertex $v$ (either at infinity or a floor) as

$$\mathrm{div}(v) = \textstyle\sum_{(s,v)\in \mathrm{Edges}} w(s,v) - \sum_{(v,t)\in \mathrm{Edges}} w(v,t)$$

("sum of weights of incoming edges minus sum of weights of outgoing edges").

We adopt the following graphical convention: floors will be represented by ellipses, and vertices at infinity will not be represented; edges will be represented by vertical line segments, oriented from the bottom to the top. The weight of an edge is indicated only if it is greater than 1, as a circled integer.

**(1.2) Definition** (floor diagram)**.** *Let $d, g$ be nonnegative integers. A (plane) floor diagram of degree $d$ and genus $g$ is a connected, acyclic, weighted oriented graph of the class described in (1.1), such that the following conditions hold:*

*(i) the first Betti number $b_1(\mathcal{D}_{\mathrm{top}})$ equals $g$;*

*(ii) every floor has divergence $1$;*

*(iii) $\sum_{e\in \mathrm{Edges}^\infty} w(e) = -d$.*

**(1.3) Example.** All possible floor diagrams of degree 3 and genus 1 are pictured in Figure 1 below. When computing the genus, i.e., $b_1(\mathcal{D}_{\mathrm{top}})$, one should keep in mind that ellipses represent floors, which are vertices, hence the real ellipse does not contribute to $b_1$.[1] In the three diagrams below, the only cycle contributing to $b_1$ is that formed by the two lower floors and the two edges between them (it consists of two edges joining the two same vertices).

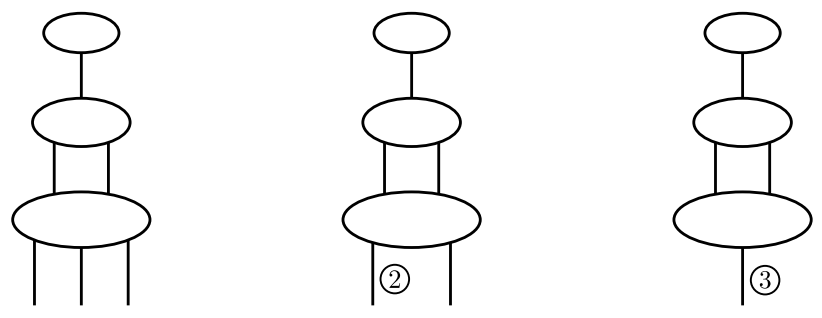

Figure 1: Floor diagrams of degree 3 and genus 1

**(1.4)** We let $\mathsf{D}$ be the abstract set $\mathrm{Vert}^{\mathrm{fin}} \coprod \mathrm{Edges}$, whose elements are all floors and edges. Since $\mathcal{D}$ is acylic, the set $\mathsf{D}$ inherits a partial ordering from the orientation.

Before we state the next definition, it is worth noting the following relations.

**(1.5) Lemma.** *Let $\mathcal{D}$ be a floor diagram of degree $d$ and genus $g$. The two following equalities hold:*

*(i)* $\mathrm{Card}(\mathrm{Vert}^{\mathrm{fin}}) = d$*;*

*(ii)* $\mathrm{Card}(\mathsf{D}) = 2d + g - 1 + \mathrm{Card}(\mathrm{Vert}^{-\infty})$*.*

*Proof.* One has

$$\begin{aligned} 0 = \textstyle\sum_{v\in \mathrm{Vert}} \mathrm{div}(v) &= \textstyle\sum_{v\in \mathrm{Vert}^{\mathrm{fin}}} \mathrm{div}(v) + \sum_{v \in \mathrm{Vert}^{-\infty}} \mathrm{div}(v) \\ &= \mathrm{Card}(\mathrm{Vert}^{\mathrm{fin}}) - d, \end{aligned}$$

[1] To make this a rational choice (pun intended), think of the ellipses as rational curves (after all ellipses are conics!), hence as having genus 0.

which proves (i). The Euler characteristic of $\mathcal{D}_{\text{top}}$ is

$$1 - g = \text{Card(Vert)} - \text{Card(Edges)},$$

so

$$\begin{aligned}\text{Card}(\mathsf{D}) &= \text{Card}(\text{Vert}^{\text{fin}}) + \big(g - 1 + \text{Card(Vert)}\big) \\ &= 2\,\text{Card}(\text{Vert}^{\text{fin}}) + g - 1 + \text{Card}(\text{Vert}^{-\infty}),\end{aligned}$$

which together with (i) proves (ii). □

**(1.6) Definition** (marking). *Let $\mathcal{D}$ be a floor diagram of degree $d$ and genus $g$, and $\alpha, \beta \in \underline{\mathbf{N}}^2$, such that $I\alpha + I\beta = d$. A* marking of type $(\alpha, \beta)$ *of $\mathcal{D}$ is a strictly increasing function*

$$m : \big\{-|\alpha| + 1, \ldots, 1, \ldots, 2d + g - 1 + |\beta|\big\} \to \mathsf{D}$$

*such that:*

*(i)* $\forall k \in \mathbf{N}^*$, $\forall i \in \llbracket -|\alpha| + (\sum_{j\leqslant k-1} \alpha_j) + 1,\ -|\alpha| + \sum_{j\leqslant k} \alpha_j \rrbracket$, $m(i) \in \text{Edges}^\infty$ *and* $w(m(i)) = k$*;*

*(ii)* $\forall k \in \mathbf{N}^*$*, there are exactly $\beta_k$ infinite edges of weight $k$ in* $m\big(\llbracket 1, 2d + g - 1 + |\beta| \rrbracket\big)$.

A floor diagram $\mathcal{D}$ enhanced with a marking $m$ of type $(\alpha, \beta)$ is called a *marked floor diagram* $(\mathcal{D}, m)$ of type $(\alpha, \beta)$.

Two marked floor diagrams of the same type are *equivalent* if there exists an isomorphism of oriented graphs $\phi : \mathcal{D} \cong \mathcal{D}'$ such that $w = w' \circ \phi$ and $m = \phi \circ m'$.

Markings will be indicated on graphic representations by plain numbers. Beware not to confuse them with weights, which are circled integers.

**(1.7) Examples.** Up to equivalence, there is a unique marked floor diagram of degree 3, genus 1, and type $(0, 3)$.

*(1.7.1) Marked floor diagrams of degree* 3 *and genus* 0 *of type* $(0, 3)$. Below I give all possibilities for such diagrams. For the leftmost diagram there is a unique marking (up to equivalence). For the next one, there are 5: the label $i$ may be chosen arbitrarily in $\llbracket 1, 5 \rrbracket$, and this choice determines uniquely the four labels not indicated on the picture. For the rightmost diagram, there are 3: the index $i$ may be chosen arbitrarily in $\llbracket 6, 8 \rrbracket$, and then the two remaining labels are uniquely determined.

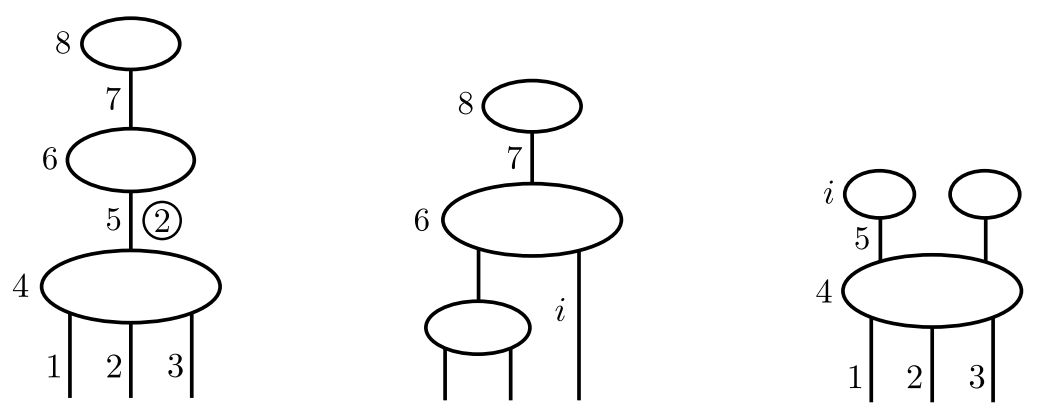

Figure 2: Marked floor diagrams of degree 3 and genus 0 of type $(0, 3)$.

*(1.7.2) Marked floor diagrams of degree* $3$ *and genus* $0$ *of type* $(0,[1,1])$. All possible marked diagrams of this type are displayed below. The leftmost diagram has two possible markings, determined by the choice of $i \in [\![1,2]\!]$. The central diagram has four possible markings, determined by the choice of $i \in [\![1,4]\!]$. The rightmost diagram has six possible markings, determined first by the choice of $i \in [\![1,2]\!]$, which fixes $i'$, and then by the choice of $j \in [\![5,7]\!]$ which fixes the rest.

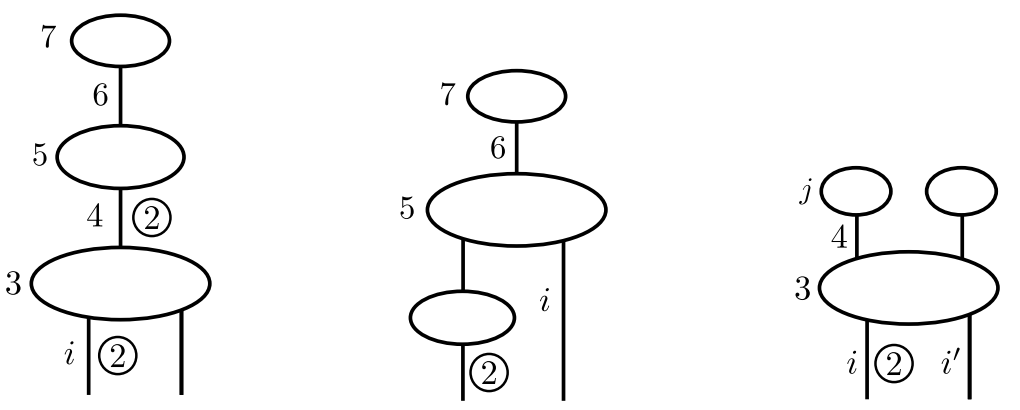

Figure 3: Marked floor diagrams of degree 3 and genus 0 of type $(0,[1,1])$.

*(1.7.3) Marked floor diagrams of degree* $3$ *and genus* $0$ *of type* $([0,1],1)$. All possible marked diagrams of this type are displayed below. It is instructive to make comparisons with the type $(0,[1,1])$ considered in the previous paragraph (1.7.2). The leftmost diagram has a unique possible marking. The central diagram has three possible markings, determined by the choice of $i \in [\![1,3]\!]$. The rightmost diagram has three possible markings, determined by the choice of $j \in [\![4,6]\!]$.

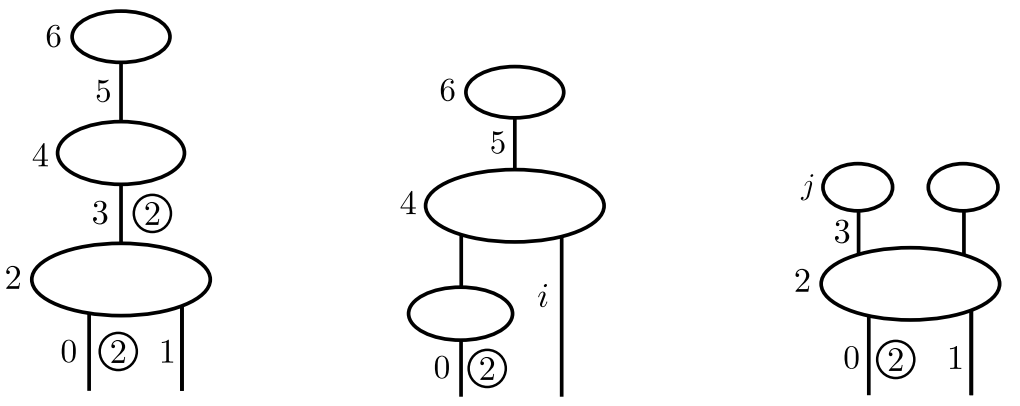

Figure 4: Marked floor diagrams of degree 3 and genus 0 of type $([0,1],1)$.

## 1.2 – The enumerative problem

We shall now state precisely the enumerative problem we are interested in. The definitions below are specializations of those given in [III].

Assume a line $\Lambda \subseteq \mathbf{P}^2$ is fixed. Let $d > 0$ and $g \geqslant 0$ be integers, and $\alpha, \beta \in \underline{\mathbf{N}}$, such that $I\alpha + I\beta = d$. Let $\Omega$ be a general set of $\alpha$ points on $\Lambda$. We let $V_g^{d,\mathrm{irr}}(\alpha,\beta)(\Omega)$ be the Zariski closure in the projective space $|\mathcal{O}_{\mathbf{P}^2}(d)|$ of the localy closed set $\mathring{V}_g^{d,\mathrm{irr}}(\alpha,\beta)(\Omega)$ parametrizing irreducible curves $C \subseteq \mathbf{P}^2$ of degree $d$, genus $g$, and having intersection pattern with $\Lambda$ prescribed by $\alpha, \beta$ and $\Omega$ as follows: write $\Omega = \big(\{p_{i,j}\}_{1 \leqslant j \leqslant \alpha_i}\big)_{i \geqslant 1}$, and let $\nu : \bar{C} \to C$ be the normalization of $C$; we ask that there exist $\alpha$ points $q_{i,j} \in C$, $1 \leqslant j \leqslant \alpha_i$, and $\beta$ points $r_{i,j} \in C$, $1 \leqslant j \leqslant \beta_i$, such that

$$\forall\, 1 \leqslant j \leqslant \alpha_i : \quad \phi(q_{i,j}) = p_{i,j} \qquad \text{and} \qquad \phi^*\Lambda = \sum_{1 \leqslant j \leqslant \alpha_i} i\, q_{i,j} + \sum_{1 \leqslant j \leqslant \beta_i} i\, r_{i,j}.$$

The variety $V_g^{d,\mathrm{irr}}(\alpha,\beta)(\Omega)$ has dimension $2d + g - 1 + |\beta|$. We call it a logarithmic Severi variety of irreducible[2] plane curves, relative to the line $\Lambda$.

[2] in the definition given in [III], the curve $C$ was not required to be irreducible; the superscript 'irr' indicates this additional requirement.

We let $N_g^{d,\mathrm{irr}}(\alpha,\beta)$ be the degree of $V_g^{d,\mathrm{irr}}(\alpha,\beta)(\Omega)$ as a projective subvariety of $|\mathcal{O}_{\mathbf{P}^2}(d)|$; it is independent of the choice of $\Lambda$ and $\Omega$. This is the number of plane curves of degree $d$, genus $g$, and having contact of the required type with $\Lambda$, that pass through $2d+g-1+|\beta|$ general points on $\mathbf{P}^2$.

## 1.3 – The Correspondence Theorem

The Correspondence Theorem (1.9) is the central result of this text. It will be proved in the subsequent sections. It states that one can compute the numbers $N_g^{d,\mathrm{irr}}(\alpha,\beta)$ by enumerating floor diagrams. To do so, one needs to count the latter with the appropriate multiplicity.

**(1.8) Definition.** *Let $(\mathcal{D},m)$ be a marked floor diagram of type $(\alpha,\beta)$. Its multiplicity is*

$$\mu^{\mathbf{C}}(\mathcal{D},m) = I^\beta \prod_{e\in \mathrm{Edges}^{\mathrm{fin}}} w(e)^2.$$

Thus, the multiplicity of a marked floor diagram is the product of the multiplicities of its edges where:
(a) a finite edge of weight $w$ has multiplicity $w^2$;
(b) an infinite edge of weight $w$ corresponding to a mobile boundary condition[3] (i.e., such that its label in the marking is positive) has multiplicity $w$;
(c) an infinite edge of weight $w$ corresponding to a fixed boundary condition (i.e., such that its label in the marking is non-positive) has multiplicity 1.

One may simply write $\mu^{\mathbf{C}}(\mathcal{D})$ in place of $\mu^{\mathbf{C}}(\mathcal{D},m)$ since the multiplicity only depends on the type and not on the marking. We leave the original superscript '**C**' as a reminder that this is the correct multiplicity only in order to count *complex* plane curves.

**(1.9) Theorem** (Brugallé–Mikhalkin [3], see also [1])**.** *Let $d>0$ and $g\geqslant 0$ be integers, and $\alpha,\beta\in\underline{\mathbf{N}}$, such that $I\alpha+I\beta=d$. One has*

$$N_g^{d,\mathrm{irr}}(\alpha,\beta) = \sum_{(\mathcal{D},m)} \mu^{\mathbf{C}}(\mathcal{D},m),$$

*where the sum on the right-hand-side is taken over all equivalence classes of marked floor diagrams of degree $d$, genus $g$, and type $(\alpha,\beta)$.*

**(1.10) Example.** Let us compute the numbers $N_0^{3,\mathrm{irr}}(\alpha,\beta)$ in the cases considered in Example (1.7), by enumerating floor diagrams.

For $(\alpha,\beta)=(0,3)$, the diagram of the leftmost shape has multiplicity 4, and all the other have multiplicity 1. Thus

$$N_0^{3,\mathrm{irr}}(0,3) = 1\times 4+5+3 = 12$$

(the first number of each summand indicates the number of possible markings for a diagram of the given shape, and the second one indicates the multiplicity; it is omitted if it is one).

For $(\alpha,\beta)=(0,[1,1])$, the diagrams of the leftmostshape have multiplicity 8, and all the other have multiplicity 2. Thus

$$N_0^{3,\mathrm{irr}}(0,[1,1]) = 2\times 8+4\times 2+6\times 2 = 36.$$

For $(\alpha,\beta)=([0,1],1)$, the diagram of the leftmost shape has multiplicity 4, and all the other have multiplicity 1. Thus

$$N_0^{3,\mathrm{irr}}(0,[1,1]) = 1\times 4+3+3 = 10.$$

[3] the term 'boundary' comes from the definition (2.2), see Paragraph (4.1).

## 2 – Log-GW invariants of $\mathbf{P}^2$ and Hirzebruch surfaces

Despite the definition of logarithmic Severi varieties having been given in terms of embedded curves in Section 1.2, it will be necessary in the context of these notes to consider their counterparts as spaces of stable maps. Indeed we will need to consider stable maps that are not birational. Moreover, we will degenerate $\mathbf{P}^2$ to the unions of $\mathbf{P}^2$ and various copies of $\mathbf{F}_1$, and therefore we need to define these spaces not only for the projective plane. The results stated in this section initially appeared in [4] and [10].

**(2.1) Definition.** *Let $X$ be a smooth variety. A* stable map *to $X$ is a map $f : C \to X$ with source $C$ a nodal curve, such that $C$ has only finitely many automorphisms commuting with $f$.*

The finiteness of automorphisms condition boils down to the following requirement: for all irreducible components $C_i$ of $C$ such that $p_a(C_i) \leqslant 1$ and $C_i$ is contracted by $f$, if $p_a(C_i) = 0$ (resp. $p_a(C_i) = 1$), then $C_i$ has at least three (resp. one) intersection points with the other components of $C$.

**(2.2) Stable maps relative to a curve.** Let $S$ be a smooth projective surface, and $R \subseteq S$ be a fixed smooth curve (not necessarily connected). The curve $R$ is often called the *boundary* of the pair $(S, R)$. Let $g \in \mathbf{Z}_{\geqslant 0}$, $\xi \in \mathrm{NS}(S)$, and $\alpha, \beta \in \underline{\mathbf{N}}$, such that $I\alpha + I\beta = \xi \cdot R$. Let $n \in \mathbf{N}$. Let $\Omega = \big(\{p_{i,j}\}_{1\leqslant j\leqslant \alpha_i}\big)_{i\geqslant 1}$ be a set of $\alpha$ general points on $R$. We define $M_g^{\xi,\mathrm{con}}(\alpha,\beta)(\Omega)$ as the moduli space of stable maps $f : C \to S$ such that:

(i) $C$ is a connected projective curve of arithmetic genus $g$;
(ii) $f_*[C] = \xi \in \mathrm{H}^2(S, \mathbf{Z})$;
(iii) $f(C)$ and $R$ have no irreducible component in common;
(iv) there exist $\alpha$ points $q_{i,j} \in C$, $1 \leqslant j \leqslant \alpha_i$, and $\beta$ points $r_{i,j} \in C$, $1 \leqslant j \leqslant \beta_i$, such that

$$f(q_{i,j}) = p_{i,j} \quad \forall\, 1 \leqslant j \leqslant \alpha_i, \qquad \text{and} \qquad f^*R = \sum_{1\leqslant j\leqslant \alpha_i} i\, q_{i,j} + \sum_{1\leqslant j\leqslant \beta_i} i\, r_{i,j}.$$

If $x = \{x_1, \ldots, x_n\}$ is a set of $n$ general points on $S$, we define $M_g^{\xi,\mathrm{con}}(\alpha,\beta)(\Omega, x)$ as the moduli subspace of $M_g^{\xi,\mathrm{con}}(\alpha,\beta)(\Omega)$ determined by the further requirement:

(v) $f(C) \supset x$ as sets.

Similarly, for all components $M$ of $M_g^{\xi,\mathrm{con}}(\alpha,\beta)(\Omega)$, $M(x)$ is the subspace of $M$ determined by condition (v) above.

The superscript 'con' stands for 'connected'; note that $C$ is required to be connected but may be reducible. We do not require that $f$ be birational.

An isomorphism between two stable maps $f_1 : C_1 \to X$ and $f_2 : C_2 \to X$ is an isomorphism $\phi : C_1 \to C_2$ such that $f_1 = f_2 \circ \phi$. The existence of moduli spaces of stable maps is a highly non-trivial result. We refer to [5] for an exposition in the case of plain stable maps; the result for relative stable maps (i.e., stable maps with additional features as in (2.2)) is due to Jun Li, see [7, 9].

**(2.3)** We shall consider instances of the space $M_g^{\xi,\mathrm{irr}}(\alpha,\beta)(\Omega)$ in the two following situations.

*(2.3.1)* $S = \mathbf{P}^2$, $\xi = dH$ with for some integer $d > 0$, and $R = \Lambda$ a line; we may then write $M_g^{d,\mathrm{con}}(\alpha,\beta)(\Omega)$ or even $M_g^{d,\mathrm{con}}(\alpha,\beta)$ in place of $M_g^{\xi,\mathrm{con}}(\alpha,\beta)(\Omega)$.

*(2.3.2)* $S = \mathbf{F}_1$, $\xi = \lambda H + kF$ for some nonnegative integers $\lambda, k$, and $R = \Lambda_E \sqcup \Lambda_H$, the disjoint union of two curves $\Lambda_E \in |E|$ and $\Lambda_H \in |H|$ (recall Notation (0.1)). In this case it is better to impose directly contact conditions with $R$ that are decomposed in two disjoint set of contact conditions $(\alpha_E, \beta_E, \Omega_E)$ and $(\alpha_H, \beta_H, \Omega_H)$ with $\Lambda_E$ and $\Lambda_H$ respectively. One thus obtains a set $M_g^{\xi,\mathrm{con}}(\alpha_E, \alpha_H, \beta_E, \beta_H)(\Omega_E, \Omega_H)$ of (isomorphism classes of) stable maps; we may abuse the notation of (2.2) and still denote this space by $M_g^{\xi,\mathrm{con}}(\alpha, \beta)(\Omega)$, with $\alpha = (\alpha_E, \alpha_H)$, $I\alpha = I\alpha_E + I\alpha_H$ and $|\alpha| = |\alpha_E| + |\alpha_H|$, and similarly for $\beta = (\beta_E, \beta_H)$. One has

$$I\alpha_E + I\beta_E = \xi \cdot E = k, \qquad I\alpha_H + I\beta_H = \xi \cdot H = \lambda + k,$$

and $-(K_S + R) = 2F$.

We shall be interested only in those irreducible components $M$ of $M_g^{\xi,\mathrm{con}}(\alpha, \beta)(\Omega)$ that are *enumeratively relevant*, i.e., such that for a set $x$ of $n = \dim(M)$ general points on $S$, the subspace $M(x)$ of those stable maps in $M$ whose image contains $x$ is not empty. The other irreducible components of $M$ are said to be *enumeratively irrelevant*; their enumerative contribution is zero, and therefore they will always be ignored. Note that if the general member of $M$ is birational onto its image, then $M$ is enumeratively relevant.

**(2.4) Theorem.** *Assume we are either in situation (2.3.1) above, or in situation (2.3.2) and $\lambda > 0$. Let $M$ be an irreducible component of $M_g^{\xi,\mathrm{con}}(\alpha, \beta)(\Omega)$. The family of those curves in $S$ that are the image of a stable map in $M$ has dimension at most*

$$-(K_S + R) \cdot \xi + g - 1 + |\beta| = \begin{cases} 2d + g - 1 + |\beta| & \text{in situation (2.3.1)} \\ 2\lambda + g - 1 + |\beta_E| + |\beta_H| & \text{in situation (2.3.2).} \end{cases}$$

*The component $M$ is enumeratively relevant if and only if equality holds, and in this case its general member $[f : C \to S]$ enjoys the following properties:*
*(i) the curve $C$ is smooth (hence irreducible);*
*(ii) $f$ is an immersion, birational onto its image;*
*(iii) $f(C)$ intersects $R$ at smooth points;*
*(iv) for all finite subset $y \subseteq S$, if $f$ is general with respect to $y$, then $f(C)$ does not intersect $y$.*

In situation (2.3.1), the statement on the dimension and birationality follow from [III, Prop. **], and the rest from [III, Thm. **]. In situation (2.3.2) the proof is similar; the first part is proven in [10, Prop. 3.4], and the rest again follows from [III, Thm. **].

**(2.5)** In the cases not covered by (2.4), i.e., in situation (2.3.2) when $\xi = kF$ $(k > 0)$, let $M$ be an irreducible component of $M_g^{kF,\mathrm{con}}(\alpha, \beta)(\Omega)$. The family of those curves in $\mathbf{F}_1$ that are the image of a stable map in $M$ has dimension 0 or 1, always less or equal than

$$-(K_S + R) \cdot \xi + g - 1 + |\beta| = g - 1 + |\beta_E| + |\beta_H|,$$

and $M$ is enumeratively relevant if and only if equality holds. The space $M_g^{kF,\mathrm{con}}(\alpha, \beta)(\Omega)$ has an enumeratively relevant irreducible component if and only if $g = 0$, and either

$$\alpha_E = \beta_H = [0, \ldots, 0, 1] \quad \text{and} \quad \beta_E = \alpha_H = 0 \quad \text{or viceversa,}$$

or

$$\beta_E = \beta_H = [0, \ldots, 0, 1] \quad \text{and} \quad \alpha_E = \alpha_H = 0.$$

In those cases there is a unique enumeratively relevant component, which parametrizes degree $k$ maps from $\mathbf{P}^1$ to a fibre $F$ of $\mathbf{F}_1$, with total ramification over the two points $F \cap \Lambda_E$ and $F \cap \Lambda_H$.

This is a consequence of the basic observation that a connected member of the linear system $|kF|$ necessarily consists of one single fibre $F$ with multiplicity $k$, hence is uniquely determined by its passing through one single point $x_1 \in \mathbf{F}_1$. When $g > 0$, irreducible components of $M_g^{kF,\mathrm{con}}(\alpha,\beta)(\Omega)$ even have dimension strictly larger than the expected $g-1+|\beta|$, which itself is already strictly larger than the dimension of the corresponding families of curves in $\mathbf{F}_1$.

Note that condition (iv) of Thm. (2.4) holds for an enumeratively relevant component of $M_0^{kF,\mathrm{con}}(0,\beta)(\Omega)$.

**(2.6) Log-GW invariants.** In favourable circumstances (see, e.g., [10, Subsec. 2.3]), and in particular in the situations we are concerned with in this text, one may define the *logarithmic Gromov–Witten invariant* (log-GW invariant for short) of the pair $(S,R)$ in genus $g$ and class $\xi$, with boundary conditions $(\alpha,\beta)$ as the length, for a general choice of a set $\Omega$ of $\alpha$ points on $R$ and of a set $x$ of

$$n = -(K_S+R)\cdot\xi + g - 1 + |\beta|$$

points on $S$, of the 0-dimensional space $M_g^{\xi,\mathrm{con}}(\alpha,\beta)(\Omega,x)$. The moduli space $M_g^{\xi,\mathrm{con}}(\alpha,\beta)(\Omega)$ is endowed with the structure of a stack, and each point $[f : C \to S] \in M_g^{\xi,\mathrm{con}}(\alpha,\beta)(\Omega,x)$ should be counted with the multiplicity

$$\mu(f,x) = \frac{1}{\mathrm{Card}\,(\mathrm{Aut}\,f)} \prod_{1\leqslant i\leqslant r} \mathrm{Card}\, f^{-1}(x_i).$$

**(2.7)** In situations (2.3.1) and (2.3.2), it follows from (2.4) and (2.5) that the space $M_g^{\xi,\mathrm{con}}(\alpha,\beta)(\Omega,x)$ is finite, and only those enumeratively relevant components of $M_g^{\xi,\mathrm{con}}(\alpha,\beta)(\Omega)$ contribute to the corresponding log-GW invariant. It holds also true that the cardinality of $M_g^{\xi,\mathrm{con}}(\alpha,\beta)(\Omega,x)$ is independent of the choice of $\Omega$ and $x$ (as long as they are general). Moreover, conditions (ii) and (iii) in (2.4) together imply that $\mu(f,x)=1$ for all $[f : C \to S] \in M_g^{\xi,\mathrm{con}}(\alpha,\beta)(\Omega,x)$. Correspondingly, in situation (2.3.1) say, there is a 1 : 1 correspondence between $M_g^{\xi,\mathrm{con}}(\alpha,\beta)(\Omega,x)$ and $V_g^{d,\mathrm{irr}}(\alpha,\beta)(\Omega) \cap \bigcap_{1\leqslant i\leqslant r} x_i^\perp$, where $x_i^\perp$ denotes the hyperplane of $|\mathcal{O}_{\mathbf{P}^2}(d)|$ consisting of those divisors that contain $x_i$, so that

$$\sum_{f\in M_g^{\xi,\mathrm{con}}(\alpha,\beta)(\Omega,x)} \mu(f,x) = N_g^{d,\mathrm{irr}}(\alpha,\beta).$$

A similar result holds in situation (2.3.2) with $\lambda > 0$, see [10].

**(2.8) Enumeratively relevant components of dimension at most 1.** It will be useful in the proof of Proposition (3.3) to have a list of the enumeratively relevant components of dimension at most 1 of $M_g^{\xi,\mathrm{con}}(\alpha,\beta)(\Omega)$ in situation (2.3.2), i.e., for the ruled surface $\mathbf{F}_1$. For such a component, we must have $\lambda \leqslant 1$. If $\lambda = 0$, we are in the situation of (2.5): enumerative relevance imposes $g = 0$, and all enumeratively relevant components have dimension at most 1. The corresponding log-GW invariants are computed in (2.9) below. If $\lambda = 1$, we must have $g = 0$ and $\beta = 0$. Then the space $M_0^{H+kF,\mathrm{irr}}(\alpha,0)(\Omega)$ is isomorphic to the sublinear system of $|H+kF|$ defined by the contact conditions imposed at $\Omega$: indeed the arithmetic genus of $H+kF$ is 0, and the contact conditions at fixed points are linear. Therefore for a general point $x \in \mathbf{F}_1$, $M_0^{H+kF,\mathrm{irr}}(\alpha,0)(\Omega,x)$ has exactly one member, which has multiplicity 1.

Note moreover that in situation (2.3.1), i.e., for the projective plane, there is no enumeratively relevant components of of dimension 0 of $M_g^{\xi,\mathrm{con}}(\alpha,\beta)(\Omega)$, as follows from similar considerations.

**(2.9) Multiplicity for multiple fibres.** Now let us compute the log-GW invariants in the cases considered in (2.5) above, i.e., we consider the spaces $M_0^{kF,\mathrm{irr}}(\alpha,\beta)$ for the pair $(\mathbf{F}_1, \Lambda_E +$ $\Lambda_H$), and either (a) $\alpha_E = \beta_H = [0, \ldots, 0, 1]$ and $\beta_E = \alpha_H = 0$ (or viceversa), or (b) $\beta_E = \beta_H = [0, \ldots, 0, 1]$ and $\alpha_E = \alpha_H = 0$.

In case (a) $\Omega$ consists of one point, and $M_0^{kF,\mathrm{irr}}(\alpha, \beta)(\Omega)$ has a unique enumeratively meaningful component, of dimension 0, and in case (b) $\Omega$ is empty, the unique enumeratively meaningful component of $M_0^{kF,\mathrm{irr}}(\alpha, \beta)$ has dimension 1, hence one should impose the passing through 1 point $x$. In either case there is a unique stable map to consider, i.e., $M_0^{kF,\mathrm{irr}}(\alpha, \beta)(\Omega)$ in case (a) and $M_0^{kF,\mathrm{irr}}(\alpha, \beta)(x)$ in case (b) are supported on exactly one point, and this map is a degree $k$ cover of the form $z \in \mathbf{P}^1 \mapsto z^k \in \mathbf{P}^1$ where the target $\mathbf{P}^1$ identifies with the relevant fibre of $\mathbf{F}_1$, and $f(0)$ and $f(\infty)$ are on $\Lambda_E$ and $\Lambda_H$ respectively. The automorphism group of this map is thus isomorphic to $\mathbf{Z}/k\mathbf{Z}$.

On the other hand, there is no point $x$ in case (a), and the unique point $x$ in case (b) has $k$ different points in its preimage by $f$. Therefore the log-GW invariants are $\frac{1}{k}$ in case (a), and 1 in case (b).

## 3 – Degeneration of $\mathbf{P}^2$ to a chain of ruled surfaces of maximal length

This section contains the heart of Brugallé's proof of the Correspondence Theorem. We introduce the necessary degeneration of the projective plane, and provide (in Proposition (3.3)) the description of the limit enumerative problem that will be encoded in terms of floor diagrams in the next section.

**(3.1) Degeneration of the plane.** We shall construct a sequence $\pi_r : \mathcal{P}_r \to \mathbf{D}$, indexed by $r \in \mathbf{N}$, of degenerations over the unit disk $\mathbf{D}$ of the projective plane together with a distinguished line, with central fibre $\pi_r^{-1}(0)$ a reduced chain of rational surfaces

$$S_{r+1} \cup S_r \cup \cdots \cup S_1 \cup S_0$$

with $S_{r+1} = \mathbf{P}^2$ and $S_i = \mathbf{F}_1$ for $0 \leqslant i \leqslant r$: for all $i = 1, \ldots, r$, the surface $S_i$ intersects only $S_{i+1}$ and $S_{i-1}$, along two smooth rational curves, respectively, while $S_{r+1}$ intersects only $S_r$ and $S_0$ intersects only $S_1$, along one smooth rational curve. We set $\Lambda_i = S_i \cap S_{i-1}$ for all $i = 1, \ldots, r+1$, and require that it is linearly equivalent to $H$ on $S_i$ and to $E$ on $S_{i-1}$. We require moreover that the limit as $t \in \mathbf{D}^*$ tends to 0 of the distinguished lines in the planes $\pi_r^{-1}(t)$, which we will denote by $R(t)$, be a smooth rational curve $\Lambda_0$ in $S_0$, linearly equivalent to $H$.

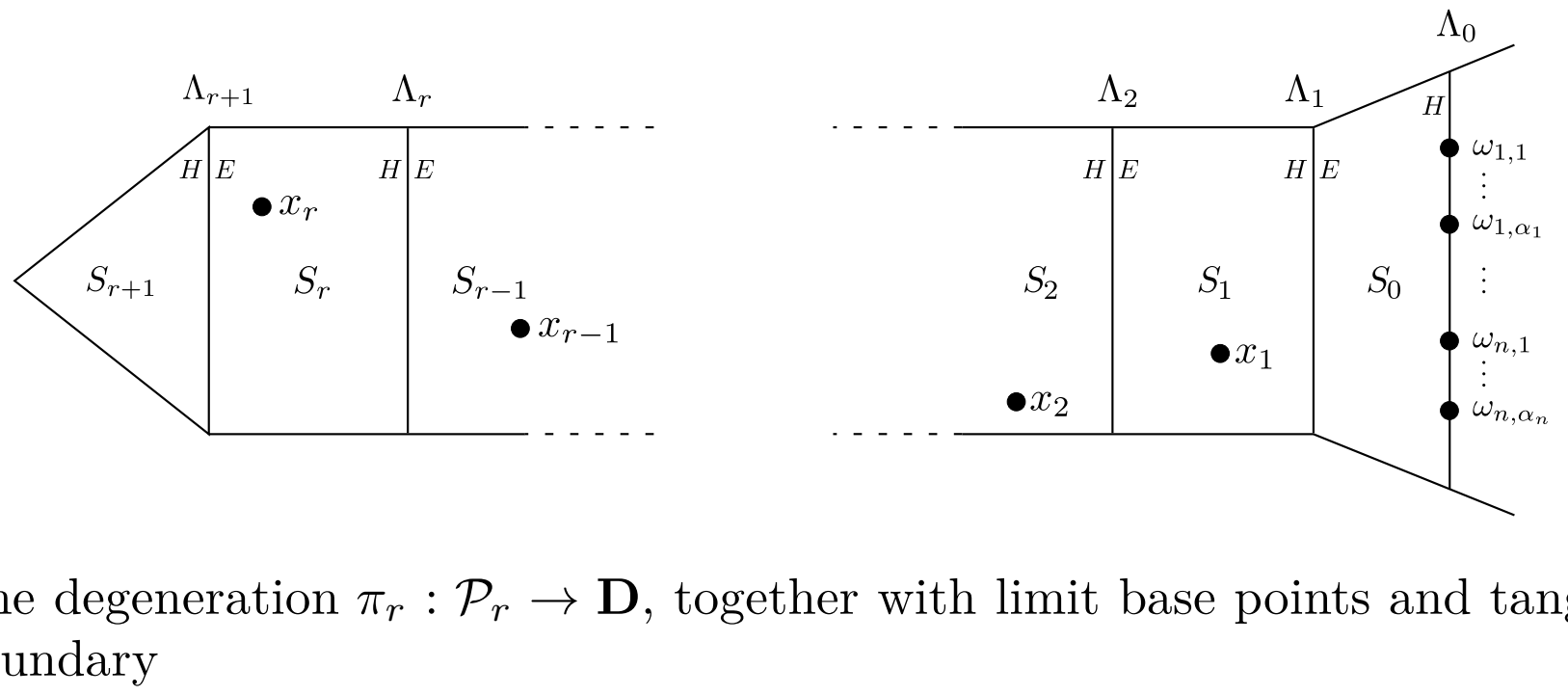

Figure 5: The degeneration $\pi_r : \mathcal{P}_r \to \mathbf{D}$, together with limit base points and tangency points along the boundary

The construction is by induction on $r \geqslant 0$. To construct $\mathcal{P}_0$, consider the trivial family $\mathbf{P}^2 \times \mathbf{D}$ together with a trivial family of lines $\Lambda \times \mathbf{D}$, and let $\varepsilon_0 : \mathcal{P}_0 \to \mathbf{P}^2 \times \mathbf{D}$ be the blow-up of $\Lambda \times \{0\}$, the distinguished line in the central fibre, and let $\pi_0 := \mathrm{pr}_2 \circ \varepsilon_0$; then $S_1 \cong \mathbf{P}^2$ is the proper transform of the initial central plane $\mathbf{P}^2 \times \{0\}$, $S_0 \cong \mathbf{F}_1$ is the exceptional divisor of $\varepsilon_0$, and $\Lambda_0$ is the intersection with $S_0$ of the proper transform of the initial family of lines $\Lambda \times \mathbf{D}$.

The inductive construction is pretty much identical. Assuming $\pi_r : \mathcal{P}_r \to \mathbf{D}$ is constructed, let $\varepsilon_{r+1} : \mathcal{P}_{r+1} \to \mathcal{P}_r$ be the blow-up of $\Lambda_0 \subseteq \mathcal{P}_r$, and $\pi_{r+1} := \pi_r \circ \varepsilon_{r+1}$. In $\mathcal{P}_{r+1}$ the surface $S_0$ is the exceptional divisor of $\varepsilon_{r+1}$, and the curve $\Lambda_0$ is the central fibre of the proper transform of the family of lines in $\mathcal{P}_r$ (or equivalently of the proper transform by $\varepsilon_0 \circ \cdots \varepsilon_{r+1}$ of the initial family of lines $\Lambda \times \mathbf{D}$). For $i = 1, \ldots, r+2$, the surface $S_i$ in $\mathcal{P}_{r+1}$ is the proper transform of $S_{i-1}$ in $\mathcal{P}_r$.

**(3.2) Degenerations of stable maps.** We now fix $d, g \in \mathbf{Z}_{\geqslant 0}$ and $\alpha, \beta \in \underline{\mathbf{N}}$ such that $I\alpha + I\beta = d$, and a set $\Omega$ of $\alpha$ points on $\Lambda$. We want to look at degenerations of the space $M_g^{d,\mathrm{irr}}(\alpha,\beta)(\Omega,x)$ for the pair $(\mathbf{P}^2, \Lambda)$ as the projective plane degenerates as in (3.1) above, for some specific degeneration of a set $x$ of $r = 2d + g - 1 + |\beta|$ points, as follows.

We consider a set $x = \{x_1, \ldots, x_r\}$ of sections of $\pi_r : \mathcal{P}_r \to \mathbf{D}$ such that for all $t \in \mathbf{D}^*$, the set $x(t)$ consists of $r$ general points in the plane $\pi_r^{-1}(t)$, and for $t = 0$, $x(0) = \{x_1(0), \ldots, x_r(0)\}$ with $x_i(0)$ a general point of $S_i \subseteq \pi_r^{-1}(0)$ for all $i \in [\![1, r]\!]$. We choose in addition a set $\Omega$ of $\alpha$ sections of $\pi_r : \mathcal{P}_r \to \mathbf{D}$, such that for all $t \in \mathbf{D}$, $\Omega(t)$ is a set of $\alpha$ general points on the line $R(t)$.

The points $x_1(0), \ldots, x_r(0)$ are indicated on Figure 5 (where they are simply indicated as $x_1, \ldots, x_r$ to lighten the picture), as well as the points of $\Omega(0)$, indicated on the picture as $\omega_{1,1}, \ldots, \omega_{1,\alpha_1}, \ldots, \omega_{n,1}, \ldots, \omega_{n,\alpha_n}$. Note that $S_{r+1}$, the proper transform of the central fibre of $\mathbf{P}^2 \times \mathbf{D}$, does not support anything, $S_i$ supports the base point $x_i(0)$ for all $i = 1, \ldots, r$, and $S_0$ supports all the fixed contact points on the boundary $\Lambda_0$.

We then let $M_g^{d,\lim}(\alpha,\beta)(\Omega,x)$ be the limit as $t \in \mathbf{D}^*$ tends to 0 of the moduli spaces of stable maps $M_g^{d,\mathrm{irr}}(\alpha,\beta)(\Omega(t),x(t))$ for the pairs $(\pi_r^{-1}(t), R(t))$. It is a moduli space of stable maps $f : C \to \pi_r^{-1}(0)$ satisfying the following conditions:

(i) $C$ is a connected nodal curve of arithmetic genus $g$;

(ii) $f(C) \supset x(0)$;

(iii) $f(C)$ does not contain $\Lambda_0$, and there exist $\alpha$ points $q_{i,j} \in C$, and $\beta$ points $r_{i,j} \in C$ such that
$$\big\{\phi(q_{i,j})\big\} = \Omega(0) \qquad \text{and} \qquad \phi^*\Lambda_0 = \sum_{1 \leqslant j \leqslant \alpha_i} i\, q_{i,j} + \sum_{1 \leqslant j \leqslant \beta_i} i\, r_{i,j};$$

(iv) if $p \in f^{-1}(S_{i+1} \cap S_i)$ does not lie on a component of $C$ which is entirely mapped to $S_{i+1} \cap S_i$, then it is a node of $C$, the intersection of two components $C'$ and $C''$ mapped to $S_{i+1}$ and $S_i$ respectively, and appears with the same multiplicity $\mu_p$ in $(f|_{C'})^*\Lambda_{i+1}$ and $(f|_{C''})^*\Lambda_{i+1}$.

That such a limit of moduli spaces of stable maps indeed exists is a highly non-trivial fact which is due to Jun Li, see [8, 9].

The following result is the core of the proof by degeneration of the Correspondence Theorem (1.9).

**(3.3) Proposition** (Brugallé). *The moduli space $M_g^{d,\lim}(\alpha,\beta)(\Omega,x)$ (defined in (3.2) above) is finite, and depends on $(\Omega,x)$ only by way of $(\Omega(0),x(0))$.*

*Let $\big[f : C \to \pi_r^{-1}(0)\big] \in M_g^{d,\lim}(\alpha,\beta)(\Omega,x)$. Then all the irreducible components $C_1, \ldots, C_n$ of $C$ are smooth rational curves, none of them is entirely mapped to one of the curves $\Lambda_i$, $0 \leqslant i \leqslant r+1$, and the following conditions hold after suitable relabelling:*

*(i) $C_1, \ldots, C_d$ are pairwise disjoint smooth rational curves, and there exist $1 \leqslant a_1 < \cdots < a_d \leqslant r$ such that for each $i \in [\![1, d]\!]$ the curve $f(C_i)$ is a member of the linear system*

*$|H + k_i F|$ on $S_{a_i}$ (for some nonnegative integer $k_i$) and passes through the point $x_{a_i}(0) \in S_{a_i}$;*

*(ii) the remaining curves form $r - d$ pairwise disjoint chains of smooth rational curves*

$$\bar{C}_j = C_{d+\ell_1+\cdots+\ell_{j-1}+1} \cup \cdots \cup C_{d+\ell_1+\cdots+\ell_{j-1}+\ell_j}$$

*(i.e., for all $1 \leqslant s < s' \leqslant \ell_j$, $C_{d+\ell_1+\cdots+\ell_{j-1}+s}$ and $C_{d+\ell_1+\cdots+\ell_{j-1}+s'}$ intersect if and only if $s' = s + 1$, and in that case transversally in exactly one point, and $\bar{C}_j \cap \bar{C}_{j'}$ is empty if $j \neq j'$); there are integers $b_1, \ldots, b_{r-d} \in [\![0, r]\!]$ and $\mu_1, \ldots, \mu_{r-d} \geqslant 1$ such that the divisor $f_*(C_{d+\ell_1+\cdots+\ell_{j-1}+s})$, $1 \leqslant j \leqslant r - d$ and $1 \leqslant s \leqslant \ell_j$, belongs to the linear system $|\mu_j F|$ on the ruled surface $S_{b_j+s-1}$; for all $j = 1, \ldots, r - d$, the chain $\bar{C}_j$ intersects exactly 1 or 2 curves among $C_1, \ldots, C_d$, necessarily at points mapped to $\Lambda_r \cup \cdots \cup \Lambda_1$ by $f$, and $f(\bar{C}_j)$ passes through exactly one of the points in $(\Omega(0), x(0))$.*

**(3.4) Example.** Depicted in Fig. 6a is an example of a source curve $C$ of a member $\big[f : C \to \pi_7^{-1}(0)\big]$ of $M_0^{3,\mathrm{lim}}(0, [1, 1])$, together with its image by $f$ in Fig. 6b. The circled numbers on Fig. 6b indicate that the restrictions of $f$ to $C_4, C_5, C_6$ have degree 2; the restrictions to all other irreducible components of $C$ are birational.

The picture on $S_3$ means that $f(C_1)$ is the only curve in $S_3$: it is irreducible, linearly equivalent to $H + 2F$.[4] On the other hand, there are two irreducible components of $f(C)$ in $S_6$, linearly equivalent to $H$ and $F$, respectively. Note that $f$ normalizes the one node indicated by a square dot on $S_6$.

To help the reader get their way around, I indicated on Fig. 6a one point in the pre-image of each element of $x(0)$; there is only one point in $f^{-1}(x_i(0))$ if $i \neq 1$, but there are two points in $f^{-1}(x_1(0))$, both in $C_5$. In this example one has

$$(a_1, a_2, a_3) = (3, 6, 7) \qquad \bar{C}_1 = C_4 \cup C_5 \cup C_6 \qquad \bar{C}_3 = C_{10} \cup C_{11}$$
$$(b_1, b_2, b_3, b_4) = (0, 0, 4, 4), \qquad \bar{C}_2 = C_7 \cup C_8 \cup C_9 \qquad \bar{C}_4 = C_{12} \cup C_{13} \cup C_{14}.$$

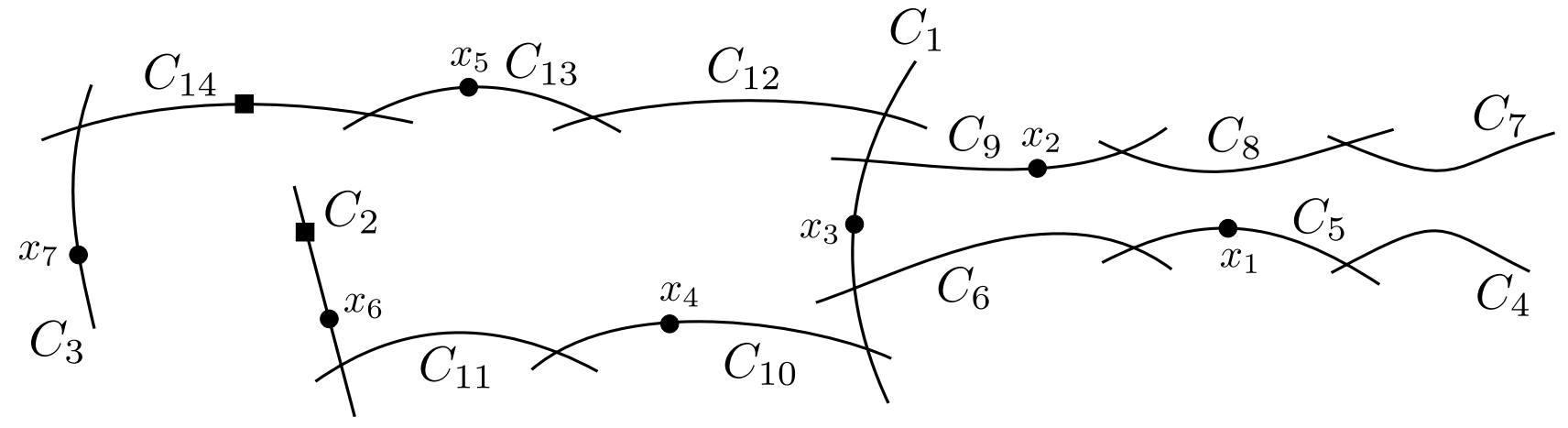

Figure 6a: Source curve of a member of $M_0^{3,\mathrm{irr}}(0, [1, 1])(x(0))$

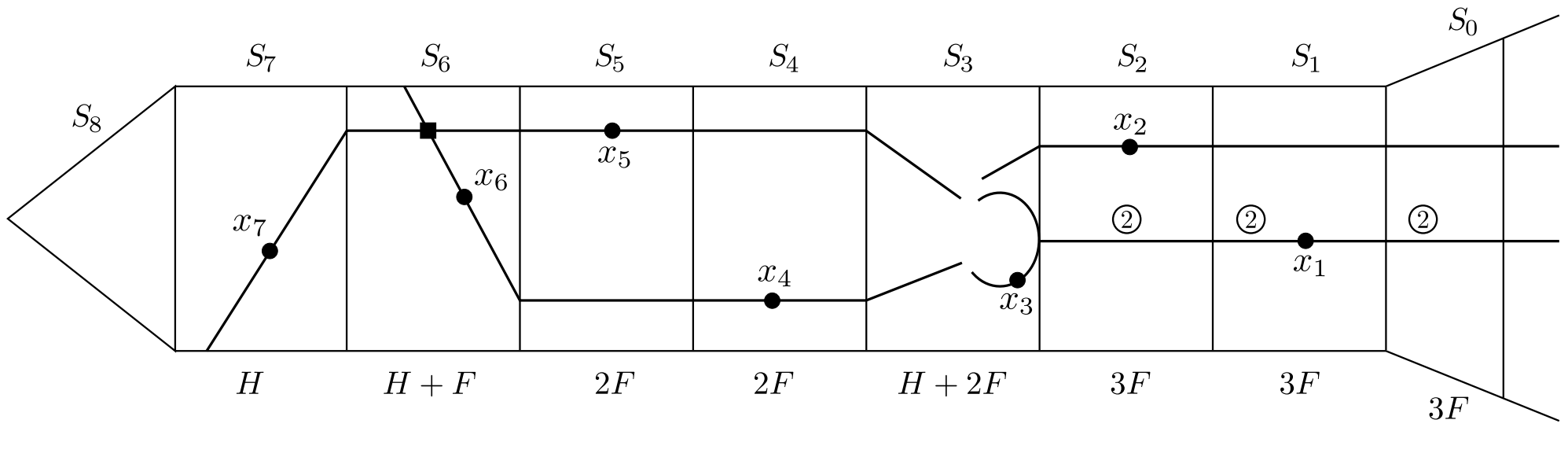

Figure 6b: Image in $\pi_7^{-1}(0)$

[4] the picture is slightly incorrect, as $f(C_1)$ should be unisecant to all fibres of the ruling of $S_3$; we will only be able to fix this when we use tropical representations of curves, see [X, Section 8]. There is the same issue with the two curves linearly equivalent to $H$ on $S_6$ and $S_7$, respectively.

**(3.5) Proof of Proposition (3.3).** Let $\big[f : C \to \pi_r^{-1}(0)\big] \in M_g^{d,\lim}(\alpha,\beta)(\Omega,x)$. The fact that no component of $C$ is entirely mapped to one of the curves $\Lambda_i$ follows from [6, Example 11.4 and Lemma 14.6], and we will take it for granted. The other assertions essentially follow from a series of dimension counts, which are now going to expound.

Let $C_1,\dots,C_n$ be the irreducible components of $C$. For all $j=1,\dots,n$, we let $f_j : C_j \to S_{s(j)}$ be the restriction of $f$ to $C_j$, $g_j$ the arithmetic genus of $C_j$, $\xi_j$ the linear equivalence class of $f_*(C_j)$ on the surface $S_{s(j)}$ (note that a priori $\xi_j$ may vanish, in case $C_j$ is contracted), $\zeta_j \in \{0,1\}$ the number of points in $x(0)$ contained in $f_j(C_j)$, and $\beta_j^E, \beta_j^H \in \underline{\mathbf{N}}$ the types of contact of $f_j(C_j)$ with $\Lambda_{s(j)+1}$ and $\Lambda_{s(j)}$ respectively, pretending that all contacts are at unassigned points. Thus

$$\big[f_j\big] \in M_{g_j}^{\xi_j,\mathrm{irr}}(0,0,\beta_j^E,\beta_j^H)$$

as a moduli space of stable maps to $S_{s(j)}$ with boundary $R = \Lambda_{s(j)} \coprod \Lambda_{s(j)+1}$, where the set $\Omega$ is empty, and where we choose to forget the datum of the possible point of $x(0)$ contained in $f(C_j)$. (This all comes with a harmless abuse of notation in case $s(j)=r+1$; in fact we shall see that this case never occurs).

Then we consider $M_j \subseteq M_{g_j}^{\xi_j,\mathrm{irr}}(0,0,\beta_j^E,\beta_j^H)$ the smallest irreducible subvariety in which $[f_j]$ moves as we let $x$ vary, maintaining the requirements of (3.2). We let $CM_j$ be the corresponding family of curves in $S_{s(j)}$. By generality of the point $x_{s(j)}(0)$ in $S_{s(j)}$, we must have $\dim(CM_j) \geqslant \zeta_j$; let $\varepsilon_j = \dim(CM_j) - \zeta_j$. By Theorem (2.4) and Paragraph (2.5) we also have:

$$\zeta_j \leqslant \dim(CM_j) = \zeta_j + \varepsilon_j \leqslant 2F\cdot\xi_j + g_j - 1 + |\beta_j^E| + |\beta_j^F|. \qquad (3.5.1)$$

On the other hand, the various $f(C_j)$'s are subject to matching conditions at the boundary of the various $S_{s(j)}$, in order for condition (iv) of (3.2) to hold. (A priori there may be other matching conditions if the $C_j$'s intersect at points mapped off $\Lambda_{r+1}\cup\dots\cup\Lambda_1$ by $f$, but we don't need to consider them for the moment; we shall see that in fact this does not happen). These matching conditions are in the number of

$$a = \mathrm{Card}\big(f^{-1}(\Lambda_{r+1}\cup\dots\cup\Lambda_1)\big),$$

and they are independent by Theorem (2.4), (2.5), and [III]: consider some curve $f(C_{j_0})$, and assume all the other components $f(C_j)$ are fixed; then the matching conditions imposed to $f(C_{j_0})$ amount to turn some of the contact conditions $\beta_j^E, \beta_j^H$ at unassigned points to contact conditions at assigned points, and it follows from [III, $**$] that each assigned contact point makes the dimension drop by 1. Note that the assigned contact points may not be in general position; fortunately this does not cause any trouble, as we have observed in [III].

Lastly, the assigned contact points $\Omega(0)$ on $\Lambda_0$ impose $|\alpha|$ further conditions, independent from all the rest for the same reasons as before. Therefore, if we let $CM$ be the smallest irreducible variety in which $f(C)$ moves as we let $x$ vary, then $CM$ has codimension at least $a+|\alpha|$ in $\prod_j CM_j$.

In the upshot, we must have:

$$\textstyle\sum_j \dim(CM_j) \geqslant \big(\sum_j \zeta_j\big) + a + |\alpha| = r + a + |\alpha|. \qquad (3.5.2)$$

On the other hand, summing the inequalities (3.5.1) for all $j=1,\dots,n$, we get

$$\sum_j \dim(CM_j) \leqslant \sum_j \Big(2F\cdot\xi_j + g_j - 1 + |\beta_j^E| + |\beta_j^F|\Big) \stackrel{\star}{=} 2d + g - 1 + |\beta| + |\alpha| + a - a' = r + |\alpha| + a - a', \qquad (3.5.3)$$

where $a'$ is the number of intersection points of the $C_j$'s mapped off $\Lambda_{r+1}\cup\cdots\cup\Lambda_1$ by $f$; equality $(\star)$ is given by the three relations of Lemma (3.6) below. Comparing with (3.5.2), we get that

$$a' = 0 \quad \text{and} \quad \textstyle\sum_j \dim(CM_j) = r + |\alpha| + a,$$

which in particular implies that the family of curves corresponding to $M_g^{d,\lim}(\alpha,\beta)(\Omega,x)$ is finite and depends only on $(\Omega(0),x(0))$.

The inequalities (3.5.2) and (3.5.3) also imply that for all $j = 1,\dots,n$ one has

$$\dim(CM_j) = \zeta_j + \varepsilon_j = 2F\cdot\xi_j + g_j - 1 + |\beta_j^E| + |\beta_j^F|, \tag{3.5.4}$$

and moreover $\sum_j \varepsilon_j = |\alpha| + a$. By Theorem (2.4) and (2.5), Equation (3.5.4) implies that each $M_j$ is an enumeratively relevant irreducible component of $M_{g_j}^{\xi_j,\mathrm{irr}}(0,0,\beta_j^E,\beta_j^H)$. In particular, the finiteness of the family of curves in $\bigcup_j S_j$ corresponding to $M_g^{d,\lim}(\alpha,\beta)(\Omega,x)$ implies the finiteness of $M_g^{d,\lim}(\alpha,\beta)(\Omega,x)$ itself.

These finiteness results imply that for all $j$, if we fix the intersection points of $f(C_j)$ with $\Lambda_{s(j)+1}$ and $\Lambda_{s(j)}$ (again with a slight abuse of notation if $s(j) = r+1$), then $[f_j]$ moves in dimension at most $\zeta_j$. By (2.8), this implies that for all $j$ one has $g_j = 0$, $s(j) < r+1$, and

$$\xi_j = \begin{cases} k_jF & \text{if } \zeta_j = 0 \\ k_jF \text{ or } H + k_jF & \text{if } \zeta_j = 1. \end{cases}$$

By relation (3.6.3) below, there are exactly $d$ indices $j$ such that $\xi_j$ is of the form $H + k_jF$ and, relabelling if necessary, we will assume that they are $j = 1,\dots,d$. Then, the indices $s(1),\dots,s(d)$ are pairwise distinct and in the interval $[\![1,r]\!]$, because there is only one point of $x(0)$ on each of $S_1,\dots,S_r$.

Condition (iv) of (3.2) imposes that the curves $C_{d+1},\dots,C_n$ (i.e., those $C_j$ such that $\xi_j = k_jF$) organize in chains of fibres of the same multiplicity as stated in Proposition (3.3). The vanishing of $a'$ imposes that each chain meets exactly one or two among $C_1,\dots,C_d$, at points mapped to $\Lambda_{r+1}\cup\cdots\cup\Lambda_1$. The finiteness results above imply that each chain must pass through exactly one point of $(\Omega(0),x(0))$.

It remains to show that two curves among $C_1,\dots,C_d$ cannot intersect each other. Assume by contradiction that $f(C_1)\subseteq S_1$, $f(C_2)\subseteq S_2$, and $C_1$ and $C_2$ intersect at a point mapped to $y \in S_1\cap S_2 = \Lambda_2$. Then we would find a 1-dimensional family of stable maps in $M_g^{d,\lim}(\alpha,\beta)(\Omega,x)$ by letting $y$ move along $\Lambda_2$, a contradiction. This ends the proof of Proposition (3.3). ☐

**(3.6) Lemma.** *In the above notation, we have the following relations:*

$$\textstyle\sum_j \big(|\beta_j^E| + |\beta_j^H|\big) = 2a + |\alpha| + |\beta|; \tag{3.6.1}$$

$$g = \big(\textstyle\sum_j g_j\big) + 1 - n + a + a'. \tag{3.6.2}$$

$$\textstyle\sum_j (\xi_j\cdot F)_{S_{s(j)}} = d; \tag{3.6.3}$$

*Proof.* The points corresponding to those $\beta_j^E$ such that $0 \leqslant s(j) \leqslant r$, and equally those corresponding to those $\beta_j^H$ such that $1 \leqslant s(j) \leqslant r+1$, are exactly the points in the singular locus of $\pi_r^{-1}(0)$ at which the matching conditions referred to in the proof of Proposition (3.3) above happen, which are in the number $a$. Each such point is counted twice, once in some number $\beta_{s(j)}^E$, and once in the corresponding number $\beta_{s(j)+1}^H$. On the other hand, those $\beta_j^H$ such that $s(j) = 0$ sum up as the number of points in $\Omega(0)$, which is $|\alpha| + |\beta|$. This proves (3.6.1).

Relation (3.6.2) is the expression of the arithmetic genus of the nodal curve $C$, which has $n$ irreducible components $C_1\dots,C_n$ of respective genera $g_1,\dots,g_n$, and $a + a'$ nodes, where $a'$ is the number of intersection points among the $C_j$'s that are mapped off $\Lambda_{r+1}\cup\cdots\cup\Lambda_1$ by $f$.

Relation (3.6.3) comes with an abuse of notation, namely that for those $j$ such that $s(j) = r+1$, $F$ should be replaced by $H$. Then the curve class on the limit surface $\pi_r^{-1}(0)$ consisting of the class $F$ on all $S_i$, $i = 0, \dots, r+1$ (which thus has to be taken as $F$ on all $S_i$, $i = 0, \dots, r$, and $H$ on $S_{r+1}$), is a limit of the class $H$ on the projective planes $\pi_r^{-1}(t)$ as $t \in \mathbf{D}^*$ tends to $0 \in \mathbf{D}$. Relation (3.6.3) is then a consequence of the fact that $\big[f : C \to \pi_r^{-1}(0)\big] \in M_g^{d,\lim}(\alpha,\beta)(\Omega,x)$ is a limit of stable maps realizing the curve class $dH$ on the projective planes $\pi_r^{-1}(t)$, $t \in \mathbf{D}^*$. See also (5.1) for another point of view on this argument. □

Proposition (3.3) has the following corollary, which we will translate in Section 4 below as saying that the number of plane curves $N_g^{d,\mathrm{irr}}(\alpha,\beta)$ may be computed by enumerating floor diagrams.

**(3.7) Corollary.** *In the notation of (3.1) to (3.3), one has*

$$N_g^{d,\mathrm{irr}}(\alpha,\beta) = \sum_{[f:C\to\pi_r^{-1}(0)]\in M_g^{d,\lim}(\alpha,\beta)(\Omega,x)} \mu(f)$$

*where $\mu(f)$ is defined as*

$$\mu(f) = \Big(\prod_{p\in f^{-1}(\Lambda_{r+1}\cup\cdots\cup\Lambda_1)} \mu_p\Big)\Big(\prod_{p\in x(0)} \mathrm{Card}\, f^{-1}(p)\Big)\Big(\prod_{1\leqslant j\leqslant n} \mathrm{Card}(\mathrm{Aut}(f|_{C_j}))\Big)^{-1}$$

*(recall in particular the definition of $\mu_p$ in Condition* (iv) *of (3.2)).*

*Proof.* By Proposition (3.3) and Jun Li's Formula, see [8, 9], for all $t \in \mathbf{D}^*$, one has

$$\sum_{f\in M_g^{\xi,\mathrm{irr}}(\alpha,\beta)(\Omega(t),x(t))} \mu(f,x(t)) = \sum_{[f:C\to\pi_r^{-1}(0)]\in M_g^{d,\lim}(\alpha,\beta)(\Omega,x)} \mu(f).$$

The result follows, since the left-hand-side of the equation above equals $N_g^{d,\mathrm{irr}}(\alpha,\beta)$, as we have seen in (2.7). □

**(3.8)** It will be useful to have a ready-to-use expression for the multiplicity $\mu(f)$ of a stable map $[f] \in M_g^{d,\lim}(\alpha,\beta)(\Omega,x)$ as in Proposition (3.3). Only the chains of multiple fibers $\bar{C}_1,\dots,\bar{C}_{r-d}$ may contribute non-trivially to the multiplicity $\mu(f)$, and their respective contributions may be computed in a similar way as in (2.9).

Let $\bar{C}_j$ be such a chain, with multiplicity $\mu_j$ and length $\ell_j$ (the number of its irreducible components). There are the following three cases:

(a) no irreducible components of $\bar{C}_j$ is mapped to $S_0$: then $f(\bar{C}_j)$ passes through exactly one point $p$ of $x(0)$, for which $f^{-1}(p)$ has cardinality $\mu_j$, there are $\ell_j+1$ points $p \in f^{-1}(\Lambda_{r+1}\cup\cdots\cup\Lambda_1)$ lying on $\bar{C}_j$, all with multiplicity $\mu_p = \mu_j$, and for each irreducible component of $\bar{C}_j$, there is an automorphism group of cardinality $\mu_j$. Thus the contribution of the chain $\bar{C}_j$ to the multiplicity $\mu(f)$ is
$$\mu_j^{\ell_j+1}\cdot\mu_j\cdot\mu_j^{-\ell_j} = \mu_j^2;$$

(b) there is an irreducible components of $\bar{C}_j$ mapped to $S_0$, and $f(\bar{C}_j)$ passes through exactly one point $p$ of $x(0)$ and no point of $\Omega(0)$: then there are only $\ell_j$ points $p \in f^{-1}(\Lambda_{r+1}\cup\cdots\cup\Lambda_1)$ lying on $\bar{C}_j$, but otherwise the situation is the same as in case (a). The contribution of the chain $\bar{C}_j$ to the multiplicity $\mu(f)$ is thus
$$\mu_j^{\ell_j}\cdot\mu_j\cdot\mu_j^{-\ell_j} = \mu_j;$$

(c) there is an irreducible components of $\bar{C}_j$ mapped to $S_0$, and $f(\bar{C}_j)$ passes through exactly one point $p$ of $\Omega(0)$ and no point of $x(0)$: then there are $\ell_j$ points $p \in f^{-1}(\Lambda_{r+1} \cup \cdots \cup \Lambda_1)$ lying on $\bar{C}_j$ as in case (b), but now there is no point $p \in x(0)$ contributing to the multiplicity. The contribution of the chain $\bar{C}_j$ to the multiplicity $\mu(f)$ is thus

$$\mu_j^{\ell_j} \cdot 1 \cdot \mu_j^{-\ell_j} = 1.$$

# 4 – Combinatorial translation

We are now ready to give the proof of [2] of the Correspondence Theorem (1.9). In fact, most of the work has already been done in Section 3, and it only remains to interpret Proposition (3.3) and Corollary (3.7) in terms of floor diagrams.

**(4.1) Proof of the Correspondence Theorem (1.9).** We associate each stable map $[f : C \to \pi_r^{-1}(0)] \in M_g^{d,\lim}(\alpha,\beta)(\Omega,x)$ as in Proposition (3.3) with a marked floor diagram of genus $g$, degree $d$, and type $(\alpha,\beta)$ as follows.

We first identify the vertices. To each irreducible component $C_i$ of $C$ such that $f_*C_i = H+kF$ on some $S_{s(i)}$ we associate a floor, and to each point in $f(C) \cap \Lambda_0$ we associate a vertex at infinity. We provide the vertices with an ordering, which in turn will determine the orientation of the edges: each vertex $v$ corresponds to an object lying on a unique ruled surface $S_{s(v)}$, and for all vertices $v$ and $v'$, we declare that $v$ precedes $v'$ if $s(v) \leqslant s(v')$.

Now let us define the edges: they are in bijective correspondence with the chains $\bar{C}_j$ described in (ii) of Proposition (3.3), and the weight of the edge corresponding to $\bar{C}_j$ is the positive integer $\mu_j$. The edge corresponding to $\bar{C}_j$ is adjacent to the floor corresponding to $C_i$ (resp. the vertex at infinity corresponding to $p \in \Lambda_0$) if $\bar{C}_j$ and $C_i$ intersect at a point (resp. if $f(\bar{C}_j)$ passes through $p$).

It follows from Proposition (3.3) that the object constructed so far is a floor diagram $\mathcal{D}$ of genus $g$ and degree $d$, and the corresponding set $\mathsf{D} = \mathrm{Vert}^{\mathrm{fin}} \coprod \mathrm{Edges}$ is in bijection with $\Omega(0) \coprod x(0)$. The set $x(0)$ naturally comes with a bijection with $[\![1,r]\!]$ (remember that each $x_i(0)$ lies on the surface $S_i$). By appropriately labelling the points in $\Omega(0)$ with indices in $[\![-|\alpha|+1,0]\!]$, we obtain a marking of $\mathcal{D}$ of type $(\alpha,\beta)$. There are several possible ways of labelling the points in $\Omega(0)$, but they are equivalent up to automorphisms of the marked floor diagram.

Conversely, from any marked floor diagram of type $(\alpha,\beta)$ one can construct a stable map as in Proposition (3.3), and this is a point in $[f : C \to \pi_r^{-1}(0)] \in M_g^{d,\lim}(\alpha,\beta)(\Omega,x)$. It only remains to observe that the multiplicity of the map $[f : C \to \pi_r^{-1}(0)] \in M_g^{d,\lim}(\alpha,\beta)(\Omega,x)$ in the formula of Corollary (3.7) equals that of the corresponding floor diagram, which follows from (3.8) and Definition (1.8).

**(4.2) Example.** I shall now illustrate the correspondence of (4.1) between floor diagrams and limit curves on various examples pertaining to plane cubic curves. I will only consider curves with trivial boundary conditions, i.e., $(\alpha,\beta) = (0,3)$, so nothing relevant happens on $S_0$ (and on $S_{r+1}$ either, as always), and thus I will only represent the $S_i$'s for $i = 1,\ldots,r$.

Under each ruled surface $S_i$, I have indicated the total curve class supported on this component, i.e., the sum of the $\xi_j = f_*C_j$ for all $j$ such that $C_j$ is mapped to $S_i$. The limit base points $x_1(0),\ldots,x_r(0)$ on $S_1,\ldots,S_r$ respectively are indicated by black disks.

Note that (obviously, I hope) to obtain a proper graphic correspondence between our representations of floor diagrams and limit curves, the diagrams should be rotated by $+\frac{\pi}{2}$.

*(4.2.1) Limit of a smooth cubic.* Below is displayed the limit curve of a smooth cubic, corresponding to the unique possible marked floor diagram of degree 3 and genus 1 (and type $(0, 3)$). The two curves in $S_7$ and $S_4$ are irreducible. Later on (in Lecture [X]) we will have a better way of representing these curves, coming from tropical geometry.

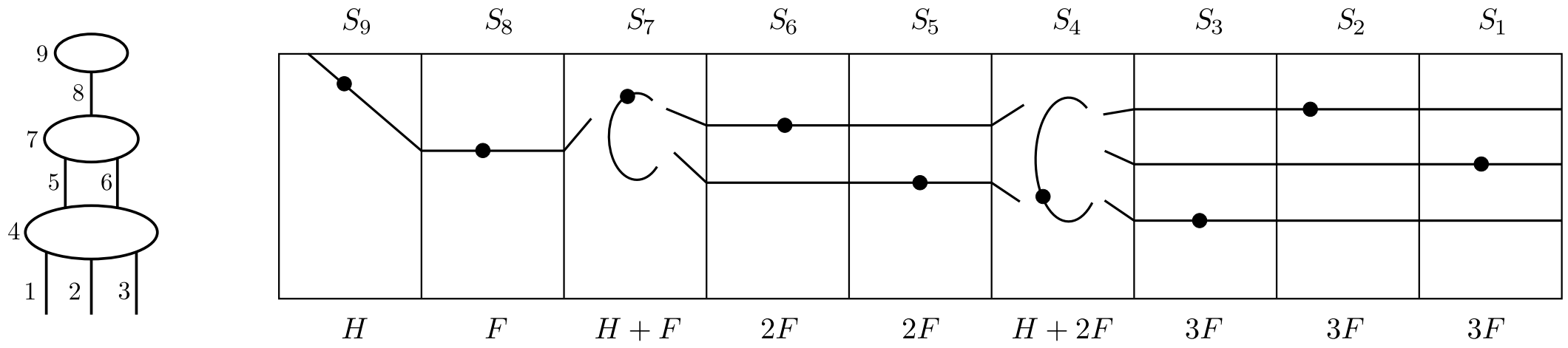

Figure 7: Limit of a smooth cubic

*(4.2.2) Various limits of rational cubics.* Below I display limit curves corresponding to the various diagrams of degree 3 and genus 0 (and type $(0, 3)$), see Example (1.7.1). For the diagram shape at the center of Figure 2, I consider two different markings (Figures 9 and 10 below), so that one can better see what is encoded in the marking.

On Figure 8, the two curves in $S_6$ and $S_4$ are irreducible, and tangent to $S_6 \cap S_5$ and $S_5 \cap S_4$ respectively; the curve on $S_5$ is a fibre with multiplicity 2.

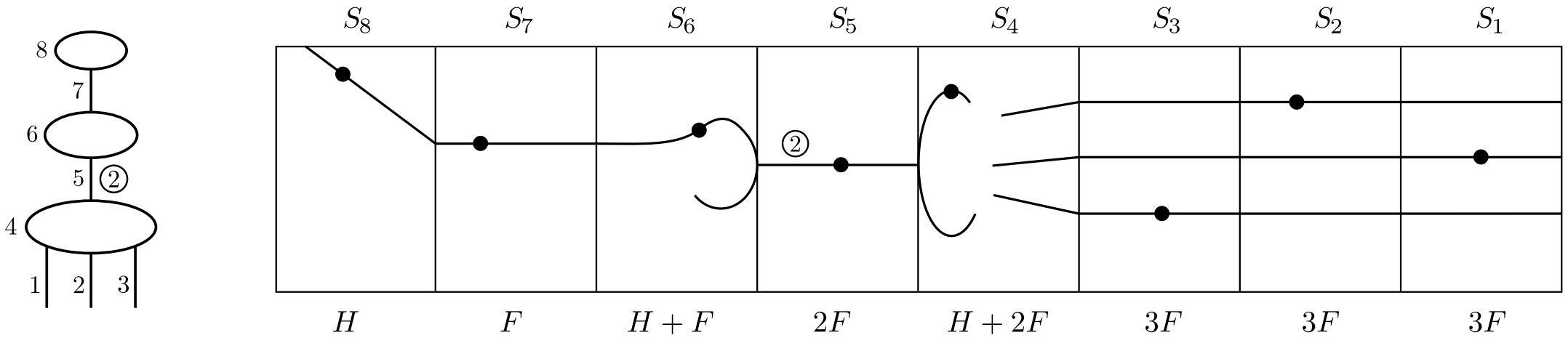

Figure 8: Limit of rational cubics (a)

On Figure 9, the curve in $S_6$ is irreducible, and that in $S_4$ has two components of classes $H + F$ and $F$ respectively, which intersect at one point represented by a full black square. On Figure 10 the situation is similar.

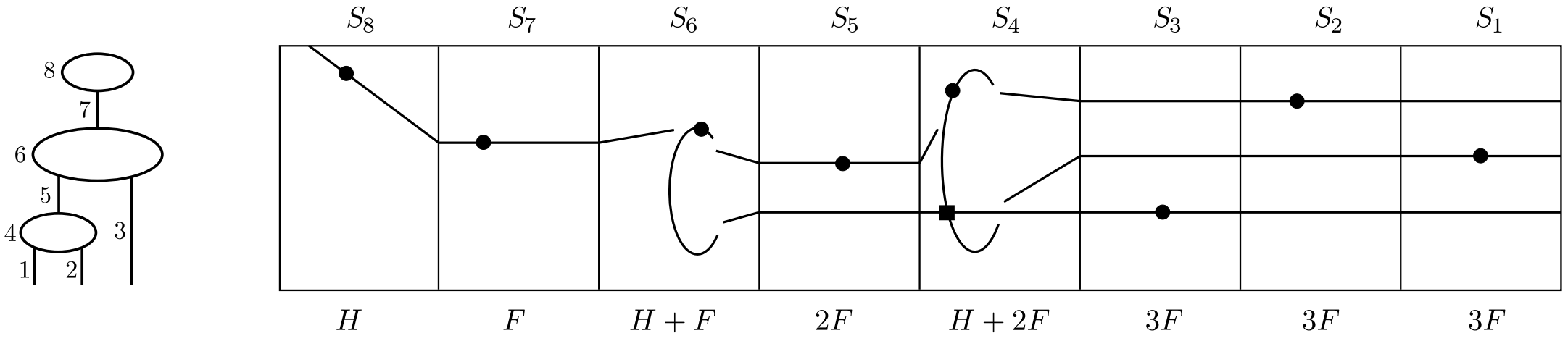

Figure 9: Limit of rational cubics (b)

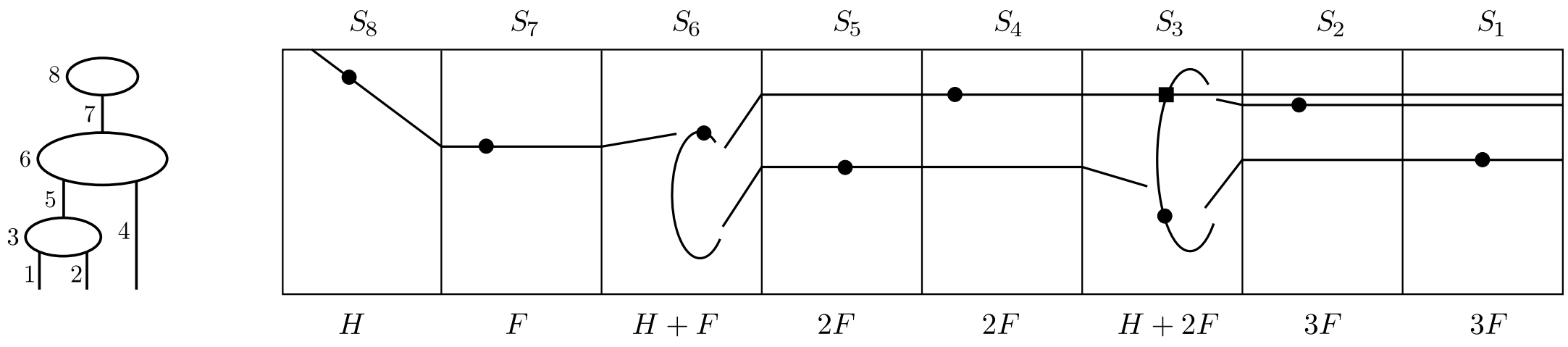

Figure 10: Limit of rational cubics (b')

On Figure 11, the curve in $S_4$ is irreducible, and that in $S_7$ has two components of classes $H$ and $F$ respectively, which intersect at one point represented by a black square.

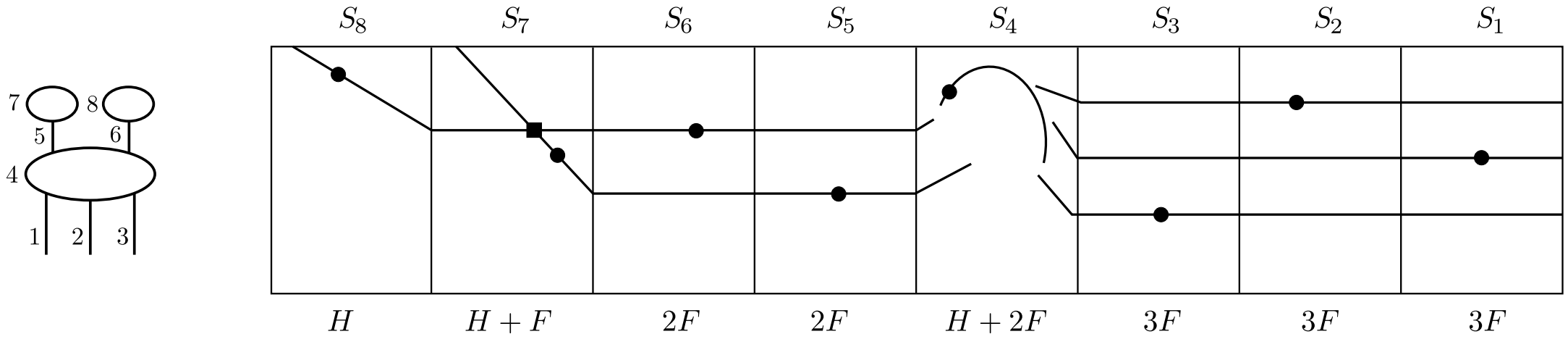

Figure 11: Limit of rational cubics (c)

I have tried to maintain the black dots representing the points of $x(0)$ in the series of Figures 8–11 at a fixed position. I have only cheated a little bit on Figure 11, on which I have moved the base point on $S_7$, for otherwise the picture was outrageously unclear.

# 5 – Projective geometric interpretation

In this section, I propose a formulation of the degenerations of stable maps considered in this lecture in terms of degenerations of nodal curves, in the framework of the theory expound in [I].

**(5.1) Twists.** Let us consider the degeneration $\pi_r : \mathcal{P}_r \to \mathbf{D}$ of the projective plane introduced in (3.1). It comes with a map $\bar{\varepsilon}_r : \mathcal{P}_r \to \mathcal{P} := \mathbf{P}^2 \times \mathbf{D}$ which is the composition of blow-ups $\bar{\varepsilon}_r = \varepsilon_r \cdots \varepsilon_0$. Consider the line bundle

$$\mathcal{L}_{r,d} = \bar{\varepsilon}_r^* \operatorname{pr}_1^*\big(\mathcal{O}_{\mathbf{P}^2}(d)\big);$$

it restricts to $\mathcal{O}_{\mathbf{P}^2}(d)$ on the fibres $\pi_r^{-1}(t)$, $t \in \mathbf{D}^*$, and thus its restriction to the central fibre $\pi_r^{-1}(0)$ defines a possible linear system for the linear system of plane degree $d$ curves to degenerate to within the degeneration $\pi_r : \mathcal{P}_r \to \mathbf{D}$. In this paragraph I will describe the various *twists* $\mathcal{L}_{r,d}(-W)$ that have the same property: the requirements are that $W$ is an effective divisor of $\mathcal{P}_r$ supported on the central fibre, such that the restriction $\mathcal{L}_{r,d}(-W)|_{\pi_r^{-1}(0)}$ defines a linear system on the central fibre that has the same dimension as $|\mathcal{O}_{\mathbf{P}^2}(d)|$.

It is elementary to prove that these twists all have the form

$$\mathcal{L}_{r,d}\big(-\textstyle\sum_{i=0}^{r} m_i \bar{S}_i\big),$$

where $\bar{S}_i = \sum_{i' \leqslant i} S_{i'}$, and the $m_i$'s are non-negative integers subject to the inequalities

$$d \geqslant m_r \geqslant \cdots \geqslant m_0.$$

As a sideremark, note that this gives an alternative way of proving the relation (3.6.3). Since $\mathcal{L}_{r,d}$ restricts to $dH$ on $S_{r+1}$ and to $dF$ on the other $S_i$'s, the line bundle $\mathcal{L}_{r,d}\big(-\sum_{i=0}^{r} m_i \bar{S}_i\big)$ restricts to the following linear equivalence classes on the various irreducible components of the central fibre of $\pi_r$:

— on $S_{r+1}$: $(d-m_r)H$;
— on $S_i$, $0 < i \leqslant r$: $dF + m_i E - m_{i-1}H = (m_i - m_{i-1})H + (d-m_i)F$;
— on $S_0$: $dF + m_0 = m_0 H + (d-m_0)F$.

Now, the choice of the limiting positions of the base points that has been made in (3.2), i.e., the choice that $x_i(0)$ lies on $S_i$ for all $i = 1,\dots,r$, implies that the twists that are relevant for the enumeration of limit curves are subject to the further restrictions:

$$m_r = d; \quad m_0 = 0; \quad m_i - m_{i-1} \leqslant 1 \quad \text{for all } i = r,\dots,1.$$

This corresponds to the fact that the classes $\xi_j$, in the notation of Proposition (3.3), are necessarily of the form $H + k_j F$ or $k_j F$.

*(5.1.1) Example.* The twist giving the limiting linear system in which the curve represented on Figure 6b lives is the following:

$$\mathcal{L}_{7,3}\big(-3\bar{S}_7 - 2\bar{S}_6 - \bar{S}_5 - \bar{S}_4 - \bar{S}_5\big).$$

**(5.2) Projective models of the degeneration.** I now describe the projective model of the degeneration $\pi_r : \mathcal{P}_r \to \mathbf{D}$ given by a twist which is relevant for the enumeration of limit curves with the limit base points as in (3.2).

By what we have seen above, such a twist gives the classes $H, H+F, \dots, H+(d-1)F$ respectively on the ruled surfaces $S_{r_0}, S_{r_1}, \dots, S_{r_{d-1}}$, with $r = r_0 > r_1 > \cdots > r_{d-1}$, and multiples of the class $F$ on the other ruled surfaces (specifically, it gives the class $kF$ on the surfaces $S_i$ for $r_{k-1} > i > r_k$, for all $k = 1,\dots,d$, setting $r_{-1} = r+2$ and $r_d = -1$ for convenience, and with the abuse of considering $S_{r+1}$ as a ruled surface).

It follows that the restriction of such a twist to the central fibre maps $S_{r_0}, S_{r_1}, \dots, S_{r_{d-1}}$ to rational normal scrolls of types $\Sigma(0,1), \Sigma(1,2), \dots, \Sigma(d-1,d)$ respectively, and contracts all the other irreducible components to rational normal curves (with the notation that a rational normal scroll of type $\Sigma(k_1,k_2)$ is built on two rational normal curves of degrees $k_1$ and $k_2$ which span two supplementary linear spaces). Two successive scrolls $\Sigma(k,k+1)$ and $\Sigma(k+1,k+2)$ in the previous list are glued along the rational normal degree $k+1$ curves $\Lambda_{r_k}$ and $\Lambda_{r_{k+1}-1}$.

The upshot is that the restriction to the fibres of $\pi_r : \mathcal{P}_r \to \mathbf{D}$ of the twist of $\mathcal{L}_{r,d}$ under consideration gives a degeneration of the $d$-th Veronese image of $\mathbf{P}^2$ to a chain of rational normal scrolls. Correspondingly, we see limits of plane degree $d$ curves as hyperplane sections of this chain of rational normal scrolls.

As a sanity check, let us verify that the degrees and embedding dimension of such a projective model are correct. The total degree is

$$\sum_{k=0}^{d-1} \deg\big(\Sigma(k,k+1)\big) = \sum_{k=0}^{d-1}(2k+1) = d^2,$$

and the dimension of the projective linear space spanned by our chain of scrolls is

$$\sum_{k=0}^{d-1} h^0(\mathcal{O}_{\mathbf{P}^1}(k+1)) = \sum_{k=0}^{d-1}(k+2) = 2d + \sum_{k=1}^{d-1} k = \sum_{k=1}^{d+1} k - 1,$$

which indeed coincide respectively to the degree and embedding dimension of $d$-th Veronese image of $\mathbf{P}^2$.

We will be able to give a nice graphic representation of this degeneration of the $d$-th Veronese to a chain of rational normal scrolls when we know some toric and tropical geometry, see [X, (8.4)].

**(5.3) Limits of nodes.** It is a central theme of this set of notes (in particular, see [I, $**$]) that the limit of a node may be of the two following types.

(a) It may be a node on one of the irreducible components of the degeneration. In our situation the irreducible components of the central fibre are rational normal scrolls, so the only way a hyperplane section has a node is that it splits off a fibre.

This is indeed what we have observed in the previous examples. See in particular (4.2.2): nodes as above are indicated by a black square. There is such a node on Figures 9, 10, and 11.

(b) It may be a tacnode along the transverse intersection curve of two irreducible components of the central fibre (or, more generally, $\mu - 1$ nodes may degenerate together to a $\mu$-tacnode). In our situation in this lecture, this phenomenon occurs when two curves of respective classes $H + kF$ and $H + (k+1)F$ are linked by a chain of fibres of multiplicity $\mu$.

This happens for instance on Figure 8. In the projective model given by the appropriate twist for the curve displayed on Figure 8, $S_5$ is contracted to a conic, and the two scrolls images of $S_6$ and $S_4$ intersect along this conic; the two curves on $S_6$ and $S_4$ are then cut out by a hyperplane which is tangent to this conic. Thus, there is an ordinary tacnode, and one can find the correct multiplicity 4 for this limit curve as follows.

A multiplicity $\mu$ chain of fibres (of type (a) in the notation of (3.8), say) contributes enumeratively with multiplicity $\mu^2$, whereas a $\mu$-tacnode contributes with multiplicity $\mu$ in the count of $(\mu - 1)$-nodal curves. The principle to obtain the multiplicity $\mu^2$ in the tacnodal context is the following: the points at which one can have a $\mu$-tacnode should be searched within a $g_\mu^{\mu-1}$ on a smooth rational curve, and there are $\mu$ divisors of the form $\mu.p$ in such a linear series (for instance, this may be seen by applying the de Jonquières formula, see [XII, Section $**$]). Then, each of those $\mu$ possibilities counts with multiplicity $\mu$, thus giving the total contribution $\mu^2$ in the end, as required.

It is nevertheless not obvious how to understand in full generality the limits of the nodes along the degeneration considered in the present lecture in terms of the theory presented in [I]. For instance, one ought to deal with situations as in Figure 12 below, where the multiple chain of fibres is not entirely contracted in the appropriate projective model. The curve on Figure 12 is a limit of rational plane quartics, hence carries the limits of three nodes. One node is amounted for by the tacnode, and the two others by the intersection point of the double fibre with the curve of class $H$ on $S_7$. In order to understand this situation within the theory of [I], it is necessary to consider an alternative model of the degeneration, in such a way that one gets rid of the multiple curves.

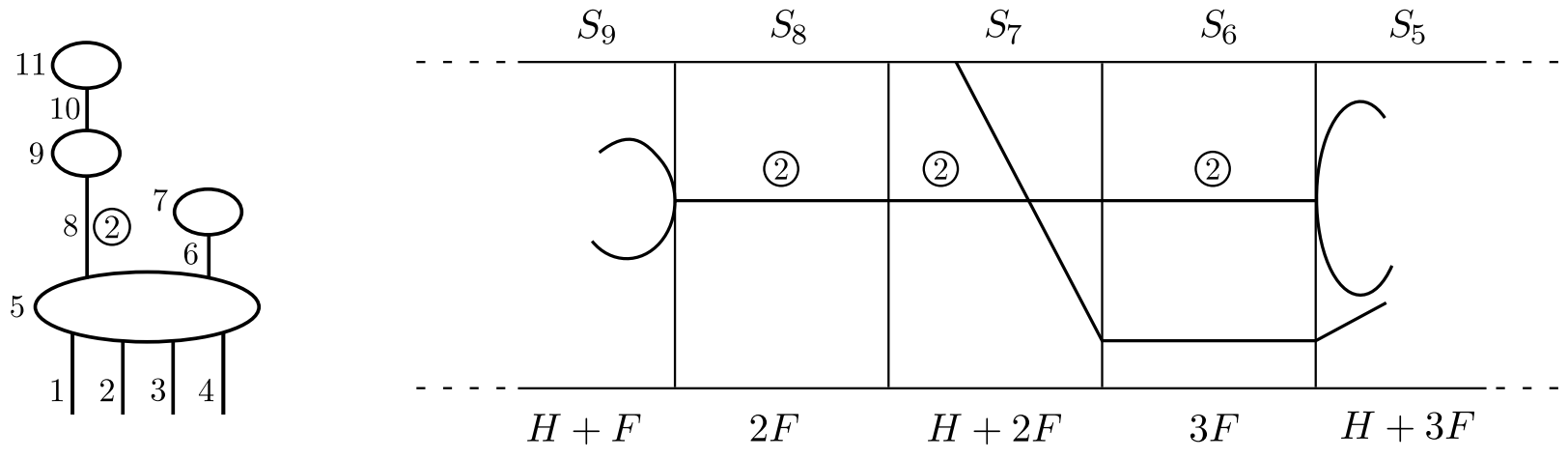

Figure 12: Limit of a rational plane quartic with multiple fibres (detail)

## References

### Chapters in this Volume[5]

### External References

[5] see `https://www.math.univ-toulouse.fr/~tdedieu/#enumTV`

# Lecture X
# Enumeration of curves in toric surfaces, tropical geometry, and floor diagrams

## by Thomas Dedieu

The theme of this text is an enumerative problem for curves on toric surfaces, analogous to the one studied in previous chapters for the projective plane, and its interplay with related enumerative problems in different contexts. Specifically, we will discuss an equivalent enumerative problem in the framework of tropical geometry, as well as its translation, under favourable technical circumstances, into a purely combinatorial enumerative problem involving ad hoc objects called floor diagrams. The text also includes an introduction to both toric varieties and tropical geometry, intended for readers with absolutely no prior knowledge of these topics (e.g., the author himself).

Toric surfaces are compact surfaces containing a torus $(\mathbf{C}^*)^2$, the complement of which is an anticanonical divisor made of a cycle of rational curves. The enumerative problem we are interested in consists in counting curves of given genus $g$ and Newton polygon $\Delta$, with prescribed contact pattern with the anticanonical cycle, both at assigned and unassigned points, and passing through the appropriate number of general base points; the shape of the polygon $\Delta$ determines which toric surface the curves live in, while the lengths of its sides determine their linear equivalence class. Problems of this kind naturally arise when one studies an enumerative problem by degeneration as discussed in [I], much in the same way that logarithmic Severi varieties of the pair $(\mathbf{P}^2, \text{line})$ show up when one studies plain Severi varieties of the projective plane, as explained in [VI].

Plane tropical[1] curves are piecewise linear objects in the real plane $\mathbf{R}^2$ (equipped with the integral structure determined by $\mathbf{Z}^2$), which appear as a sort of projected shadows of algebraic curves in $(\mathbf{C}^*)^2$, following seminal ideas of Kapranov [12] and others. Of course much information is lost in the projection, but it is quite miraculous how good tropical curves are at representing algebraic curves. In particular, it has been shown by Mikhalkin, following a suggestion by Kontsevich (see [20]), that the aforementioned enumerative problem for curves in toric surfaces is equivalent to its tropical shadow.

An advantage of tropical curves over algebraic curves is that they are much more hands-on: it is possible to play around with them using a ruler and a sheet of gridded paper without losing any piece of information. The tropical enumerative problem is therefore more tractable than the original toric problem, and in principle it can be solved with elementary considerations only: we shall see several instances of this in the text. Yet, of course, when the size of the problem gets too large, it is no longer possible to solve it with bare hands.

It has been realized by Brugallé and Mikhalkin that stretching the base points of the tropical enumerative problem along a chosen direction in $\mathbf{R}^2$ helps structuring its solution. In particular, under mildly favourable technical conditions, if the base points are far enough one from another along this direction, then for all the tropical curves one wants to count, the base points are distributed one by one on certain "components"; these "components" are of two types, floors and elevators, and their incidence relations may be neatly encoded in a combinatorial device called a floor diagram. Then, the tropical enumerative problem translates into a problem of enumeration of floor diagrams enhanced with a marking, the latter encoding the distribution of the base points on the various components. The tropical enumerative problem may also be solved using the so-called lattice path algorithm, which had been given by Mikhalkin already in [23], but we will not discuss it here.

To conclude this introduction, let me mention some aspects and developments that will not be discussed any further in the text. Quite remarkably, tropical geometry and floor diagrams also allow the study of enumerative problems in real algebraic geometry. In particular they have permitted the computation of Welschinger invariants, which count real rational curves with a sign depending on the real nature of their singularities (see the nice survey [5]). It has to be noted that the correspondence between (algebraic) curves and floor diagrams counting can be established without resorting to tropical geometry, by degeneration, as has been discovered by Brugallé [4] (this is explained in [IX]). On the other hand, tropical geometry gives a framework in which it is possible to define new enumerative invariants. For instance, one can define tropical Welschinger invariants in positive genus, whereas it is proven that the corresponding real algebraic enumerative problem is not invariant under the choice of the general base points. One can also count plane tropical curves with multiplicities that are rational fractions in the (quantum) indeterminate $q$, which leads to the so-called Block–Göttsche invariants; the latter are polynomials in $q$ interpolating between Welschinger (given by $q = -1$) and Gromov–Witten (given by $q = 1$) invariants; the enumerative interpretation of their coefficients is still unclear. About these new invariants, see [5] and the references therein.

The formulation as an enumerative problem about floor diagrams has enabled Fomin and Mikhalkin [13] to prove various nice results about complex algebraic plane curves. In particular, they have obtained a proof of the Di Francesco–Itzykson conjecture, to the effect that for all $\delta \geqslant 0$, there exists a degree $2\delta$ polynomial $P_\delta \in \mathbf{Q}[d]$ such that for all large enough $d$, $P_\delta(d)$ is the number of degree $d$ plane curves with $\delta$ nodes passing through the appropriate number of general base points. This is a particular case of the Göttsche conjecture, and it is related to the results discussed in [XIII] and [XIV].

[1] The term 'tropical' is empty from mathematical intuition: it refers to Imre Simon who resided in Brazil, see [21, p. 283].

Floor diagrams have been used by Brugallé and Mikhalkin [6] to give an equivalent formulation of the problem of enumerating rational curves in a complex projective space of arbitrary dimension. Also, adapted versions of floor diagrams have been applied to the study of Gromov–Witten invariants with (stationary) descendant insertions of Hirzebruch surfaces, see [3] and [9] (the descendant insertions mean that we are imposing more elaborate conditions than the passing through base points in order to arrive at a finite number of curves).

For this presentation I have drawn from various sources, to which I owe a great intellectual debt. The general presentation of tropical geometry follows rather closely a mini course given by Ilia Itenberg, *Géométrie tropicale énumérative*, at the 2010 GAGC meeting at CIRM. The presentation of plane and abstract tropical curves, and of the tropical enumerative problem, exists only thanks to personal lectures offered by Erwan Brugallé to the author in Jussieu, 2016. The presentation of floor diagrams and their enumeration follows Brugallé and Mikhalkin's article [8]. Finally, I wish to mention a wonderful talk by Kontsevich [19], which in principle was targeted for master students, but contains illuminating insights on tropical geometry and its interplay with algebraic geometry.

There are several texts written by experts of tropical geometry which the reader can consult: [7] is written for a broad audience of mathematicians, and does not discuss enumerative geometry; [16] does discuss enumerative geometry, and targets an audience of graduate students; [5] and [18] are intended for geometers. I also recommend the somehow historical [27].

The text starts by a quick introduction to toric surfaces (in Section 1), based on my personal down-to-earth understanding of it (but supported by the reference [11]) and guaranteed free of fans, after what we state the enumerative problem for curves on toric surfaces (in Section 2.2). Section 2.3 gives several examples, showing how the toric enumerative problem in question arises when one studies, say, the enumerative geometry of $K3$ surfaces by degeneration to unions of planes. Plane tropical curves are introduced in Section 5 in the way we need them in order to state and study the relevant tropical enumerative problem. Before that, I briefly explain in Section 3 the ideas that led to the definition of these objects. Section 4 is a parenthesis about the Legendre–Fenchel transformation (often merely called Legendre transformation), which is a useful tool in the study of corner loci and tropical hypersurfaces. It is possible to skip Sections 3 and 4, as Section 5 should be self-contained, and sufficient to understand the remainder of the text. In Section 6 we state the tropical enumerative problem (using the definitions of Section 5) and explain the correspondence with the toric enumerative problem. Finally, Section 7 is dedicated to the translation of the tropical enumerative problem in terms of floor diagrams. The last Section 8 is a sort of appendix, in which we develop the idea that tropical curves are great for the graphic representation of algebraic curves; we will see in particular that the graphic representation of the degeneration procedure of [IX] (which is an alternative path to floor diagrams) makes for striking revelations.

**Thanks.** I am grateful to Thomas Blomme and Ilia Itenberg for their comments on a preliminary version of this text, and to Erwan Brugallé for his essential help in the conception of this text.

**(0.1) General notation.** We use the notation of [III, ∗∗] and [VI, ∗∗]. In particular $\underline{\mathbf{N}}$ denotes the set of sequences $\alpha = [\alpha_1, \alpha_2, \ldots]$ of non-negative integers with only finitely many $\alpha_i$ non-zero, and for all $\alpha \in \underline{\mathbf{N}}$,

$$|\alpha| = \sum_i \alpha_i, \qquad I\alpha = \sum_i i\,\alpha_i, \qquad \text{and} \qquad I^\alpha = \prod_i i^{\alpha_i}.$$

For any vector $u$ in the real plane $\mathbf{R}^2$, we denote by $u^\perp$ the vector obtained by rotating $u$ by $+\pi/2$. A vector $u = (a, b) \in \mathbf{Z}^2$ is *primitive* if $a$ and $b$ are relatively prime.

On the projective plane $\mathbf{P}^2$, we let $H$ denote the line class. On Hirzebruch ruled surfaces $\mathbf{F}_e$, we let $E$ be the exceptional section, $F$ the class of fibres of the ruling, and $H = E + eF$.

# 1 – Toric projective surfaces

## 1.1 – The prototypical examples

Let $S$ be a surface equipped with a system of (pluri) homogeneous coordinates, and fix a (multi) degree $d$. The projectivization of the space $V_d$ of all homogeneous polynomials of degree $d$ is a complete linear system of curves on $S$. The notion of *Newton polygon* enables one to define linear subsystems of $\mathbf{P}V_d$ in a remarkably meaningful way, as follows.

Degree $d$ monomials form a basis of $V_d$, and may be represented as lattice points $\mathsf{p}_1, \dots, \mathsf{p}_{N_d}$ in some real affine space (by *lattice points*, or *integral points*, we mean points with integral coordinates). Consider a convex polygon $\Delta$, the vertices of which belong to $\{\mathsf{p}_1, \dots, \mathsf{p}_{N_d}\}$. One lets $V_\Delta$ be the subspace of $V_d$ generated by those monomials whose corresponding lattice point $\mathsf{p}$ lies in $\Delta$ (thus, we see $\Delta$ as a 2-dimensional convex domain). We will see that the subsystem $\mathbf{P}V_\Delta$ of $\mathbf{P}V_d$ may alternatively be defined in a natural way by base points, possibly infinitely near. This is best understood by looking at a couple of examples.

**(1.1) Example** ($\mathbf{P}^2$)**.** It comes with a system of homogeneous coordinates $(x : y : z)$. Fix a degree $d \in \mathbf{N}$; to each monomial $x^a y^b z^c$ of degree $d$, we associate the lattice point $(a, b) \in \mathbf{R}^2$ (of course this determines $c = d - a - b$). Degree $d$ monomials then correspond bijectively to points with integral coordinates in the triangle

$$T_d \stackrel{\text{def.}}{=} \{(a, b) \in \mathbf{R}^2 : a, b \geqslant 0 \text{ and } a + b \leqslant d\}.$$

We may for instance consider the two Newton polygons pictured below.

(a) 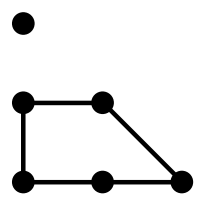 (b) 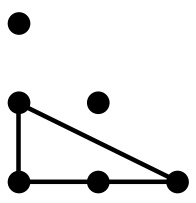

In case (a), the elements of $V_\Delta$ are the quadratic polynomials of the form

$$f_2(x, z) + y f_1(x, z),$$

with $f_2$ and $f_1$ homogeneous polynomials of degrees 2 and 1, respectively. Thus, $\mathbf{P}V_\Delta$ is the linear system of conics passing through the point $(0 \!:\! 1 \!:\! 0)$.

In case (b), $V_\Delta$ consists of the polynomials of the form

$$f_2(x, z) + \alpha \cdot yz,$$

with $f_2$ a degree 2 homogeneous polynomial, and $\alpha \in \mathbf{C}$. Thus, $\mathbf{P}V_\Delta$ is the linear system of conics passing by the point $(0 \!:\! 1 \!:\! 0)$, and tangent to the line $z = 0$ there (if you set $z = 0$ in the above expression, only the $x^2$ term does not vanish).

**(1.2) Example** ($\mathbf{P}^1 \times \mathbf{P}^1$)**.** It comes with a system of bihomogeneous coordinates $((x_0 : x_1), (y_0 : y_1))$. Fix a bidegree $d = (d_1, d_2) \in \mathbf{N}^2$; to each monomial $x_0^{a_0} x_1^{a_1} y_0^{b_0} y_1^{b_1}$ of bidegree $d$, we associate the lattice point $(a_1, b_1) \in \mathbf{R}^2$. Bidegree $d$ monomials then correspond bijectively to points with integral coordinates in the rectangle

$$R_{d_1,d_2} \overset{\text{def.}}{=} \{(a,b) \in \mathbf{R}^2 : a, b \geqslant 0,\ a \leqslant d_1,\ \text{and } b \leqslant d_2\}.$$

On the picture below, the polygon (a) (resp. (b)) corresponds to the linear system of hyperplane sections (resp. of sections by quadrics) of $\mathbf{P}^1 \times \mathbf{P}^1$ embedded as a smooth quadric surface in $\mathbf{P}^3$.

(a) 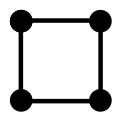 (b) 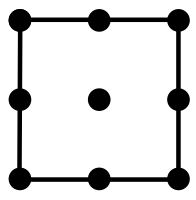

Observe that the polygon (a) may also be realized inside the triangle $T_2$ of Example (1.1). In this incarnation, it corresponds to the linear system of plane conics with the two base points $(1:0:0)$ and $(0:1:0)$. The map associated to this linear system is birational from $\mathbf{P}^2$ to a smooth quadric surface. It follows that the two incarnations of the polygon (a) give the same linear system of curves.

## 1.2 – The general construction

**(1.3)** Let $\Delta \subseteq \mathbf{R}^2$ be a convex polygon with vertices at integral points (this is called a *lattice polygon*), and let $\mathsf{p}_i = (a_i, b_i)$, $i = 1, \dots, N_\Delta$, be the integral points of $\mathbf{R}^2$ lying inside $\Delta$. We associate the following map to $\Delta$:

$$\begin{aligned} \phi_\Delta : \quad & (\mathbf{C}^*)^2 \to \mathbf{P}^{N_\Delta - 1} \\ & (t_1, t_2) \mapsto \big(t_1^{a_1} t_2^{b_1} : \dots : t_1^{a_{N_\Delta}} t_2^{b_{N_\Delta}}\big). \end{aligned}$$

Then, we define $X_\Delta \subseteq \mathbf{P}^{N_\Delta - 1}$ to be the Zariski closure in $\mathbf{P}^{N_\Delta - 1}$ of the image of $\phi_\Delta$; it is a projectively normal surface, see [11, Thm. 2.4.1, Def. 2.3.14, and Cor. 2.2.13].

Thus, the object associated to $\Delta$ is a normal surface in a given projective embedding; equivalently, it is a surface together with a linear system of Cartier divisors. It is possible to build a combinatorial object from $\Delta$ that forgets about the embedding, which is called a *fan*; it is usually the preferred object in toric geometry. In this text however, our primary interest resides in linear systems (hence in polytopes), and we shall not discuss fans.

Quite remarkably, the surface $X_\Delta$ thus constructed may always be equipped with a system of homogeneous coordinates, see [11, Chapter 5], so that we can play the game of Section 1.1 with it, and associate the complete linear system $|\mathcal{O}_{X_\Delta}(1)|$ with a polygon. The latter polygon is the very $\Delta$ we started with.

**(1.4) Remark.** One may rephrase the above definition in the terms of Section 1.1. For instance, any $\Delta$ as above may be realized inside a triangle $T_d$ as in (1.1), for sufficiently large $d$ (note that the definition of $\phi_\Delta$ is invariant under translation of $\Delta$). Then, $\Delta$ determines a linear system of degree $d$ plane curves of dimension $N_\Delta - 1$, and $\phi_\Delta$ equals the restriction to

$$(\mathbf{C}^*)^2 \simeq \big\{(t_0 : t_1 : t_2) : t_0 t_1 t_2 \neq 0\big\} \subseteq \mathbf{P}^2$$

of the rational map $\mathbf{P}^2 \dashrightarrow \mathbf{P}^{N_\Delta - 1}$ associated with this linear system.

Similarly, $\Delta$ may always be realized inside a rectangle $R_{d_1,d_2}$ as in (1.2), and $\phi_\Delta$ identified with the restriction to

$$(\mathbf{C}^*)^2 \simeq \big\{(s_1 : t_1), (s_2 : t_2)) : s_1 t_1 s_2 t_2 \neq 0\big\} \subseteq \mathbf{P}^1 \times \mathbf{P}^1$$

of the rational map $\mathbf{P}^1 \times \mathbf{P}^1 \dashrightarrow \mathbf{P}^{N_\Delta - 1}$ associated with the corresponding linear system.

We thus see that the construction in Section 1.1 does not depend on the choice of the surface $S$ at the beginning, but only on the polygon $\Delta$.

**(1.5) Equations.** One may compute the equations of $X_\Delta$ in $\mathbf{P}^{N_\Delta-1}$ by elimination from a parametrization of the affine cone $\widehat{X}_\Delta$ over $X_\Delta$. To do so, beware that the map

$$\widehat{\phi}_\Delta : (t_1,t_2) \in (\mathbf{C}^*)^2 \mapsto \big(t_1^{a_1}t_2^{b_1} : \ldots : t_1^{a_{N_\Delta}}t_2^{b_{N_\Delta}}\big) \in (\mathbf{C}^*)^{N_\Delta}$$

of course does not parametrize (an open subset of) the affine cone $\widehat{X}_\Delta$. To obtain such a parametrization, one may apply the trick described in [11, top of p. 58], and consider instead the map

$$(t_0,t_1,t_2) \in (\mathbf{C}^*)^3 \mapsto \big(t_0t_1^{a_1}t_2^{b_1} : \ldots : t_0t_1^{a_{N_\Delta}}t_2^{b_{N_\Delta}}\big) \in (\mathbf{C}^*)^{N_\Delta}.$$

Then, the equations of $X_\Delta$ in $\mathbf{P}^{N_\Delta-1}$ may be obtained by eliminating $(t_0,t_1,t_2)$ from the ideal

$$\big(x_1 - t_0t_1^{a_1}t_2^{b_1}, \ldots, x_{N_\Delta} - t_0t_1^{a_{N_\Delta}}t_2^{b_{N_\Delta}}\big) \subseteq \mathbf{C}[t_0,t_1,t_2] \otimes \mathbf{C}[x_1,\ldots,x_{N_\Delta}].$$

Alternatively, one may consider the map

(1.5.1) $$(t_0,t_1,t_2) \in (\mathbf{C}^*)^3 \mapsto \big(t_0^{d-a_1-b_1}t_1^{a_1}t_2^{b_1} : \ldots : t_0^{d-a_{N_\Delta}-b_{N_\Delta}}t_1^{a_{N_\Delta}}t_2^{b_{N_\Delta}}\big) \in (\mathbf{C}^*)^{N_\Delta},$$

corresponding to the extension of $\phi_\Delta$ to a rational map $\mathbf{P}^2 \dashrightarrow \mathbf{P}^{N_\Delta-1}$ discussed in Remark (1.4). It does parametrize an open subset of the affine cone $\widehat{X}_\Delta$ by [11, Prop. 2.1.4], since the linear functional on monomials $t_0^ct_1^at_2^b \mapsto a+b+c$ takes the value $d$ on all monomials appearing in the map (1.5.1). Then, one may obtain the equations of $X_\Delta$ in $\mathbf{P}^{N_\Delta-1}$ by eliminating $(t_0,t_1,t_2)$ from the ideal

$$\big(x_1 - t_0^{d-a_1-b_1}t_1^{a_1}t_2^{b_1}, \ldots, x_{N_\Delta} - t_0^{d-a_{N_\Delta}-b_{N_\Delta}}t_1^{a_{N_\Delta}}t_2^{b_{N_\Delta}}\big)$$

as above.

Similarly, one may consider the map

(1.5.2) $$(s_1,s_2,t_1,t_2) \in (\mathbf{C}^*)^4 \mapsto \big(s_1^{d_1-a_1}s_2^{d_2-b_1}t_1^{a_1}t_2^{b_1} : \ldots : s_1^{d_1-a_{N_\Delta}}s_2^{d_2-b_{N_\Delta}}t_1^{a_{N_\Delta}}t_2^{b_{N_\Delta}}\big) \in (\mathbf{C}^*)^{N_\Delta}$$

corresponding to the extension of $\phi_\Delta$ to a rational map $\mathbf{P}^1 \times \mathbf{P}^1 \dashrightarrow \mathbf{P}^{N_\Delta-1}$ discussed in Remark (1.4), which parametrizes an open subset of the affine cone $\widehat{X}_\Delta$ by [11, Prop. 2.1.4], since the two linear functionals $s_1^{a'}s_2^{b'}t_1^at_2^b \mapsto a+a'$ and $s_1^{a'}s_2^{b'}t_1^at_2^b \mapsto b+b'$ take the values $d_1$ and $d_2$, respectively, on all monomials appearing in the map (1.5.2).

**(1.6) Example.** The trapezium

$$\mathrm{Tz}_{a,b}^r = \big\{(x,y) \in \mathbf{Z}^2 : 0 \leqslant y \leqslant a \text{ and } 0 \leqslant x \leqslant b + (r-y)a\big\}$$

corresponds to the complete linear system $|aH + bF|$ on the Hirzebruch surface $\mathbf{F}_r$, where $F$ is the class of a fibre, and $H = E + rF$ with $E$ the class of the section with self-intersection $-r$.

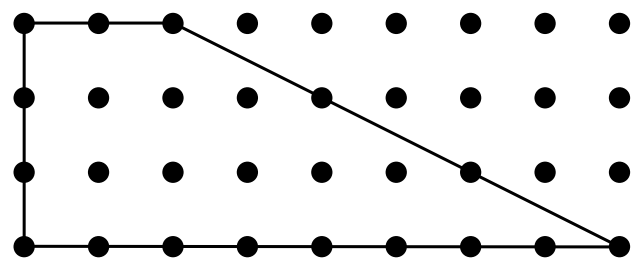

Figure 1: The trapezium $\mathrm{Tz}_{3,2}^2$: $|3H + 2F|$ on $\mathbf{F}_2$

## 1.3 – Geometric properties and numerical invariants

Consider a surface $X_\Delta$ constructed from a convex integral polygon $\Delta \subseteq \mathbf{R}^2$ as in (1.3).

**(1.7) Anticanonical cycle.** The complement $X_\Delta \setminus \phi_\Delta((\mathbf{C}^*)^2)$ is a cycle of rational curves, and it is a Weil divisor belonging to the anticanonical linear equivalence class $-K_{X_\Delta}$, see [11, Theorem 8.2.3]; we shall call it the *anticanonical cycle* of $X_\Delta$, and denote it by $\daleth_\Delta$.

Being a cycle of rational curves means that $\daleth_\Delta$ is a union of rational curves $C_1, \ldots, C_n$, where each $C_i$ intersects $C_{i-1}$ and $C_{i+1}$ transversely[2] at one point (with the convention that $C_0 = C_n$ and $C_{n+1} = C_1$) and no other curve.

The curves $C_1, \ldots, C_n$ correspond to the edges $\Sigma_1, \ldots, \Sigma_n$ (labelled in a cyclic way) of the polygon $\Delta$, and the intersection points $C_i \cap C_{i+1}$ to its vertices. In the projective embedding given by $\Delta$, the curve $C_i$ is a rational normal curve of degree $\ell_{\mathbf{Z}}(\Sigma_i)$, which we define as the *integral length* of the segment $\Sigma_i$, i.e., the number of integral points lying on $\Sigma_i$ minus one.

**(1.8) Degree and number of sections.** We let $H_\Delta$ be the class $c_1(\mathcal{O}_{X_\Delta}(1))$, corresponding to the embedding defined by $\Delta$. The dimension $h^0(\mathcal{O}_{X_\Delta}(1))$ is the number $N_\Delta$ of integral points lying in $\Delta$ minus one. One may think of $N_\Delta$ as an integral measure of the area of $\Delta$, in analogy with the integral length $\ell_{\mathbf{Z}}$, but beware that the integral area of $\Delta$ is usually defined as $2\mathscr{A}(\Delta)$, in the notation introduced in (1.9) below.

The anticanonical degree of $H_\Delta$ is

$$-K_{X_\Delta} \cdot H_\Delta = \daleth_\Delta \cdot H_\Delta = \textstyle\sum_{i=1}^n \ell_{\mathbf{Z}}(\Sigma_i) = \mathrm{Card}(\partial\Delta_{\mathbf{Z}}),$$

where $\partial\Delta_{\mathbf{Z}}$ stands for the set of integral points in the boundary of $\Delta$, see [11, Proposition 10.5.6].

**(1.9) Arithmetic genus.** The arithmetic genus of the divisors belonging to the linear system $|\mathcal{O}_{X_\Delta}(1)|$ is

$$p_a(\Delta) \overset{\text{def.}}{=} p_a(H_\Delta) = \mathrm{Card}(\mathring{\Delta}_{\mathbf{Z}}),$$

where $\mathring{\Delta}_{\mathbf{Z}}$ stands for the set of integral points in the interior of $\Delta$, see [11, Proposition 10.5.8].

Equivalently, one has

$$(H_\Delta)^2 = 2\mathscr{A}(\Delta),$$

where $\mathscr{A}(\Delta)$ is the *Euclidean area* of $\Delta$, taking the area of the quadrilateral defined by any lattice basis of $\mathbf{R}^2$ as the unit, [11, Proposition 10.5.6].

With the two above identities, the adjunction formula

$$2p_a(H_\Delta) - 2 = (K_{X_\Delta} + H_\Delta) \cdot H_\Delta$$

amounts to Pick's formula, see [11, Example 9.4.4], which expresses the area of a lattice polygon in terms of its lattice points.

**(1.10) Singularities.** The surface $X_\Delta$ is smooth except (possibly) at the points corresponding to the vertices of the anticanonical cycle $\daleth_\Delta$, where it has cyclic quotient singularities, the types of which are determined as follows (see [11, Section 10.1]).

Consider a vertex $v$ of $\Delta$, and let $u'_v$ and $u''_v$ be the two primitive integral vectors giving the directions of the two edges of $\Delta$ intersecting at $v$ and pointing away from $v$, in such a way that the angle $(u'_v, u''_v)$ points inside $\Delta$. Then the two vectors $(u'_v)^\perp$ and $-(u''_v)^\perp$ are the normal vectors of the two edges at $v$, pointing towards the interior of $\Delta$ (see Figure 2 for an example).

[2]two consecutive curves $C_i$ and $C_{i+1}$ of $\daleth_\Delta$ certainly intersect transversely in $\mathbf{P}^{N_\Delta-1}$, because they span two linear subspaces $\mathbf{P}^{\ell_i}$ and $\mathbf{P}^{\ell_{i+1}}$ which intersect at one point only. However, the surface $X_\Delta$ may be singular at the point $C_i \cap C_{i+1}$, and the intersection number $C_i \cdot C_{i+1}$ differs from 1 in general.

It is always possible to find an integral basis $(e_1, e_2)$ of $\mathbf{Z}^2$ such that the two vectors $(u'_v)^\perp$ and $-(u''_v)^\perp$ write

$$e_2 \quad \text{and} \quad de_1 - ke_2,$$

where $d$ and $k$ are relatively prime integers such that $d > 0$ and $0 \leqslant k < d$. Then the singularity of the surface $X_\Delta$ at the point corresponding to the vertex $v$ is of type $\frac{1}{d}(1,k)$, i.e., it is analytically locally isomorphic to the quotient of $\mathbf{C}^2$ by the group of $d$-th roots of unity acting as $\zeta.(x,y) = (\zeta x, \zeta^k y)$.

In particular, the Euclidean area of the parallelogram defined by $u'_v$ and $u''_v$, namely

$$\mathscr{A}\big(\langle u'_v, u''_v\rangle\big) = \det(u'_v, u''_v)$$

computes the order $d$ of the cyclic group, the action of which defines the singularity at the point of $X_\Delta$ corresponding to $v$. The surface is smooth at this point if $(u'_v, u''_v)$ forms a basis of $\mathbf{Z}^2$, and it has an ordinary double point if $\det(u'_v, u''_v) = 2$.

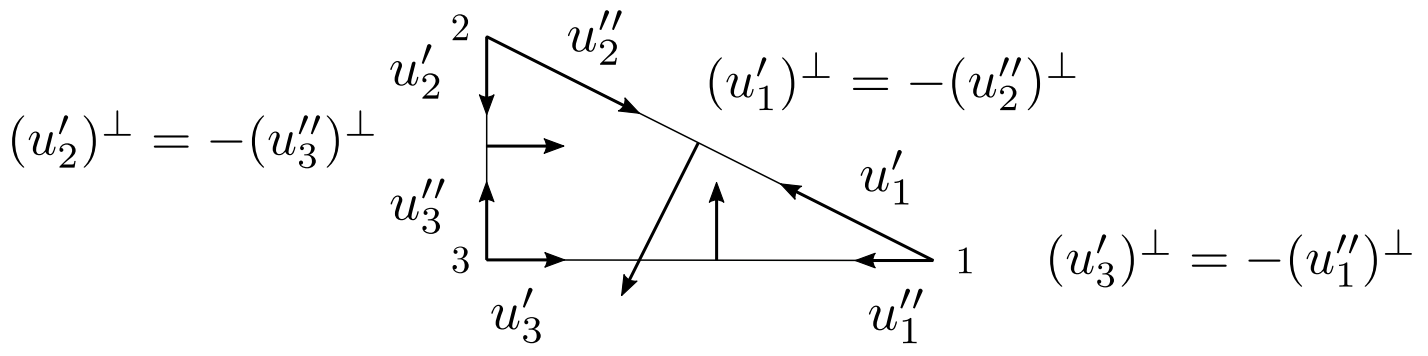

Figure 2: Determination of the types of singularities at the vertices: the triangle $\mathrm{Tz}^2_{3,0}$, $|3H|$ on $\mathbf{F}_2$

**(1.11) Example.** Consider $v$ the lower-right corner of the trapezium $\mathrm{Tz}^r_{a,b}$, see (1.6) and Figure 2 above. Then $u'_v = (-r, 1)$ and $u''_v = (-1, 0)$, hence $(u'_v)^\perp = (-1, -r)$ and $-(u''_v)^\perp = (0, 1)$, and we find that the corresponding point in the image of $\mathbf{F}_r$ by $|aH + bF|$ is a smooth point.

Consider on the other hand $v$ the upper vertex of the triangle $\mathrm{Tz}^r_{a,0}$. Then $u'_v = (0, -1)$ and $u''_v = (r, -1)$, hence $(u'_v)^\perp = (1, 0)$ and $-(u''_v)^\perp = (-1, -r)$. We may take $e_1 = (0, -1)$ and $e_2 = (1, 0)$, to the effect that

$$(u'_v)^\perp = e_2 \quad \text{and} \quad -(u''_v)^\perp = re_1 - e_2.$$

Thus the singularity at the point corresponding to $v$ is of type $\frac{1}{r}(1,1)$, which analytically locally is the vertex of the cone over a rational normal curve of degree $r$. And indeed the image of $\mathbf{F}_r$ by $|aH|$ is a generalized cone over a degree $r$ rational normal curve, and the point under consideration is its vertex.

# 2 – The toric enumerative problem

## 2.1 – Toric Cayley–Bacharach

The following may be considered as a toric Cayley–Bacharach property, and needs to be taken in consideration in order to properly state the toric enumerative problem. Later on, we will also see a tropical incarnation of this phenomenon, see Figure 10 and Figure 25.

**(2.1) Lemma.** *Let $Z$ be a subscheme of $\mathcal{X}_\Delta$ cut out by a (properly intersecting) hyperplane in $\mathbf{P}^{N_\Delta - 1}$. It imposes only* $\deg Z - 1$ *independent conditions to the linear system* $|\mathcal{O}_{X_\Delta}(1)|$*, i.e., the linear subsystem of divisors containing $Z$ has codimension* $\deg Z - 1$ *in* $|\mathcal{O}_{X_\Delta}(1)|$.

This follows from a basic cohomological argument, valid regardless of our toric setup, see also [III, Example **]. In the toric setup however, there is a nice elementary way of understanding it.

*Proof.* The lemma follows from the long exact sequence in cohomology induced by the restriction exact sequence,

$$0 \to \mathcal{O}_{X_\Delta}(K_{X_\Delta} + H) \to \mathcal{O}_{X_\Delta}(H) \to \mathcal{O}_{\partial_\Delta}(H) \to 0.$$

Indeed, $h^1(O_{X_\Delta}(K_{X_\Delta} + H)) = 0$ by Kawamata-Viehweg vanishing, and $h^0(\mathcal{O}_{\partial_\Delta}(H)) = \deg Z$ by Riemann–Roch, since $p_a(\partial_\Delta) = 1$. □

Let me now give an elementary argument showing why $Z$ cannot impose independent conditions to hyperplane sections of $X_\Delta$. Let $C_1, \ldots, C_n$ be the smooth rational curves forming the cycle $\partial_\Delta$. We shall exhibit an explicit relation between the roots of the polynomials $P|_{C_1}, \ldots, P|_{C_n}$, for all polynomials $P$ with Newton polygon $\Delta$, which implies that if $P$ vanishes at all the points of $Z$ but one, then it also vanishes at the last point.

Let $P$ be a polynomial with Newton polygon $\Delta$, and name its coefficients regarding the monomials on the boundary of $\Delta$ as indicated on the picture below.

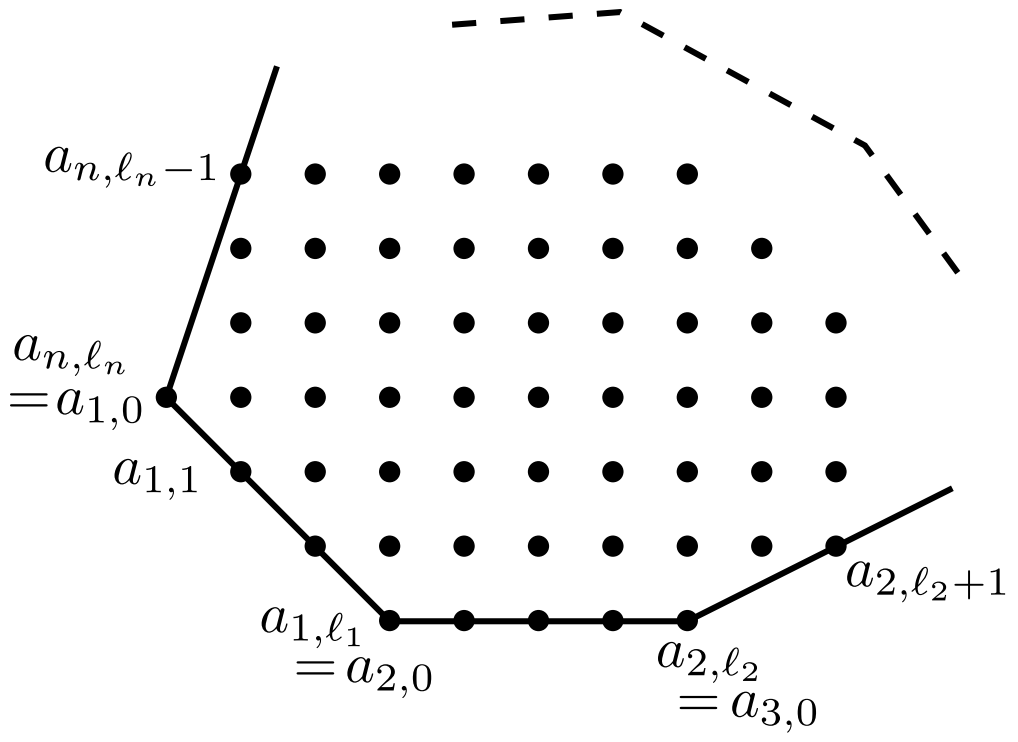

Thus, for all $i = 1, \ldots, n$, we have $P|_{C_i} = a_{i,0} + \cdots + a_{i,\ell_i} X^{\ell_i}$, hence the product of the roots of $P|_{C_i}$ equals $(-1)^{\ell_i} a_{i,0}/a_{i,\ell_i}$. Therefore, the product of the roots of $P|_{C_1}, \ldots, P|_{C_n}$ altogether equals

$$(-1)^{\ell_1+\cdots+\ell_n} \frac{a_{1,0}}{a_{1,\ell_1}} \cdots \frac{a_{n,0}}{a_{n,\ell_n}} = (-1)^{\ell_1+\cdots+\ell_n},$$

which is the announced relation.

## 2.2 – The enumerative problem

**(2.2)** Let $\Delta \subseteq \mathbf{R}^2$ be a convex integral polygon as in Section 1. We call $\Sigma_1, \ldots, \Sigma_n$ its edges, and correspondingly $C_1, \ldots, C_n$ the irreducible components of the anticanonical cycle $\partial_\Delta$.

Let $g$ be an integer such that $0 \leqslant g \leqslant p_a(\Delta) = \mathrm{Card}(\mathring{\Delta}_{\mathbf{Z}})$, and $\alpha = (\alpha_1, \ldots, \alpha_n), \beta = (\beta_1, \ldots, \beta_n) \in \underline{\mathbf{N}}^n$. We assume that for $s = 1, \ldots, n$,

$$I\alpha_s + I\beta_s = (C_s \cdot H_\Delta) = \ell_{\mathbf{Z}}(\Sigma_s),$$

and

$$I\alpha \stackrel{\text{def.}}{=} \textstyle\sum_{s=1}^n I\alpha_s < (\partial_\Delta \cdot H_\Delta) = \mathrm{Card}(\partial\Delta_{\mathbf{Z}}).$$

The problem is to find the number $N_g^\Delta(\alpha, \beta)$ of genus $g$ curves in the linear system $|H_\Delta|$ that (i) are the Zariski closure of a curve in $\phi_\Delta((\mathbf{C}^*)^2) \cong (\mathbf{C}^*)^2 \subseteq X_\Delta$, (ii) verify the contact conditions with $\mathcal{K}_\Delta$ at assigned and unassigned points fixed by $\alpha$ and $\beta$, respectively, and (iii) pass through the appropriate number of fixed general points in $X_\Delta$.

The precise setup is as follows. We fix $\underline{x} = \{x_1, \dots, x_{g-1+|\beta|}\}$ a set of

$$g - 1 + |\beta| = \mathrm{Card}(\partial\Delta_{\mathbf{Z}}) + g - 1 - I\alpha - (I\beta - |\beta|)$$

general points in $X_\Delta$, and $\Omega = (\{p_{s,i,j}\}_{1\leqslant j\leqslant \alpha_{s,i}})_{1\leqslant s\leqslant n, 1\leqslant i}$ a set of $\alpha$ general points on $\mathcal{K}_\Delta$, i.e., for $s = 1, \dots, n$, $\Omega_s = (\{p_{s,i,j}\}_{1\leqslant j\leqslant \alpha_{s,i}})_{1\leqslant i}$ is a set of $\alpha_s \in \underline{\mathbf{N}}$ general points on the curve $C_s$. We want the number of birational maps $\phi : C \to X_\Delta$, where $C$ is a reduced, nodal curve of arithmetic genus $g$, such that (i) $\phi(C)$ equals the Zariski closure of $\phi(C) \cap (X_\Delta \setminus \mathcal{K}_\Delta)$, (ii) there exist $\alpha$ points $q_{s,i,j}$ and $\beta$ points $r_{s,i,j}$ on $C$ such that $\phi(q_{s,i,j}) = p_{s,i,j}$ and

$$\phi^*\mathcal{K}_\Delta = \sum_{1\leqslant s\leqslant n} \sum_{1\leqslant i} \left( \sum_{1\leqslant j\leqslant \alpha_i} i\, q_{i,j} + \sum_{1\leqslant j\leqslant \beta_i} i\, r_{i,j} \right),$$

and (iii) $\phi(C)$ contains all $g - 1 + |\beta|$ points of $\underline{x}$.

**(2.3) Remarks.** This is exactly the same problem that we considered in [III] and subsequent chapters. In particular, [III, Theorem $**$] applies to our present situation.

The condition that $I\alpha < (\mathcal{K}_\Delta \cdot H_\Delta)$ is necessary because of Lemma (2.1): if we impose all contact points with $\mathcal{K}_\Delta$ to be assigned at general points, then there is no curve at all.

The condition (i) that the curves are the Zariski closures of their intersections with $X_\Delta \setminus \mathcal{K}_\Delta$ really means that we count curves that do not contain any irreducible component of the anticanonical cycle.

## 2.3 – Examples from degenerations of $K3$ surfaces

An emblematic situation where enumerative problems of the kind described in Section 2.2 occur is when one undertakes the enumeration of curves on $K3$ surfaces by considering type III degenerations, e.g., degenerations of $K3$ surfaces to unions of planes $X_0 \subseteq \mathbf{P}^N$ realizing a triangulation of the sphere $\mathbf{S}^2$ (by this we mean that the dual complex of $X_0$ is the dual complex of a triangulation of $\mathbf{S}^2$). Such a degeneration in general has various relevant birational models, see [I], all isomorphic one with another outside the central fibre. Typically, as we shall illustrate below, the irreducible components of the central fibres of these alternative degenerations are toric surfaces, with the anticanonical cycle cut out by the other components of the central fibre. The examples in (2.4), (2.5) and (2.6) below show up when one considers a degeneration of quartic hypersurfaces in $\mathbf{P}^3$ to a *tetrahedron*, i.e., the union of four general planes in $\mathbf{P}^3$ (equivalently, there exist homogeneous coordinates $(x:y:z:w)$ in which it is defined by the equation $xyzw = 0$). We refer to [10] for a detailed study, and to [I, Section $**$] for a brief discussion. The last example in (2.7) shows up when one considers a degeneration of octic $K3$ surfaces in $\mathbf{P}^5$ (i.e., smooth complete intersections of three quadrics) to an *octahedron*, i.e., the union of eight 2-planes in $\mathbf{P}^5$, the dual complex of which is a plain, Platonic, cube in $\mathbf{R}^3$ (equivalently, there exist homogeneous coordinates $(x:y:z:w:u:v)$ in which it is defined by the equations $xy = zw = uv = 0$).

**(2.4) A cubic surface with an anticanonical triangle.** The surface $X_\Delta \subseteq \mathbf{P}^3$ corresponding to the Newton polygon $\Delta$ below is the cubic hypersurface $xyz = w^3$ in homogeneous coordinates $(x:y:z:w)$. Its anticanonical cycle $ꓘ_\Delta$ is the triangle $xyz = w = 0$. The three sides of this triangle are the only lines contained in $X_\Delta$; the surface $X_\Delta$ has an $A_2$ double point at each vertex of $ꓘ_\Delta$, and is smooth elsewhere.

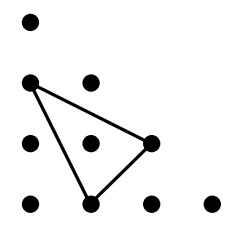

Let us verify the type of singularity at the vertex at the bottom, following the rule of (1.10). One has $u'_v = (1,1)$ and $u''_v = (-1,2)$, hence $(u'_v)^\perp = (-1,1)$ and $-(u''_v)^\perp = (2,1)$. We may take $e_1 = (0,1)$ and $e_2 = (-1,1)$, to the effect that

$$(u'_v)^\perp = e_2 \quad \text{and} \quad -(u''_v)^\perp = 3e_1 - 2e_2.$$

The singularity corresponding to this vertex is of type $\frac{1}{3}(1,2)$, i.e., an $A_2$ double point.

The linear system of hyperplane sections $|\mathcal{O}_{X_\Delta}(1)|$ may be realized, much in the same way as in Example (1.1), (b), above, as the linear system of plane cubics with six base points in the configuration indicated on the picture below. On the picture, two successive blow-ups of $\mathbf{P}^2$ at three points are drawn, with exceptional divisors $E_x, E_y, E_z$ and $E'_x, E'_y, E'_z$ respectively; base points are represented by black dots when they are proper, and arrows when they are infinitely near.

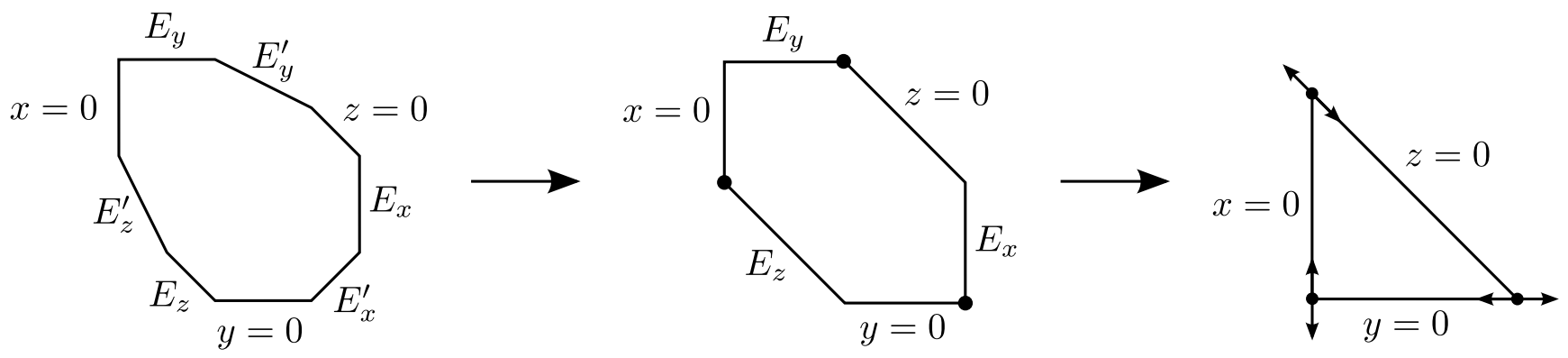

The three sides of $ꓘ_\Delta$ are the three $(-1)$-curves: $C_1 = E'_y$, $C_2 = E'_z$, and $C_3 = E'_x$. The cubic $X_\Delta$ is the main character of [I, Section $**$].

Alternatively, the surface $X_\Delta$ may be realized as a cyclic triple cover of $\mathbf{P}^2$ branched over a general triangle.

It turns out that, when considered as limits of hyperplane sections of neighbouring quartic $K3$ surfaces, sections of a tetrahedron $X_0 \subseteq \mathbf{P}^3$ by a hyperplane passing through one of its vertices $v$ are best substituted by hyperplane sections of the quartic $X_\Delta \cup \Pi_v \subseteq \mathbf{P}^3$, where $\Pi_v$ is the face of $X_0$ opposite to $v$, identified with the plane spanned by $ꓘ_\Delta$ in such a way that the triangle cut out on $\Pi_v$ by the three other faces of $X_0$ agrees with $ꓘ_v$. In order to see this, one has to realize the cubic $X_\Delta$ as some exceptional divisor living over $v$ in an appropriate blow-up of the total space of the degeneration to $X_0$. We refer to [I] for details.

Correspondingly, one wants to enumerate genus 0 hyperplane sections of $X_\Delta$ matching with the divisor cut out on the triangle $ꓘ_\Delta$ by a general line in the span $\langle ꓘ_\Delta \rangle$. Such a divisor is made of three points, one on each of the three sides of the triangle, that together impose only two independent linear conditions on $|\mathcal{O}_{X_\Delta}(1)|$, in line with Lemma (2.1).

To reduce to the enumerative problem (2.2), one may impose the passing through a fixed point on only two sides of $ꓘ_\Delta$, and let the intersection point with the third side move freely. For instance, one may set $\alpha = (1,1,0)$ and $\beta = (0,0,1)$. The upshot of Lemma (2.1) is that this third point, although free to move, will turn out to remain fixed. About this trick, see also [III, Paragraph $**$].

One may see that $N_0^\Delta(\alpha,\beta) = 3$ with the topological formula giving the number of singular members in a pencil of curves, see [XI, Section $**$], as follows. The two fixed points on $ꓘ_\Delta$ (or three aligned fixed points) define a pencil in $|\mathcal{O}_{X_\Delta}(1)|$; this is a pencil of plane cubics, hence it has twelve singular members, counted with multiplicities. There is a special member in this pencil that corresponds to $ꓘ_\Delta$ itself, incarnated as its pull-back to the minimal desingularization of $X_\Delta$, namely the 9-gon in the picture above; it counts with multiplicity 9 among the singular members of the pencil. It follows that there are $12 - 9 = 3$ irreducible 1-nodal curves in this pencil. We will confirm this number by solving the tropical counterpart of this enumerative problem in (6.11).

**(2.5)** Consider the surface $X_\Delta$ corresponding to the Newton polygon depicted below.

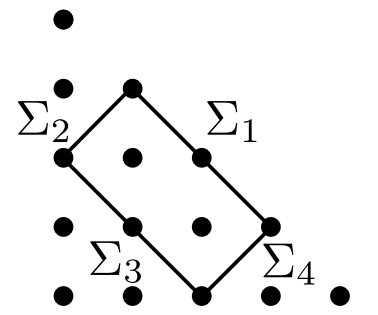

Its anticanonical cycle $\partial_\Delta$ consists of two conics $C_1, C_3$ corresponding to the sides $\Sigma_1, \Sigma_3$, and two lines $C_2, C_4$ corresponding to the sides $\Sigma_2, \Sigma_4$.

Along a degeneration of quartic $K3$ surfaces to a tetrahedron $X_0 \subseteq \mathbf{P}^3$, the sections of $X_0$ by a hyperplane containing one of its edges $W$ are best substituted by the hyperplane sections of $X_\Delta$ passing by two base points on $C_1$ and two base points on $C_3$, in order to be understood as limits of hyperplane sections of the $K3$ surfaces neighbouring $X_0$. The four base points on $C_1$ and $C_3$ are determined by the neighbouring quartic $K3$ surfaces, as in [I, Examples $**$ and $**$]; they impose four independent conditions to the linear system corresponding to $\Delta$, and therefore determine a 3-dimensional linear subsystem of $|\mathcal{O}_{X_\Delta}(1)|$. The latter can be realized as the linear system of plane quartics with base points as indicated on the picture below (with the same graphic convention as in (2.4); there are four simple and one double proper base points, and four infinitely near base points, all simple).

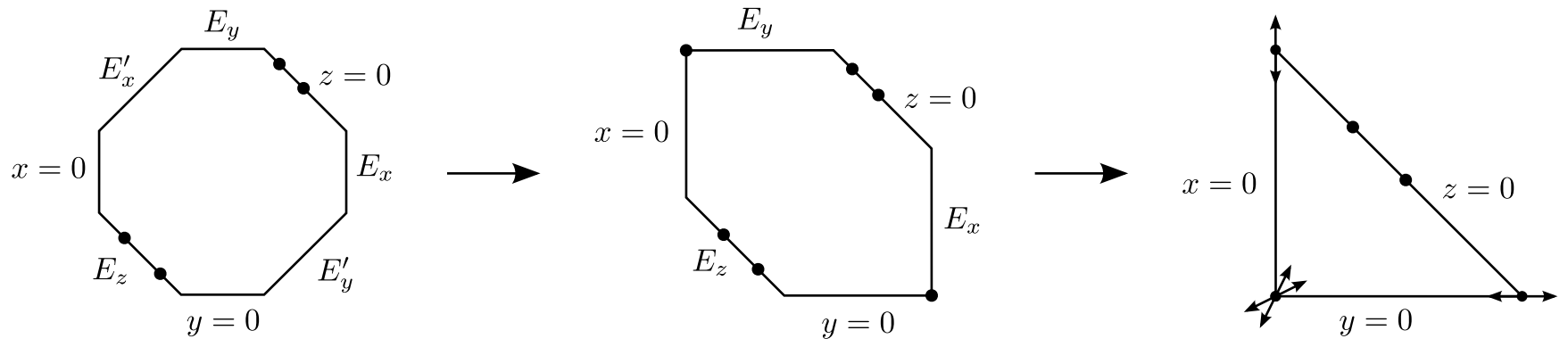

The projection map $\phi : X_\Delta \to \mathbf{P}^3$ from the four base points contracts the two conics $C_1$ and $C_3$ to two triple points, and identifies the two lines $C_2$ and $C_4$; its image is a rational quartic surface $X_W \subseteq \mathbf{P}^3$, with a double line (the image of both $C_2$ and $C_4$), two triple points on the latter (the images of $C_1$ and $C_3$, respectively), and no further singularity. The rational quartic surface $X_W \subseteq \mathbf{P}^3$ is one of the relevant degenerations of $K3$ surfaces alternative to the tetrahedron $X_0$; in this perspective, the double line of $X_W$ is reminiscent of the edge of $X_0$ opposite to $W$.

Correspondingly, one wants to enumerate genus 0 hyperplane sections of $\phi(X_\Delta) = X_W$ passing through a general point of the double line $\phi(C_2) = \phi(C_4)$. This amounts to adding two more base points on $C_2$ and $C_4$ having the same image by $\phi$. Equivalently, these two base points together with the four on $C_1$ and $C_3$ only impose five independent linear conditions to $|\mathcal{O}_{X_\Delta}(1)|$, again in line with Lemma (2.1). To set up our enumerative problem in the terms of (2.2), we fix the intersection with $C_2$, and let that with $C_4$ move freely; in other words, we set $\alpha = (2, 1, 2, 0)$ and $\beta = (0, 0, 0, 1)$.

By carefully analysing the geometry, one finds with the Plücker formulae that $N_1^\Delta(\alpha, \beta) = 10$ and $N_0^\Delta(\alpha, \beta) = 16$, see [10, Section 5.5]. We will come back to this with a tropical point of view in Example (7.18). For related floor diagram considerations, see Examples (7.5), (7.8), (7.11), and (7.17.2).

**(2.6)** Consider the 4-Veronese surface $X_\Delta = v_4(\mathbf{P}^2) \subseteq \mathbf{P}^{14}$, corresponding to the triangle $\Delta = T_4$ as in (1.1). The anticanonical cycle $\partial_\Delta$ is the image by $v_4$ of the triangle $xyz = 0$ in $\mathbf{P}^2$.

Along a degeneration of quartic $K3$ surfaces to a tetrahedron $X_0 \subseteq \mathbf{P}^3$, to the limits of hyperplane sections of the general surfaces tending to sections of $X_0 \subseteq \mathbf{P}^3$ by one of its faces $F$, there correspond hyperplane sections of $X_\Delta$ passing through twelve base points on $\text{ꓘ}_\Delta$. These base points are determined by the neighbouring surfaces as in [I, Examples ** and **], and together form the image by $v_4$ of the divisor cut out on the triangle $xyz = 0$ in $\mathbf{P}^2$ by a general quartic curve. As in [III, Example **], they impose only eleven independent linear conditions to $|\mathcal{O}_{X_\Delta}(1)|$. The projection of $X_\Delta \subseteq \mathbf{P}^{14}$ from these twelve points is a quartic surface $X_F \subseteq \mathbf{P}^3$, with a triple point as a result of the contraction of $\text{ꓘ}_\Delta$. It is one of the relevant degenerations of $K3$ surfaces alternative to the tetrahedron $X_0$.

Correspondingly, one wants to enumerate hyperplane sections of this quartic surface $X_F \subseteq \mathbf{P}^3$ of genera $2, 1, 0$, respectively, and passing through $2, 1, 0$ general points of $X_F$, respectively. Again, to set this up as in (2.2), we impose the passing through only eleven out of the twelve base points. For instance, we set $\alpha = (4, 4, 3)$ and $\beta = (0, 0, 1)$.

By degeneration, one finds

$$N_2^\Delta(\alpha, \beta) = 21; \quad N_1^\Delta(\alpha, \beta) = 132; \quad N_0^\Delta(\alpha, \beta) = 304,$$

cf. [10, Theorem 3]. Mikhalkin (private communication) has been able to find these numbers with floor diagrams, but we won't do it in these notes.

**(2.7)** Consider the surface $X_\Delta \subseteq \mathbf{P}^4$ corresponding to the polygon $\Delta$ represented below.

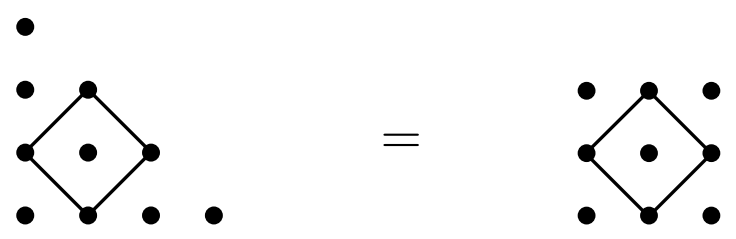

It is the quartic surface $xy - t^2 = zw - t^2 = 0$ in homogeneous coordinates $(x : y : z : w : t)$. It has four ordinary double points and no further singularity. Its anticanonical cycle $\text{ꓘ}_\Delta$ is a quadrilateral made of four lines, the vertices of which coincide with the ordinary double points of the surface.

The linear system $|\mathcal{O}_{X_\Delta}(1)|$ may be thought of either as a linear system of plane cubics with five base points — as indicated by the picture on the left —, or as a linear system of curves of bidegree $(2, 2)$ on $\mathbf{P}^1 \times \mathbf{P}^1$ with four base points — as indicated by the picture on the right. Alternatively, $X_\Delta$ may be realized as the double cover of $\mathbf{P}^1 \times \mathbf{P}^1$ branched over an anticanonical quadrilateral.

The quartic surface $X_\Delta \subseteq \mathbf{P}^4$ and the complete linear system of its hyperplane sections is the key to understand the limits, in terms of the theory expound in [I], of hyperplane sections of smooth surfaces degenerating to sections of some surface $X_0$ by a hyperplane passing through a point $v$, such that $X_0$ has four irreducible components through $v$, all 2-planes, and the face of the dual complex of $X_0$ corresponding to $v$ is a quadrilateral, e.g., $v$ is a vertex of an octahedron $X_0 \subseteq \mathbf{P}^5$. This is analogous to the cubic $X_\Delta \subseteq \mathbf{P}^3$ of (2.4) corresponding to sections of $X_0$ by a hyperplane passing through a triple point formed by three 2-planes, see [I, Section **].

In this context, the enumerative problem one needs to consider is that of counting the hyperplane sections of $X_\Delta$ of genus 0 intersecting $\text{ꓘ}_\Delta$ along a prescribed member of $|\mathcal{O}_{\text{ꓘ}_\Delta}(1)|$; such a divisor consists of four coplanar points, one on each side of $\text{ꓘ}_\Delta$. Correspondingly, one sets (for example) $\alpha = (1, 1, 1, 0)$ and $\beta = (0, 0, 0, 1)$. One finds $N_0^\Delta(\alpha, \beta) = 4$ in the same way as in (2.4). We will compute this number from the tropical perspective in (6.10), see also (6.12), and using floor diagrams in (7.17.1).

# 3 – From the torus to the tropics

In this section we present amœbæ and tropical hypersurfaces. The amœba of a hypersurface in the torus $(\mathbf{C}^*)^n$ has first been considered by Gelfand–Kapranov–Zelevinsky in [15]. Non-archimedean amœbæ have been introduced by Kapranov, who showed that they coincide with tropical hypersurfaces.

**(3.1) Amœba of a hypersurface in a torus.** Let $X$ be a hypersurface of the torus $(\mathbf{C}^*)^n$. The *amœba* of $X$ is the image of $X$ by the logarithmic moment map

$$\begin{aligned} \mathrm{Log}: \quad (\mathbf{C}^*)^n &\to \mathbf{R}^n \\ (z_1,\dots,z_n) &\mapsto \big(\log|z_1|,\dots,\log|z_n|\big). \end{aligned}$$

Pictured below is the amœba of a line in $(\mathbf{C}^*)^2$. From this you should guess whence the name 'amœba' came from.

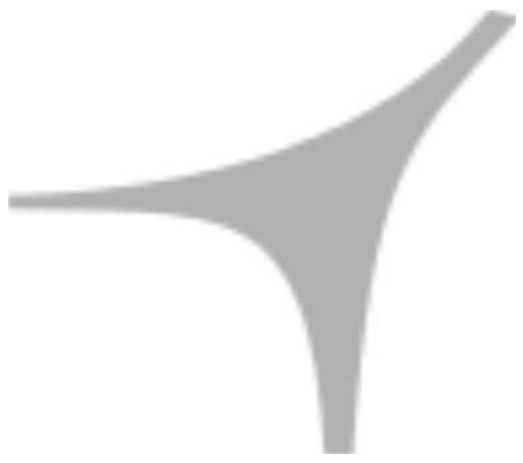

**(3.2) Limits of amœbæ.** One of the ideas at the origin of tropical geometry is to consider limits of amœbæ. To this effect, consider a hypersurface $X \subseteq (\mathbf{C}(t)^*)^n$, which may be seen as a rational family of hypersurfaces in $(\mathbf{C}^*)^n$. The $\mathbf{C}(t)$-points of $X$ are $n$-tuples $z(t)$ of rational functions $(z_1(t),\dots,z_n(t))$. Renormalizing, it is possible to consider the limit of $\mathrm{Log}(z(t))$ at $t=0$, namely

$$\lim_{t\to 0} \frac{1}{\log|t|}\big(\log|z_1(t)|,\dots,\log|z_n(t)|\big) = \big(v(z_1),\dots,v(z_n)\big),$$

where $v$ is the $t$-adic valuation, i.e., for $f(t)\in\mathbf{C}(t)$, $v(f(t))$ is the integer $a$ such that one may write $f(t) = t^a(\alpha + \sum_{n=1}^{+\infty}\alpha_n t^n)$ for some $\alpha \neq 0$.

The field $\mathbf{C}(t)$ however is not algebraically closed, so in general $X$ doesn't have enough $\mathbf{C}(t)$-points, and the previous considerations don't give a fully satisfactory definition of the limit at $t=0$. The straightforward workaround to this issue is to replace $\mathbf{C}(t)$ by its algebraic closure; another nice possibility is to work over the field of Puiseux series — the latter being the algebraic closure of the field of formal Laurent series $\mathbf{C}((t))$. Recall that a Puiseux series of the variable $t$ and with coefficients in $\mathbf{C}$ is a formal power series of the form

$$f(t) = \sum_{n=n_0}^{+\infty} \alpha_n t^{n/N}$$

for some $n_0 \in \mathbf{Z}$ and $N \in \mathbf{N}^*$.

The idea of taking limits of amœbæ is one inspiration for the definitions in the following paragraph.

**(3.3) Non-Archimedean amœbae.** Let $\mathbf{K}$ be a valuated, non-Archimedean, algebraically closed field. Being non-Archimedean means that $\mathbf{K}$ has some infinitely large elements, i.e., elements that are larger than all integers: in the context of series in $t$ observed at $t = 0$, these are the series with negative $t$-adic valuation.

Precisely, we say that $\mathbf{K}$ is *valuated non-Archimedean* if it comes with a function $v : \mathbf{K}^* \to \mathbf{R}$, the *valuation*, such that

$$\forall f, g \in \mathbf{K}^* : \quad v(fg) = v(f) + v(g), \quad \text{and} \quad v(f+g) \geqslant \min\{v(f), v(g)\} \text{ (if } f+g \neq 0).$$

This determines a *non-Archimedean norm* $|\cdot|$, by setting

$$|f| = \begin{cases} e^{-v(f)} & \text{if } f \in \mathbf{K}^* \\ 0 & \text{if } f = 0. \end{cases}$$

Now, let $X \subseteq (\mathbf{K}^*)^n$ be a hypersurface. The (non-Archimedean) *amœba* of $X$ is, by definition, its image by the valuative moment map

$$\begin{aligned} \mathrm{Val} : \quad (\mathbf{K}^*)^n &\longrightarrow \mathbf{R}^n \\ (x_1, \ldots, x_n) &\longmapsto \big(v(x_1), \ldots, v(x_n)\big). \end{aligned}$$

Pictured below is the amœba of a line in $(\mathbf{K}^*)^2$. One way to think about it is as the limit of the amœbæ of a rational family of lines in $(\mathbf{R}^*)^2$, i.e., the limit as $t$ goes to 0 of the amœbæ of a family of lines defined by an equation $a(t)x + b(t)y + c(t) = 0$, where $a, b, c$ are rational functions in the parameter $t$.

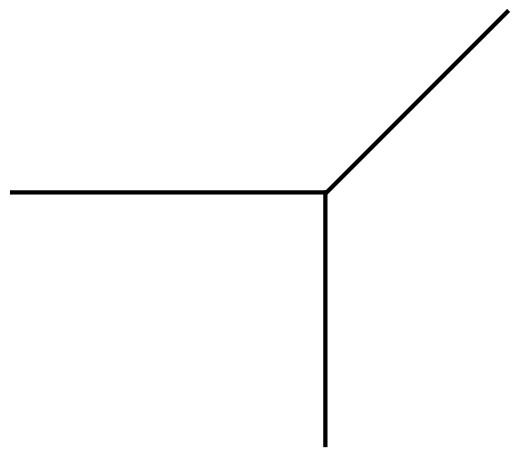

**(3.4) Theorem** (Kapranov[3], see [12]). *Let $\mathbf{K}$ be a valuated, non-Archimedean, algebraically closed field, and assume its valuation $v : \mathbf{K}^* \to \mathbf{R}$ is surjective.*

*Consider a hypersurface $X \subseteq (\mathbf{K}^*)^n$ defined by a polynomial*

$$P(x) = \sum_{I \in \mathrm{Supp}(P)} a_I x^I$$

*(*$\mathrm{Supp}(P)$ *is the support of $P$, i.e., the set of those multi-indices $I$ such that that coefficient of $P$ in the monomial $x^I$ is non-zero).*

*The amœba of $X$ is the corner locus of the piecewise linear function*

$$p \in \mathbf{R}^n \longmapsto \max_{I \in \mathrm{Supp}(P)} \{I \cdot p - v(a_I)\}.$$

The corner locus of a piecewise linear function $f : \mathbf{R}^n \to \mathbf{R}$ is the set of walls in $\mathbf{R}^n$ separating the chambers on which $f$ equals a unique linear function. We will come back to this notion in the next section.

The theorem also holds in case the valuation is not surjective, provided one takes the closure of the amœba. This typically happens when $\mathbf{K}$ is the field of Puiseux series.

[3] this theorem has apparently appeared in a preprint of Kparanov only, in 2000; this preprint is referenced in [21] for instance.

**(3.5) The tropical semi-field.** Still applying the idea of taking the limit at $t = 0$, one defines a structure of semi-field on $\mathbf{R}$ as the limit of the usual field structure.

We consider the family indexed by $t \in \mathbf{R}_+^*$ of bijective maps

$$\frac{\log}{\log t} : \mathbf{R}_+^* \to \mathbf{R},$$

and use it to degenerate the field structure on $\mathbf{R}$: One defines the two operation laws $\boxplus$ and $\boxdot$ as the limits of $+$ and $\cdot$ respectively, i.e., for all $x, y \in \mathbf{R}$

$$\begin{aligned} x \boxplus y &\overset{\text{def.}}{=} \lim_{t\to 0^+} \frac{1}{\log t} \log\big((\log t))e^x + (\log t)e^y\big) = \max\{x, y\} \\ \text{and} \quad x \boxdot y &\overset{\text{def.}}{=} \lim_{t\to 0^+} \frac{1}{\log t} \log\big((\log t))e^x \cdot (\log t)e^y\big) = x + y. \end{aligned}$$

We call $(\mathbf{R}, \boxplus, \boxdot)$ the *tropical semi-field*, and denote it by $\mathbf{T}$. A convenient notation for the operations $\boxplus$ and $\boxdot$ is to surround algebraic expressions by quotation marks to indicate that they are meant in the tropical semi-field. Thus, "$x + y$" $= x \boxplus y$ and "$xy$" $= x \boxdot y$.

**(3.6) Tropical polynomials and hypersurfaces.** The definition of a tropical polynomial (function) is pretty straightforward, one just replace the ordinary algebraic expression by its counterpart in the tropical semifield. Thus, a *tropical polynomial* is a function $P : \mathbf{T}^n \to \mathbf{T}$ defined by an expression of the form

$$\begin{aligned} P(x_1, \ldots, x_n) &= \text{``} \sum_I a_I x^I \text{''} = \boxplus_I (a_I \boxdot x^{\boxdot I}) \\ &= \max_I \{a_I + i_1 \cdot x_1 + \cdots i_n \cdot x_n\}, \end{aligned}$$

where the sum is taken on a finite set of multi-indices.

It does not make much sense to consider the zero locus of a function as above. Instead, a *tropical hypersurface* in $\mathbf{T}^n$ is by definition the corner locus of a tropical polynomial.

With these definitions, the content of Theorem (3.4) is that amœbæ of hypersurfaces in $(\mathbf{K}^*)^n$ and tropical hypersurfaces in $\mathbf{T}^n$ ($\mathbf{R}^n$, if one prefers) are the same objects.

# 4 – Parabolic polarity[4], b.k.a. Legendre–Fenchel duality

This section is about what is often called Legendre transformation. A nice reference is [17, Part E]. For the case of twice differentiable functions, a classical reference is [1, Section 14].

**(4.1) Definition.** *Let $f : \mathbf{R}^n \to \mathbf{R}$ be a function. The* Legendre–Fenchel transform *of $f$ is the function*

$$\begin{aligned} f^\vee : \quad \check{\mathbf{R}}^n &\to \mathbf{R} \\ p &\mapsto \sup_{x\in\mathbf{R}^n} \big\{p \cdot x - f(x)\big\}, \end{aligned}$$

*where $\check{\mathbf{R}}^n$ is the dual of $\mathbf{R}^n$.*

[4]about this name, imagined (but not used) by Fenchel, see [17, p. 247].

The definition is intended for all functions defined on a subset $U$ of $\mathbf{R}^n$, so in fact the supremum in the definition is to be taken on all $x \in U$ only.

In fact we don't allow $f^\vee$ to take the value $+\infty$, and thus $f^\vee$ is in general only defined on a subset of $\check{\mathbf{R}}^n$. To make sure that $f^\vee$ is defined on a non-empty subset of $\check{\mathbf{R}}^n$, one may assume that there exists an affine function $m$ such that $m \leqslant f$. We make this assumption from now on and for the rest of the section; it is true in particular if $f$ is convex.

The Legendre–Fenchel transform is always a closed convex function, see [17, Example B.2.1.3]. The Legendre–Fenchel transformation is involutive on closed convex functions. More precisely, for all functions $f$, the "bi-conjugate" $(f^\vee)^\vee$ is the closed convex hull of $f$, see [17, Theorem E.1.3.5]. By definition, the closed convex hull of $f$ is obtained by considering the epigraph $\mathrm{epi}(f)$ (i.e., the part of $\mathbf{R}^n \times \mathbf{R}$ lying above the graph of $f$), taking the closed convex hull $\overline{\mathrm{co}}\, \mathrm{epi}(f)$ of the latter, and finally taking the unique function of which this is the epigraph; in other words, it is the function

$$x \mapsto \inf\big\{y : (x, y) \in \overline{\mathrm{co}}\, \mathrm{epi}(f)\big\}.$$

**(4.2) Link with projective duality.** Consider a strictly convex function $f : \mathbf{R}^n \to \mathbf{R}$ that is twice differentiable (thus, $f'' > 0$). Then, for all linear functional $p \in \check{\mathbf{R}}^n$, the supremum of all $p \cdot x - f(x)$ is attained at a unique point $x(p) \in \mathbf{R}^n$, which may be characterized as the unique $x \in \mathbf{R}^n$ such that the tangent space to the graph of $f$ at the point $(x, f(x))$ has slope $p$; moreover, the equation of this very tangent space is

$$y = p \cdot x - f^\vee(p).$$

Therefore, the points of the graph of the Legendre–Fenchel transform $f^\vee$ naturally correspond to the tangent spaces of the graph of $f$: this is an instance of projective duality! (see [XII]).

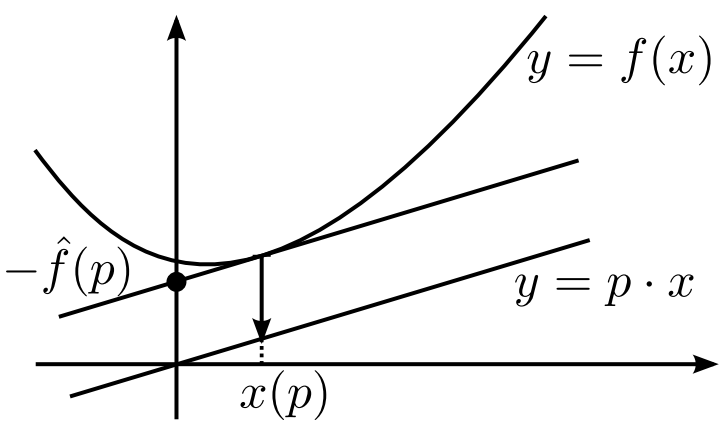

The Legendre–Fenchel transformation is also linked with classical mechanics, hence the notation '$p$' for linear functionals. Luckily enough, the letter 'p' can also stand for *pente* and *pendenza* (French and Italian for slope).

**(4.3)** For convex piecewise affine functions, there is an interpretation akin to that of (4.2). In this context, the Legendre–Fenchel transformation exchanges, loosely speaking, the roles of the corner loci and the slopes; this is perhaps best understood by looking at the figure below, which pictures the graphs of a piecewise affine linear function (defined on a finite interval) and its Legendre–Fenchel transform (defined on the whole $\mathbf{R}$). It is elementary to check that things indeed happen as on the picture, and that taking the Legendre–Fenchel transform of the function on the right gives back the function on the left; I of course encourage the reader to try it for him or herself.

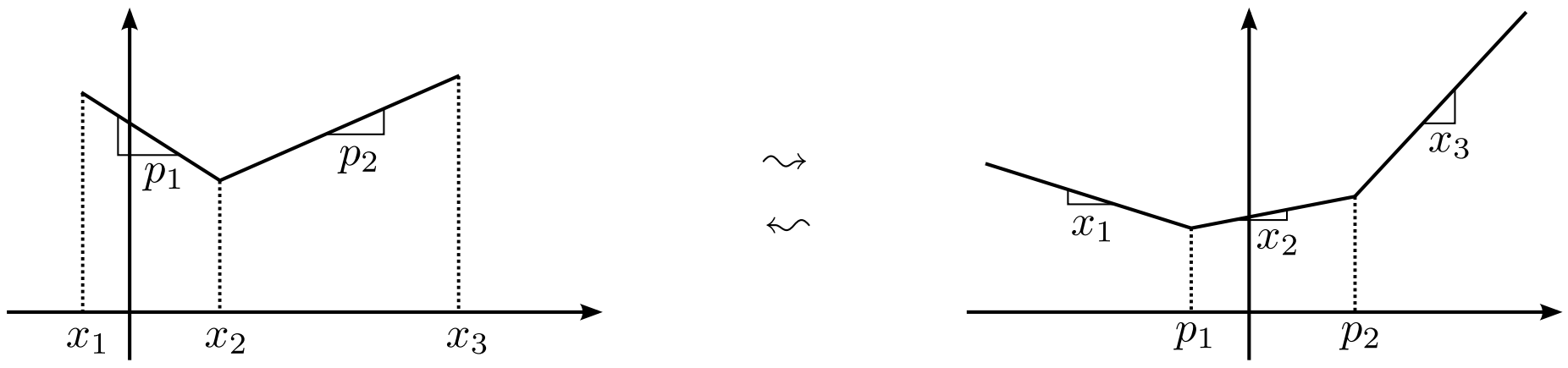

There is a similar story for functions $f$ defined at finitely many points $x_1, \dots, x_n$ only. In this case, the Legendre–Fenchel transform is a convex piecewise affine function defined on the whole dual space, and the bitransform is a convex piecewise affine function defined on the convex hull of $x_1, \dots, x_n$, and characterized by the fact that its graph is the bottom of the convex hull of the graph of $f$. Consider for instance the restriction to $\{x_1, x_2, x_3\}$ of the function on the left in the above figure. Its transform is the function on the right of the picture, and its bi-transform is the function on the left itself. This is again elementary to check.

**(4.4) Application to tropical geometry.** In the context of Theorem (3.4), the function $p \in \mathbf{R}^n \mapsto \max_{I \in \operatorname{Supp}(P)}\{I \cdot p - v(a_I)\}$ is the Legendre–Fenchel transform of the function

$$v_P : I \in \operatorname{Supp}(P) \mapsto v(a_I). \tag{4.4.1}$$

Therefore the above discussion shows that, dual to the amœba of $X$, which may be seen as a subdivision of $\mathbf{R}^n$ into convex polytopes, there is a subdivision of the Newton polytope $\Delta(P)$ of $P$ (i.e., the convex hull of $\operatorname{Supp}(P)$) into convex polytopes, given by the corner locus of the real-valued function on $\Delta(P)$ obtained by "extending by convexity" the function $v_P$, i.e., of the closed convex hull $\bar{v}_P$ of $v_P$. In other words, the tiles of this subdivision of $\Delta(P)$ are the maximal polytopes on which the function $\bar{v}_P$ is affine.

# 5 – Tropical curves

In this section we introduce plane tropical curves, following Mikhalkin, first as rational piecewise linear objects in $\mathbf{R}^2$, which we will see coincide with the tropical hypersurfaces that have been presented in Section 3, and then as images of abstract tropical curves by tropical maps. Beware not to confuse the graphs that we use here in order to support our definitions of tropical curves with those used in Section 7 to define floor diagrams; floor diagrams are purely combinatorial objects retaining only some information of plane tropical curves, the latter being geometric objects in nature.

## 5.1 – Plane tropical curves

**(5.1) Definition.** *A* plane tropical curve *is a rational, piecewise linear, real curve $C \subseteq \mathbf{R}^2$, which is a balanced topological incarnation of a finite graph $\Gamma$ equipped with a weight function $w : \operatorname{Edges}(\Gamma) \to \mathbf{N}^*$, by which we mean that conditions* (i)–(iv) *below hold.*

(i) The curve $C$ is the union of finitely many linear pieces (i.e., segments, semi-lines, or (rarely) lines) with rational slopes, and these linear pieces are in bijective correspondence with the edges of $\Gamma$. Each such linear piece has two ends, possibly at infinity (one end at infinity for a semi-line, two for a line), and the ends of the linear pieces of $C$ are in bijective correspondence with the vertices of $\Gamma$; moreover, the ends at infinity are in bijective correspondence with the 1-valent vertices of $\Gamma$. Beware that points at infinity are not intended in the usual sense of projective geometry (rather, see Remark (6.3)): lines have two ends at infinity, and disjoint parallel lines do not meet at infinity.

We may speak about the edges and vertices of a tropical curve to denote their incarnations in $\mathbf{R}^2$ or as points at infinity. We call vertices at finite distance (resp., at infinity) of $C$ the vertices of $\Gamma$ that are actually incarnated as points in $\mathbf{R}^2$ (resp., incarnated as points at infinity). Equivalently, the vertices at infinity of $C$ are the 1-valent vertices of $\Gamma$, and the vertices at finite distance are the others.

(ii) Any two linear pieces $C_1$ and $C_2$ of $C$ intersect properly and at their ends, i.e., if non-empty, the intersection $C_1 \cap C_2$ is a single point, which must be an end of both $C_1$ and $C_2$.

(iii) Any two linear pieces $C_1$ and $C_2$ of of $C$ with the same slope may intersect only if there is at least a third edge ending at $C_1 \cap C_2$.

For each vertex $v$ and edge $\Sigma$ ending at $v$ (we write $v \in \Sigma$ for this condition), we denote by $u_{v,\Sigma}$ the primitive (see (0.1)) integral vector $(a, b) \in \mathbf{Z}^2$ that has the same direction as the linear piece of $C$ corresponding to $\Sigma$, and is oriented from $v$ to the other end of $C$. We also let $u_{\Sigma,v}$ be $-u_{v,\Sigma}$.

(iv) *(Balancing condition)* For each plurivalent vertex $v$, one has

$$\sum_{\Sigma \ni v} w(\Sigma) \cdot u_{v,\Sigma} = 0. \tag{5.1.1}$$

Note that conditions (iii) and (iv) imply that the vertices of a plane tropical curve have valence either 1, or at least 3.

**(5.2) Examples.** Examples of plane tropical curves are pictured on the left in Figures 3 and 4. The weight of an edge is indicated as a circled integer, only if it is larger than one. The grids on the pictures are drawn so that one can easily read the rational slopes, but they do not mean that the vertices are integral points, which indeed is not requested in the definition of plane tropical curves.

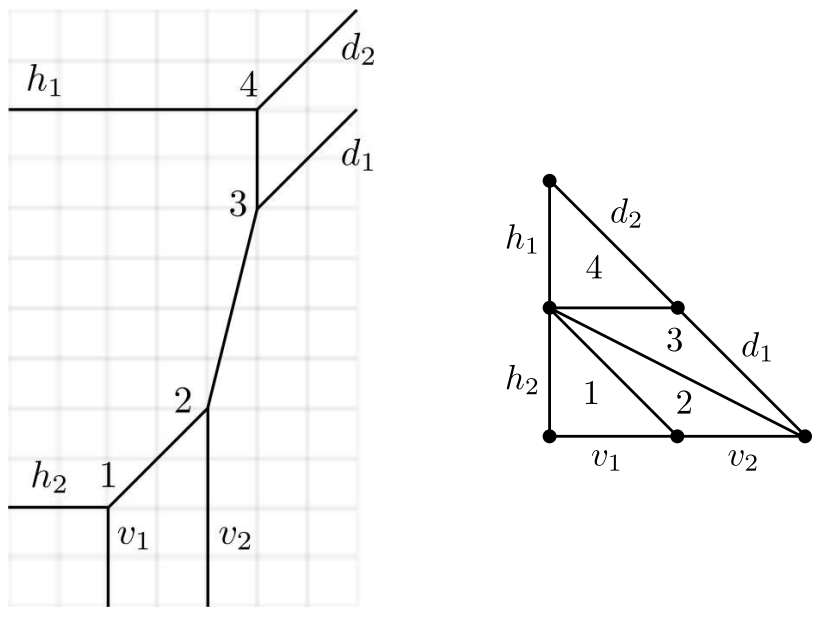

Figure 3

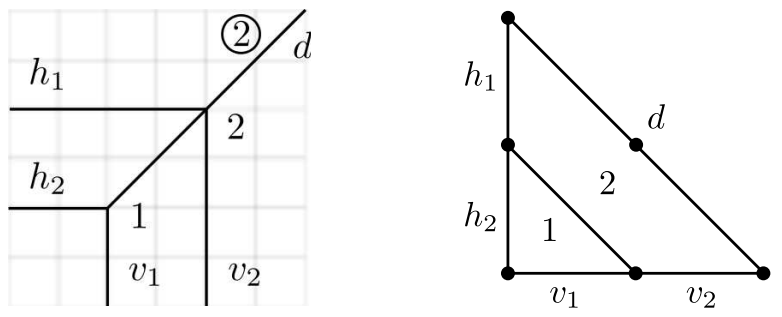

Figure 4

**(5.3) Theorem** (Mikhalkin [22, Section 2.1]). *Plane tropical curves as defined in (5.1) are exactly the hypersurfaces of the tropical plane $\mathbf{T}^2$, i.e., the corner loci of tropical polynomials, as defined in (3.6).*

The main content of this theorem is that the balancing condition (5.1.1) characterizes tropical hypersurfaces. To make this statement work, one needs to equip the corner locus of a tropical polynomial with a weighting. This has been introduced by Mikhalkin, and goes as follows. Let $\Sigma$ be an edge of the corner locus of some tropical polynomial. Then $\Sigma$ is a corner between the graphs of two affine functions $a+ix+jy$ and $b+kx+ly$ with rational slopes. By definition, the *weight* of $\Sigma$ is the integral distance between the two integral points $(i,j)$ and $(k,l)$ in $\check{\mathbf{R}}^2$, i.e., the number of integral points on the segment $\big[(i,j),(k,l)\big]$ minus one, which we have denoted by $\ell_{\mathbf{Z}}\big[(i,j),(k,l)\big]$ in Section 1.

**(5.4) Newton polygon of a plane tropical curve.** Let $C\subseteq\mathbf{R}^2$ be a tropical curve, and consider the list $(\Sigma_1,\dots,\Sigma_n)$ of its infinite edges, as a list ordered up to cyclic permutation, as follows: (a) an edge is *infinite* if it has infinite Euclidean length in $\mathbf{R}^2$, equivalently if it has at least one end that is a vertex at infinity; (b) the order on the set of infinite edges is defined by drawing a circle of sufficiently large radius, so that it intersects $C$ at exactly one point on each infinite edge and nowhere else; then the order is that of the intersection points on the circle according to the orientation of $\mathbf{R}^2$.

Starting at an arbitrary integral point of $\mathbf{R}^2$, draw successively the vectors $w(\Sigma_1).u^{\perp}_{\Sigma_1,\infty},\dots,$ $w(\Sigma_n).u^{\perp}_{\Sigma_n,\infty}$, so that when one vector is drawn, the next one starts where the previous one ended; for all infinite edges $\Sigma$, the vector $u_{\Sigma,\infty}$ is the unitary direction vector of $\Sigma$ oriented towards the end at infinity. It is a consequence of the balancing condition that this circuit of vectors is a closed path, and actually a convex integral polygon. We call it the *Newton polygon* of $C$.

*Examples.* For each example of plane tropical curve in (5.2), the corresponding Newton polygon is indicated next to it, on the right. The subdivision corresponds to the following construction.

**(5.5) Subdivision of the Newton polygon.** Going on with the same idea as in (5.4), one may enhance the Newton polygon of a plane tropical curve with a subdivision. This subdivision is dual to the subdivision of $\mathbf{R}^2$ defined by the tropical curve, as we have seen in (4.4) using the Legendre–Fenchel transformation.

Let $C\subseteq\mathbf{R}^2$ be a tropical curve, and $v$ one of its vertices. We consider the list $(\Sigma_1,\dots,\Sigma_n)$, ordered up to cyclic permutation, of the edges ending at $v$. Because of the balancing condition at $v$, drawing successively the vectors $w(\Sigma_1)\cdot u^{\perp}_{v,\Sigma_1},\dots,w(\Sigma_n)\cdot u^{\perp}_{v,\Sigma_n}$ builds an integral convex polygon $\Delta_v$.

There is a canonical way of arranging all the $\Delta_v$ together as a polygonal decomposition of the Newton polygon $\Delta$ of $C$. Indeed, each edge $\Sigma$ of $C$ has two ends $v_1$ and $v_2$, hence determines a glueing of $\Delta_{v_1}$ and $\Delta_{v_2}$ along the segments of their boundaries supported along $w(\Sigma)u_{v_1,\Sigma}$ and $w(\Sigma)u_{v_2,\Sigma}$, respectively (if $v$ is a vertex at infinity, then $\Delta_v$ is the Newton polygon $\Delta$). This is the subdivision we are interested in.

*Examples.* For each plane tropical curve in (5.2), the Newton polygon is indicated with its subdivision pictured.

**(5.6) Subdivisions and degenerations.** To a subdivision of a lattice polygon $\Delta$ into a tiling of lattice polygons as above, one may associate a “reducible toric surface” with the same degree and embedding dimension as $X_\Delta$, see below. In fact, such a reducible toric variety (in any dimension) may be realized as a toric degeneration of $X_\Delta$, see, e.g., [24, Section 3]. It is then natural (in hindsight!) to expect that plane tropical curves with the same subdivised polygon $\Delta$ are related to this degeneration. This has first been realized by Viro [27], who establishes this relation by means of his patchworking construction. Pursuing this idea, Shustin explains in [26] (see also [18, Chapter 2]) that the two processes of tropicalization and patchworking are somehow inverse to each other (one can think of them as degeneration and deformation, respectively), and as an application he gives an alternative proof of Mikhalkin’s Correspondence

Theorem (6.8). The relation of plane tropical curves with toric degenerations of the projective plane is also illustrated in Paragraph (8.4), following [4].

Now, the reducible toric variety defined by a subdivision of a lattice polytope $\Delta$ is the following: its irreducible components are the various $X_{\Delta_v}$ corresponding to the various tiles $\Delta_v$ of the subdivision of $\Delta$, and these are glued along toric varieties of smaller dimensions, according to the dual graph of the subdivision of $\Delta$. Let me simply give an example instead of making a formal definition. Consider the subdivised polygon $\Delta$ from Figure 4. The polygon $\Delta$ is the triangle $T_2$ of (1.1), which defines the 2-Veronese surface $X_\Delta \subseteq \mathbf{P}^5$, of degree 4. The trapezium in the subdivision is a $T^1_{1,2}$ in the notation of (1.6), and thus defines a degree 3 rational normal scroll $X_1 \subseteq \mathbf{P}^4$, built on a line $X_{1,2}$ and a conic spanning two supplementary linear spaces $\mathbf{P}^1$ and $\mathbf{P}^2$ in $\mathbf{P}^4$, such that the line $X_{1,2}$ corresponds to the intersection of the two tiles of the subdivision. We may view the $\mathbf{P}^4$ spanned by $X_1$ as a linear subspace in $\mathbf{P}^5$. The triangle in the subdivision defines a plane $X_2 = \mathbf{P}^2$. We may view this $\mathbf{P}^2$ as a linear subspace of $\mathbf{P}^5$, intersecting the $\mathbf{P}^4$ spanned by $X_1$ exactly along the $\mathbf{P}^1$ spanned by the line $X_{1,2}$. Then, the reducible toric surface corresponding to the subdivision is the union $X_1 \cup X_2 \subseteq \mathbf{P}^5$, where $X_1 \cap X_2 = X_{1,2}$.

## 5.2 – Abstract tropical curves and tropical maps

The definitions in this section are due to Mikhalkin [23].

**(5.7) Definition.** *An* abstract tropical curve *is the data of a finite graph $\bar\Gamma$ without isolated vertex, together with a complete metric $g$ on $\Gamma \overset{def.}{=} \bar\Gamma \setminus \Gamma^\infty$, where $\Gamma^\infty$ is the set of vertices of $\Gamma$ having valence* 1.

The condition on the completeness of the metric means that every edge of $\Gamma$ should be isometric either to a segment, a semi-line, or a line, depending on how many ends at infinity it has.

A tropical curve $(\bar\Gamma, g)$ is said to be *irreducible* if the space $\bar\Gamma$ is connected. The *arithmetic genus* $p_a(\bar\Gamma, g)$ of $(\bar\Gamma, g)$ is defined by the first Betti number of $\bar\Gamma$ (and/or its topological Euler–Poincaré characteristic) as

$$p_a(\bar\Gamma, g) \overset{\text{def.}}{=} 1 - \chi_{\text{top}}(\bar\Gamma) = b_1(\bar\Gamma) - n + 1,$$

where $n$ is the number of connected components of $\bar\Gamma$.

**(5.8) Definition.** *A* tropical morphism *from an abstract tropical curve $(\bar\Gamma, g)$ to the tropical plane is a continuous map $f : \Gamma \to \mathbf{R}^2$ such that: (i) the restriction of $f$ to each edge is affine integral; and (ii) the balancing condition holds at every vertex at finite distance.*

Some more precision is in order. Let $\Sigma$ be an edge of $\Gamma$. The restriction of $f$ to $\Sigma$ is *affine integral* if (i) the image $f(\Sigma) \subseteq \mathbf{R}^2$ is supported on a line directed by an integral unitary vector $u_\Sigma$, and (ii) the differential $df_\Sigma$ is $\mathbf{Z}$-linear with respect to the integral structure on the tangent spaces defined by $g$ and $u_\Sigma$ respectively; the absolute value of the $\mathbf{Z}$-linearity factor of $df_\Sigma$ is called the *stretching factor*, and denoted by $w_{\Sigma,f}$.

Let $v$ be a plurivalent vertex of $\Gamma$. The balancing condition at $v$ is

$$\sum_{\Sigma \ni v} w_{\Sigma,f} \cdot u_{v,\Sigma} = 0,$$

where as before $u_{v,\Sigma}$ is the unitary direction vector for $f(\Sigma)$ directed from $v$ to the other end.

**(5.9) Example.** Two abstract tropical curves are pictured on Figure 5, of genera 1 and 0, respectively. They both admit a tropical morphism to the plane tropical curve pictured on Figure 6.

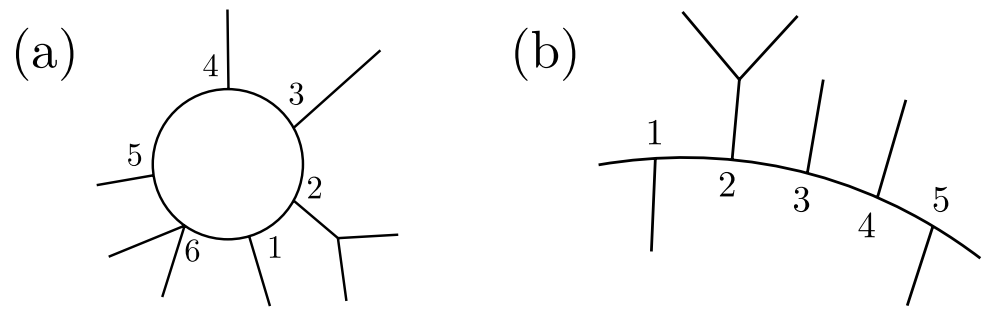

Figure 5

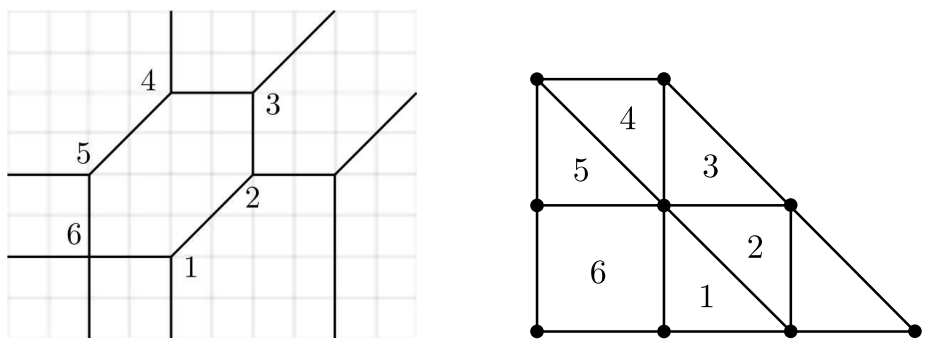

Figure 6

**(5.10)** If $C$ is an abstract tropical curve with a tropical morphism $f$ to the plane, one defines the *Newton polygon* of the map $f$, denoted by $\Delta(C, f)$, as the Newton polygon of the plane tropical curve image of $C$ by $f$; as such, it comes with a subdivision.

Actually, it is not completely straightforward how to define the image $f(C)$ as a plane tropical curve, as it may involve the adding of new vertices (at finite distance). This happens for instance in Example (5.9): the image curve must have the vertex #6 since, by definition, two edges may only intersect at a vertex, see (5.1)); thus in case (b) the graph has to be modified.

**(5.11) Geometric genus and resolutions.** So far we have defined the arithmetic genus of an abstract tropical curve; one may next want to define the geometric genus of a plane tropical curve as the arithmetic genus of its normalization, as in algebraic geometry. However, it is not clear how to appropriately define the normalization. One may instead define the *geometric genus* of a plane tropical curve $C \subseteq \mathbf{R}^2$ as the smallest possible arithmetic genus of all abstract tropical curves $\Gamma$ such that there exists a tropical map $f : \Gamma \to \mathbf{R}^2$ mapping $\Gamma$ to $C$. In general, it is arguably best to proceed as in the theory of stable maps and, instead of considering the geometric genus of a plane tropical curve itself, to consider the genus of a map from an abstract tropical curve to it, defined as the genus of the source curve.

**(5.12) Geometric genus and singularities.** Alternatively, there is in principle the possibility of defining the geometric genus of a plane curve $C$ as its arithmetic genus (defined as the number of integral points in the interior of its Newton polygon, as in (1.9)), minus some corrections determined by its singularities (the latter concept yet to be defined). It is however delicate to give a general definition along these lines, and we shall therefore only consider plane curves of a particular shape. It will follow from Proposition (6.4) that this limited range of shapes is in general sufficient in view of the tropical enumerative problem we will be interested in.

Let $C$ be a plane tropical curve, and assume that it is reduced[5], and that the tiles of the subdivision of its Newton polygon $\Delta$ are all triangles and parallelograms. The latter condition means that all vertices at finite distance of $C$ are either trivalent or look like the crossing of two lines (like vertex 6 on Figure 6). Then, define the *$\delta$-invariant* of $C$ as

$$\delta(C) = \sum_{e \text{ finite segment of } C} \big(w(e) - 1\big) + \sum_{\substack{v \text{ vertex of } C: \\ \Delta_v \text{ parallelogram}}} \mathscr{A}(\Delta_v) + \sum_{\substack{v \text{ vertex of } C: \\ \Delta_v \text{ triangle}}} p_a(\Delta_v)$$

where $w$ is the weight, $\mathscr{A}$ is the Euclidean area, and $p_a(\Delta_v)$ is the arithmetic genus of the lattice polygon $\Delta_v$, defined as in (1.9) as the number of integral points in its interior. A *finite segment* of $C$ is a maximal connected linear subset of $C$ having its two ends at finite distance; this corresponds to a maximal chain of parallelograms (in the sense of footnote 6 below) in the subdivision of $\Delta$ such that none of the extremal edges in the corresponding sequence of parallel edges lies on the boundary of $\Delta$. Note that a finite segment of $C$ is in general a chain of finitely

[5] a plane tropical curve is *reduced* if it cannot be written as the sum of two non-empty plane tropical curves $C_1$ and $C_2$ such that the support of $C_2$ is contained in the support of $C_1$.

many edges, all supported on the same line, and all having the same weight; the weight of the finite segment is this common weight.

The above definition of $\delta(C)$ is equivalent to the following. Let $v$ be a vertex of $C$, and let $\Delta_v$ be the corresponding lattice polygon in the subdivision of $\Delta$. An edge $\sigma$ of $\Delta_v$ is said to *border* $\Delta$ if it is linked to an edge in the boundary of $\Delta$ by a (possibly empty) chain of parallelograms[6]. For example, for the subdivision on Figure 6, in addition to the edges directly lying on the boundary of $\Delta$, the left edge of triangle 1 and the bottom edge of triangle 5 border $\Delta$. Recall that $\ell_{\mathbf{Z}}(\sigma)$ denotes the integral length of $\sigma$, i.e., the number of integral points on $\sigma$ minus one. Then, we let

$$\delta(v) = \begin{cases} p_a(\Delta_v) + \frac{1}{2} \sum_{\substack{\sigma \text{ edge of } \Delta_v \\ \text{not bordering } \Delta}} \big(\ell_{\mathbf{Z}}(\sigma) - 1\big) & \text{if } \Delta_v \text{ triangle} \\ \mathscr{A}(\Delta_v) & \text{if } \Delta_v \text{ parallelogram,} \end{cases}$$

and $\delta(C)$ be the sum of $\delta(v)$ over all vertices at finite distance of $C$.

Finally, we define the *geometric genus* of $C$ as its arithmetic genus $p_a(\Delta)$ minus its $\delta$-invariant $\delta(C)$. This is coherent with our previous definition by Lemma (5.15) below.

**(5.13) Examples.** With the above definitions, the plane tropical curve in Example (5.9) has arithmetic genus 1 and geometric genus 0. There are tropical maps to this plane tropical curve from the two abstract tropical curves (a) and (b), having genus 1 and 0, respectively. Note that vertex #6, which gives the only non-trivial contribution to the $\delta$-invariant of the plane curve, does not appear in the abstract curve (b).

The plane tropical curves on Figure 7 and 8 have arithmetic genus one and geometric genus zero. The only non-trivial contribution to the $\delta$-invariant of the curve on Figure 7 (resp., Figure 8) is from its weight two edge (resp., from the vertex with associated polygon the triangle with an integral point in its interior).

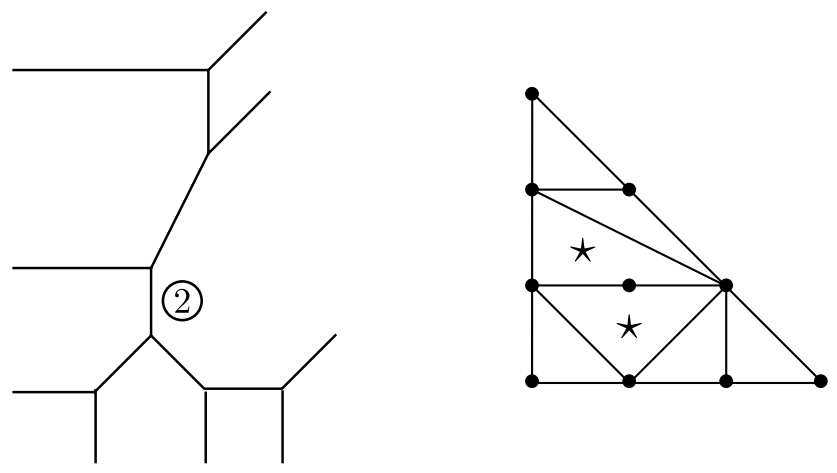

Figure 7

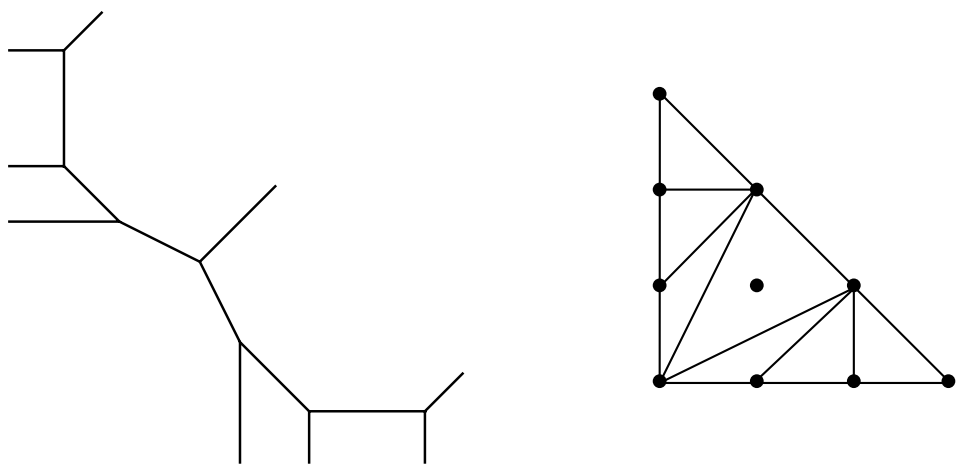

Figure 8

The plane tropical curve on Figure 9 has arithmetic genus three and geometric genus zero. The vertical weight two finite segment contributes by one to the $\delta$-invariant, and its crossing with another edge contributes by two.

**(5.14) Remark.** Remembering that tropical curves are in essence limits of families of curves, one way to look at the above definition of the $\delta$-invariant is in terms of the description of the limit of Severi varieties given in [I, Section $**$]. A plane tropical curve $C$ with Newton polygon $\Delta$ encodes data about limits $C_0$ of algebraic curves $C_t$ along a degeneration of the toric surface $X_\Delta$ into a union of the surfaces $X_{\Delta_v}$ over all vertices $v$ of $C$, as mentioned in Paragraph (5.6). For each vertex $v$, the curve $C_0 \cap X_{\Delta_v}$ may have up to $p_a(\Delta_v)$ (resp., $\mathcal{A}(\Delta_v)$) nodes off the toric

[6] two segments $\sigma$ and $\sigma'$ are *linked by a chain of parallelograms* if they are parallel and there exist parallelograms $\Delta_{v_1}, \ldots, \Delta_{v_n}$ in the subdivision of $\Delta$ such that $\sigma$ and $\sigma'$ are edges of $\Delta_{v_1}$ and $\Delta_{v_n}$ respectively and, for all $i = 1, \ldots, n-1$, $\Delta_{v_i}$ and $\Delta_{v_{i+1}}$ are adjacent along an edge parallel to $\sigma$ and $\sigma'$.

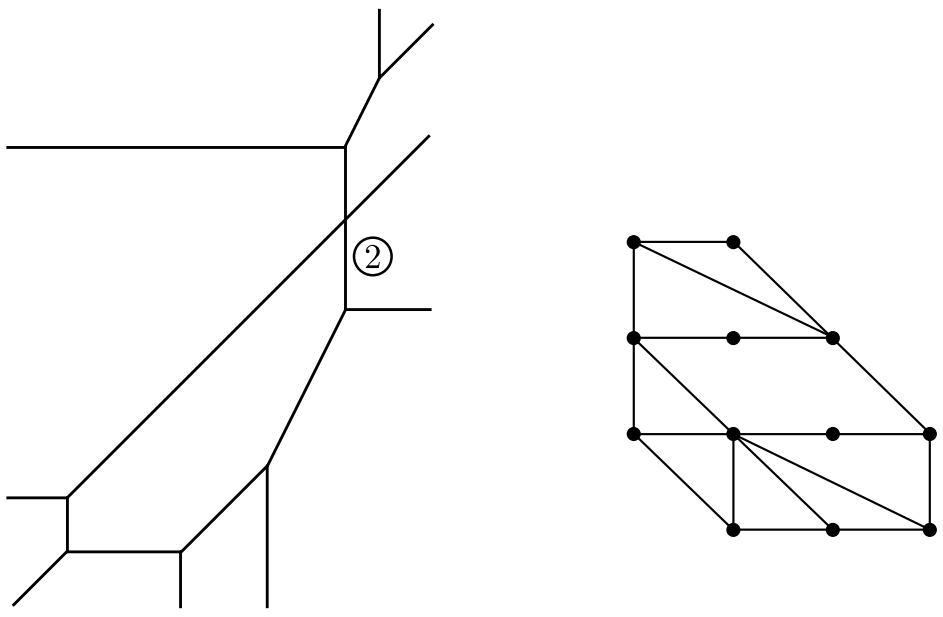

Figure 9

boundary of $X_{\Delta_v}$ if $\Delta_v$ is a triangle (resp., a parallelogram), and each of these nodes may be the limit of a node on the curves $C_t$. Moreover, for each finite segment $e$ of weight $w$ of $C$, the curve $C_0$ may have up to a $w$-tacnode along the intersection $X_{\Delta_v} \cap X_{\Delta_{v'}}$, where $v$ and $v'$ are the two ends of $e$ (if $e$ is made of several edges, this has to be taken with a grain of salt, see [IX, Paragraph (5.3)]); each such $w$-tacnode may be the limit of $w-1$ nodes on the curves $C_t$. In sum, the curve $C_0$ may be the limit of curves $C_t$ having up to $\delta(C)$ nodes.

**(5.15) Lemma.** *Let $C$ be a reduced plane tropical curve with Newton polygon $\Delta$ such that all the tiles of the subdivision of $\Delta$ are triangles and parallelograms. Then the two definitions of the geometric genus of a plane tropical curve in (5.11) and (5.12) coincide.*

The following folklore proof has been communicated to me by Erwan Brugallé.

*Proof.* For the sake of simplicity, we shall assume that $C$ is reduced and irreducible, i.e., it cannot be written as the sum of two non-empty plane tropical curves. In order to compute the geometric genus of $C$ as in (5.11), we note that abstract tropical plane curves of minimal genus with a map to $C$ are obtained by separating the "nodes" of $C$, i.e., by seeing each crossing of two line segments as the image of two disjoint line segments. Let $\Gamma$ be such an abstract tropical curve. It has only trivalent vertices besides the (univalent) vertices at infinity, in the number $t$ equal to the number of triangles in the subdivision of $\Delta$. The finite edges of $\Gamma$ (i.e., edges joining two vertices at finite distance) correspond to the finite segments of $C$, and their number is the number $a$ of unordered pairs of parallel edges $\{\sigma, \sigma'\}$ in the subdivision of $\Delta$ such that $\sigma$ and $\sigma'$ belong to two distinct triangles $\Delta_v$ and $\Delta_{v'}$ of the subdivision and are linked by a (possibly empty) chain of parallelograms (the most common situation is that $\sigma = \sigma'$ is the intersection of two adjacent triangles); note that an edge of a triangle $\Delta_v$ appears in such a pair if and only if it does not border $\Delta$, and in that case it appears in a unique pair. The topological Euler–Poincaré characteristic of $\Gamma$ is

$$\chi_{\mathrm{top}}(\Gamma) = 1 - g = t - a,$$

where $g$ is the genus of $\Gamma$, hence $g = 1 - t + a$. Note how this computation of $\chi_{\mathrm{top}}(\Gamma)$ uses the irreducibility of $C$: if $C$ contained a line joining two vertices at infinity (corresponding to two edges in the boundary of $\Delta$ linked by a chain of parallelograms), then $\Gamma$ would have a connected component not containing any vertex at finite distance, and the formula for $\chi_{\mathrm{top}}(\Gamma)$ would need to be suitably adjusted.

In order to relate this value with the genus of $C$ computed as in (5.12), we shall use Pick's formula (see (1.8) and (1.9)); it says that for all lattice polygons $D$,

$$p_a(D) = \mathscr{A}(D) - \frac{1}{2}\sum\nolimits_{\sigma \text{ edge of } D} \ell_{\mathbf{Z}}(\sigma) + 1.$$

Thus, for all vertices $v$ of $C$,

$$\delta(v) = \begin{cases} \mathscr{A}(\Delta_v) - \frac{1}{2}\sum_{\sigma \text{ edge of } \Delta_v} \ell_{\mathbf{Z}}(\sigma) + 1 + \frac{1}{2}\sum_{\substack{\sigma \text{ edge of } \Delta_v \\ \text{not bordering } \Delta}} \big(\ell_{\mathbf{Z}}(\sigma) - 1\big) & \text{if } \Delta_v \text{ triangle} \\ \mathscr{A}(\Delta_v) \qquad \text{if } \Delta_v \text{ parallelogram,} \end{cases}$$

and therefore

$$\begin{aligned} \delta(C) &= \mathscr{A}(\Delta) - \frac{1}{2} \sum_{\substack{(\sigma,\Delta_v), \\ \sigma \text{ edge of } \Delta_v, \\ \Delta_v \text{ triangle}}} \ell_{\mathbf{Z}}(\sigma) + t + \frac{1}{2} \sum_{\substack{(\sigma,\Delta_v), \\ \sigma \text{ edge of } \Delta_v, \\ \Delta_v \text{ triangle}, \\ \sigma \text{ not bordering } \Delta}} \ell_{\mathbf{Z}}(\sigma) - a \\ &= \mathscr{A}(\Delta) - \frac{1}{2} \sum_{\substack{(\sigma,\Delta_v), \\ \sigma \text{ edge of } \Delta_v, \\ \Delta_v \text{ triangle}, \\ \sigma \text{ bordering } \Delta}} \ell_{\mathbf{Z}}(\sigma) + t - a. \end{aligned}$$

Now we claim that the last sum equals the sum of the lengths of the edges of $\Delta$. Indeed, since $C$ is irreducible it cannot contain a line between two vertices at infinity. Therefore, an edge in the boundary of $\Delta$ cannot be linked by a chain of parallelograms to another edge in the boundary of $\Delta$, and thus every edge in the boundary of $\Delta$ is linked by a chain of parallelograms to an edge of a triangle of the subdivision of $\Delta$.

Finally, by Pick's formula again, $\delta(C) = p_a(\Delta) - 1 + t - a$, hence $p_a(\Delta) - \delta(C) = 1 - t + a$, which coincides with the genus $g$ computed in the first place. □

# 6 – The tropical enumerative problem and the correspondence theorem

In this section, we first state the tropical enumerative problem and describe its correspondence with the toric enumerative problem. Then we give several detailed examples.

## 6.1 – General theory

**(6.1) The tropical enumerative problem.** Let $\Delta \subseteq \mathbf{R}^2$ be a convex integral polygon, and call $\Sigma_1, \ldots, \Sigma_n$ its edges. We consider the data of an integer $g$ such that $0 \leqslant g \leqslant \operatorname{Card}(\mathring{\Delta}_{\mathbf{Z}})$, and $\alpha = (\alpha_1, \ldots, \alpha_n), \beta = (\beta_1, \ldots, \beta_n) \in \underline{\mathbf{N}}^n$ such that

$$\forall s = 1, \ldots, n, \qquad I\alpha_s + I\beta_s = \ell_{\mathbf{Z}}(\Sigma_s),$$

and

$$I\alpha < \operatorname{Card}(\partial\Delta_{\mathbf{Z}}).$$

In a nutshell, the problem is to find the number of plane tropical curves of genus $g$ and Newton polygon $\Delta$ that (i) have assigned and unassigned infinite edges as prescribed by $\alpha$ and $\beta$, respectively, and (ii) pass through the appropriate number of fixed general points of $\mathbf{R}^2$.

To set this up precisely, fix a set $\underline{x} = \{x_1, \ldots, x_{g-1+|\beta|}\}$ of $g - 1 + |\beta|$ general points of $\mathbf{R}^2$, and a set $\Omega = (\{p_{s,i,j}\}_{1 \leqslant j \leqslant \alpha_{s,i}})_{1 \leqslant s \leqslant n, 1 \leqslant i}$ of $\alpha$ general real numbers. Let $u_1, \ldots, u_n \in \mathbf{R}^2$ be the unitary vectors such that $u_1^\perp, \ldots, u_n^\perp$ give the oriented directions of $\Sigma_1, \ldots, \Sigma_n$, respectively ($\Sigma_1, \ldots, \Sigma_n$ are ordered, up to cyclic permutation, and oriented according to the orientation of the $\mathbf{R}^2$ in which $\Delta$ dwells). For all $s = 1, \dots, n$ and $p \in \mathbf{R}$, we let $\mathsf{L}_{u_s,p}$ be the line in $\mathbf{R}^2$ that is directed by the vector $u_s$ and contains the point with coordinates $p \cdot u_s^\perp$.

We then define $\mathcal{E}_g^\Delta(\alpha, \beta)(\underline{x}, \Omega)$ as the set of tropical maps $f$ from an abstract tropical curve $C$ of arithmetic genus $g$ to the tropical plane, such that the Newton polygon of $f$ equals $\Delta$, and in addition: (i) for all $1 \leqslant s \leqslant n$ and $i \geqslant 1$, the plane tropical curve $f(C)$ has $\alpha_{s,i} + \beta_{s,i}$ infinite edges $E$ of weight $\geqslant i$ such that $u_{E,\infty} = u_s$, with $\alpha_{s,i}$ of them supported on the lines $\mathsf{L}_{u_s,p_{s,i,1}}, \dots, \mathsf{L}_{u_s,p_{s,i,\alpha_{s,i}}}$, respectively, and (ii) $f(C)$ contains all the points $x_1, \dots, x_{g-1+|\beta|}$. The problem is to find the cardinality of the set $\mathcal{E}_g^\Delta(\alpha, \beta)(\underline{x}, \Omega)$, counting its elements with appropriate multiplicities; it is expected to be finite, cf. Proposition (6.4).

**(6.2) Example.** Consider the same Newton polygon $\Delta$ as in (1.1), (b), with the edges labelled as in the picture below.

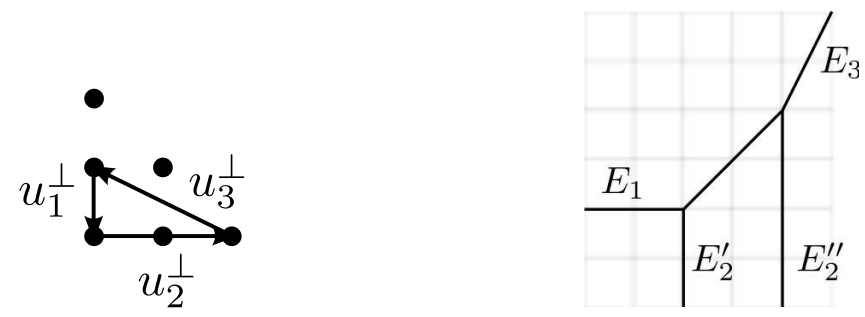

Next to the diagram is depicted a typical plane tropical curve with Newton polygon $\Delta$; it has arithmetic genus 0.

We consider the data $g = 0$, $\alpha = (0, 0, 1)$, $\beta = (1, [0, 1], 0)$, $\underline{x} = \{x_1\}$, and $\Omega = \{p_{3,1,1}\}$: this means that no condition is imposed on the infinite edge with direction $u_1$, and we require that (i) one edge with direction $u_2$ have weight $\geqslant 2$, with unassigned support, and (ii) one edge with direction $u_3$ be supported on the line $\mathsf{L}_{u_3,p_{3,1,1}}$ determined by the real number $p_{3,1,1}$; we require in addition that the curve contain the general point $x_1$.

One checks that there is a unique such tropical map, as indicated on the picture below (the fixed line $\mathsf{L}_{u_3,p_{3,1,1}}$ is represented as a dotted line).

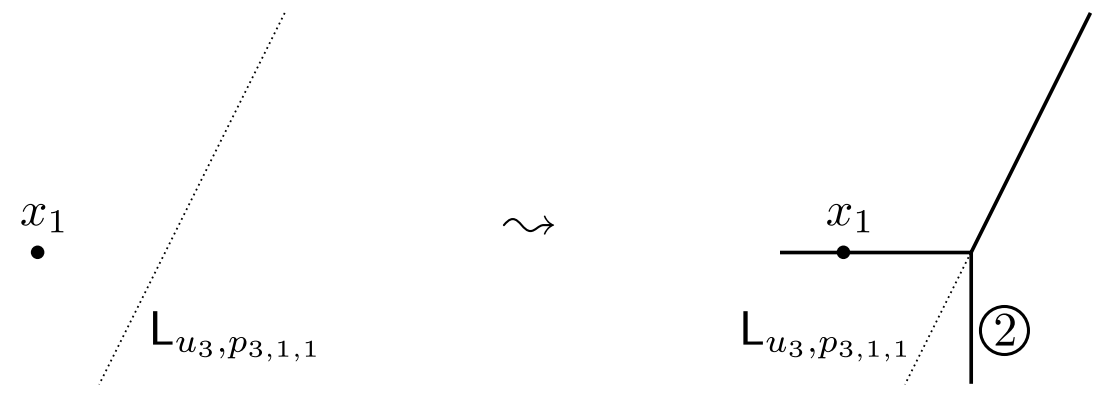

**(6.3) Remark.** The condition on infinite edges being supported by the assigned lines $\mathsf{L}$ may be formulated in terms of the passing through assigned points at infinity, for a suitable notion of points at infinity. The latter are the points one needs to add to $\mathbf{R}^2$ in order to obtain the tropical toric surface associated to $\Delta$, see [5, Section 7.3].

For instance the tropical projective plane, i.e., the tropical toric surface associated to a triangle $T_d$ as in Example (1.1), looks like the triangle $T_d$ itself, with the interior corresponding to $\mathbf{R}^2$ and the edges to the points at infinity; each edge may be identified with $\mathbf{R} \cup \{\pm\infty\}$. It is the quotient of $(\mathbf{R} \cup \{-\infty\})^3 \setminus \big\{(-\infty, -\infty, -\infty)\big\}$ by the equivalence relation $\sim$ defined by

$$(x_0, x_1, x_2) \sim (x_0', x_1', x_2') \quad \Longleftrightarrow \quad \exists \lambda \in \mathbf{R} : \ (x_0', x_1', x_2') = (x_0 + \lambda, x_1 + \lambda, x_2 + \lambda).$$

If $x_0 \neq -\infty$, then $(x_0, x_1, x_2)$ is equivalent to a unique point of the form $(0, x, y)$ with $(x, y) \in (\mathbf{R} \cup \{-\infty\})^2$, namely $(0, x_1 - x_0, x_2 - x_0)$; if $(x_0, x_1, x_2) \in \{-\infty\} \times \mathbf{R}^2$, then $(x_0, x_1, x_2)$ is equivalent to a unique point of the form $(-\infty, x, -x)$ with $x \in \mathbf{R}$, namely $(-\infty, \frac{x_1 - x_2}{2}, \frac{x_2 - x_1}{2})$.

The following result states that for a general choice of $\Omega$ and $\underline{x}$, the set $\mathcal{E}_g^{\Delta}(\alpha, \beta)(\underline{x}, \Omega)$ is indeed finite, and its members have the simplest possible shape. This is thus a tropical analogue of the main theorem in [III]. It is mainly due to Mikhalkin [23]; the article [23] does not comprise the possibility of boundary conditions (those encoded by $\alpha$ and $\beta$), even if it contains all the necessary ingredients to handle them. Boundary conditions in this context have first appeared in [14].

**(6.4) Proposition.** *Let $\Delta, g, \alpha, \beta$ be as in (6.1). For a general choice of $\underline{x}, \Omega$, the set $\mathcal{E}_g^{\Delta}(\alpha, \beta)(\underline{x}, \Omega)$ is finite. Moreover, for every $(C, f) \in \mathcal{E}_g^{\Delta}(\alpha, \beta)(\underline{x}, \Omega)$, the following hold:*

*(i) the graph $\Gamma = \bar{\Gamma} \setminus \Gamma^{\infty}$ underlying $C = (\bar{\Gamma}, g)$ is trivalent;*

*(ii) the map $f$ is a topological immersion;*

*(iii) the preimage of any point $x \in \mathbf{R}^2$ has at most two elements;*

*(iv) the infinite edges of $f(C)$ have minimal weight;*[7]

*(v) the set of vertices of the plane curve $f(C)$ is disjoint from the set $\underline{x}$;*

*(vi) an edge of the plane curve $f(C)$ contains at most one point from the set $\underline{x}$.*

Condition (ii) asserts in particular that no neighbourhood of a trivalent vertex of $\Gamma$ is ever mapped to a segment. Conditions (i) and (iii) together imply that the tiles of the subdivision of the Newton polygon of $f(C)$ are all triangles and parallelograms.

One should beware that Proposition (6.4) does not say that the shapes of the plane tropical curves corresponding to tropical maps in the set $\mathcal{E}_g^{\Delta}(\alpha, \beta)(\underline{x}, \Omega)$ are independent of the choice of $\underline{x}$ and $\Omega$, and this is indeed not true in general. For instance, their respective subdivisions of $\Delta$ may vary as $\underline{x}$ moves, cf. Example (7.14) and Figure 18.

**(6.5) Remark.** It most often happens that some members of $\mathcal{E}_g^{\Delta}(\alpha, \beta)(\underline{x}, \Omega)$ have their corresponding plane tropical curve that has finite edges of weight larger than one (see Example (7.14) for instance). This is analogous (see Remark (5.14)) to the fact that when considering degenerations of (classical) nodal algebraic curves, one more or less inevitably ends up with some tacnodal curves, which is one of the central themes of this whole volume.

In order to correspond to the counting of algebraic curves, the counting of tropical curves must be made with appropriate multiplicities. Besides, the mere cardinality (i.e., without multiplicities) of $\mathcal{E}_g^{\Delta}(\alpha, \beta)(\underline{x}, \Omega)$ does in general depend on the choice of $\underline{x}$ and $\Omega$.

**(6.6) Tropical multiplicity.** Let $f$ be a tropical morphism from an abstract tropical curve $C = (\bar{\Gamma}, g)$ to the tropical plane. Assume that every vertex $v \in \Gamma$ (i.e., every vertex at finite distance) is trivalent.

For each such vertex $v$, we define its multiplicity to be twice the Euclidean area of the triangle $\Delta_v$ dual to $v$ in the subdivision of the Newton polygon of $f$, i.e.,

$$\mu^{\mathbf{C}}(v) \stackrel{\text{def.}}{=} 2\mathscr{A}(\Delta_v) = \left|\det\bigl(w(\Sigma) \cdot u_{v,\Sigma},\, w(\Sigma') \cdot u_{v,\Sigma'}\bigr)\right|, \tag{6.6.1}$$

where in the last equality $\Sigma, \Sigma'$ are any two edges ending at $v$ (out of three, since $v$ is assumed to be trivalent; the balancing condition implies that the multipicity is independent of the choice of those two edges).

The *complex tropical multiplicity* of the map $(C, f)$ is

$$\mu^{\mathbf{C}}(C, f) \stackrel{\text{def.}}{=} \prod_{v \in \Gamma} \mu^{\mathbf{C}}(v).$$

[7] i.e., for all $s = 1, \dots, n$ and $i \geqslant 1$, the plane tropical curve $f(C)$ has exactly $\alpha_{s,i} + \beta_{s,i}$ infinite edges of weight precisely $i$ in the oriented direction $u_i$.

The enumeration of plane tropical curves may also be related with the enumeration of real algebraic plane curves, albeit with a different multiplicity: this is why we keep an index **C** in the notation.

**(6.7) Example.** The plane tropical curve in Example (5.9), seen as the image of a tropical map from the genus zero curve (b), has multiplicity one. Note that the vertex #6 of the plane curve has two points in its preimage, and does not correspond to a vertex of the abstract curve.

The plane tropical curve in Figure 7, seen as the image of a tropical map from a genus zero curve, has multiplicity four. Indeed, the two triangles marked with a star in the decomposition of the Newton polygon correspond to vertices $v$ of the abstract curve such that $\mu^{\mathbf{C}}(v) = 2$. See [IX, (5.3)] for a discussion (in the context of floor diagrams) of this multiplicity four in relation with the multiplicity of components of limit Severi varieties.

The following result, known as Mikhalkin's Correspondence Theorem, states the correspondence between the algebraic and tropical enumerative problems stated in (2.2) and (6.1) respectively. It is mainly due to Mikhalkin [23]; as for Proposition (6.4), the original result does not mention boundary conditions, which have first appeared in [14].

**(6.8) Theorem.** *Let $\Delta, g, \alpha, \beta$ be as in (6.1), and choose $\underline{x}, \Omega$ general enough. The following equality holds:*

$$\sum_{(C,f)\in\mathscr{E}_g^{\Delta}(\alpha,\beta)(\underline{x},\Omega)} \mu^{\mathbf{C}}(C,f) = I^{\alpha} \cdot N_g^{\Delta}(\alpha,\beta).$$

The starting point of the proof is the idea presented in Section 3 that tropical hypersurfaces are limits of toric hypersurfaces.

## 6.2 – Examples

**(6.9)** In the situation of Example (6.2) above, there is exactly one element $(C, f)$ in $\mathscr{E}_g^{\Delta}(\alpha,\beta)(\underline{x},\Omega)$. The corresponding decomposition of $\Delta$ is the trivial one, consisting only of the triangle $\Delta$ itself, which has Euclidean area $\mathscr{A}(\Delta) = 1$; one therefore has $\mu^{\mathbf{C}}(C, f) = 2$.

As the Correspondence Theorem tells us, this gives $N_g^{\Delta}(\alpha,\beta) = 2$, which indeed is the number of conics tangent to the line $z = 0$ at $(0 : 1 : 0)$, with prescribed general second order osculatory behaviour at this point, tangent to the line $y = 0$ at an unprescribed point, and passing through a general point $x_1$.

This is pictured below: two successive blow-ups of $\mathbf{P}^2$ are drawn, with exceptional divisors $E$ and $E'$; base points are represented by black dots or arrows when they are infinitely near, and the unassigned tangency condition is represented by an arrow with a dash.

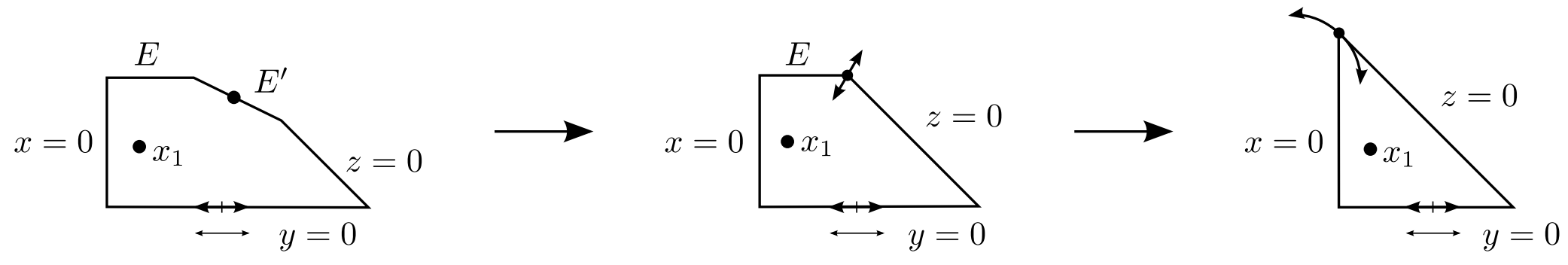

**(6.10)** Let us return to the enumerative problem (2.7): we let $\Delta$ be the Newton polygon indicated below, and fix $\alpha = (1, 1, 1, 0)$, $\beta = (0, 0, 0, 1)$. The set $\Omega$ consists of three real values, fixing three lines $\mathsf{L}_1, \mathsf{L}_2, \mathsf{L}_3$ as indicated on the picture below.

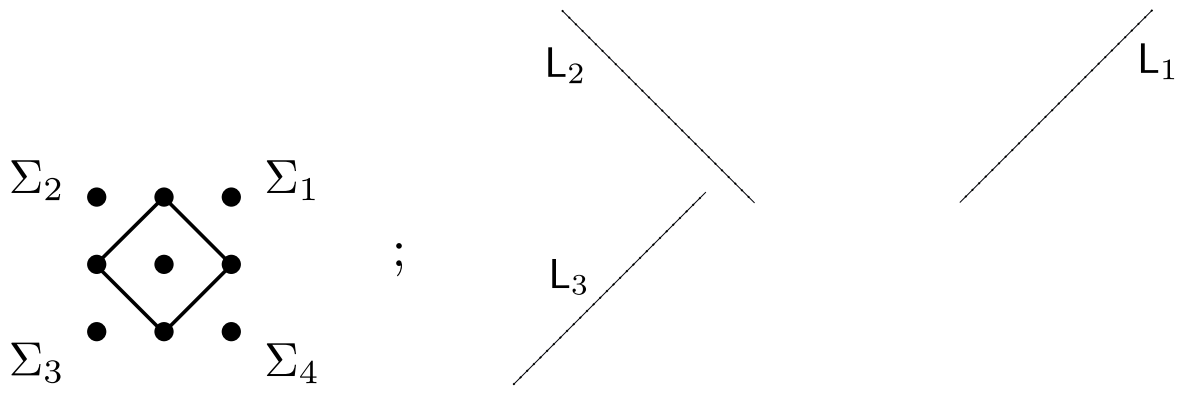

The tropical curves defined by this data move in a 1-dimensional family, in the way pictured below (two curves are pictured, one in green and one in blue): the deformations of the curves amount to change the size of the central rectangle by sliding one of its corners along a line $\mathsf{L}_i$, which lets all the rest follow. As the curves move in this family, one notices that the semi-line corresponding to the edge $\Sigma_4$ of the Newton polygon is always supported on the same line, which we shall call $\tilde{\mathsf{L}}_4$: this is the tropical incarnation of the toric Cayley–Bacharach property stated in Lemma (2.1).

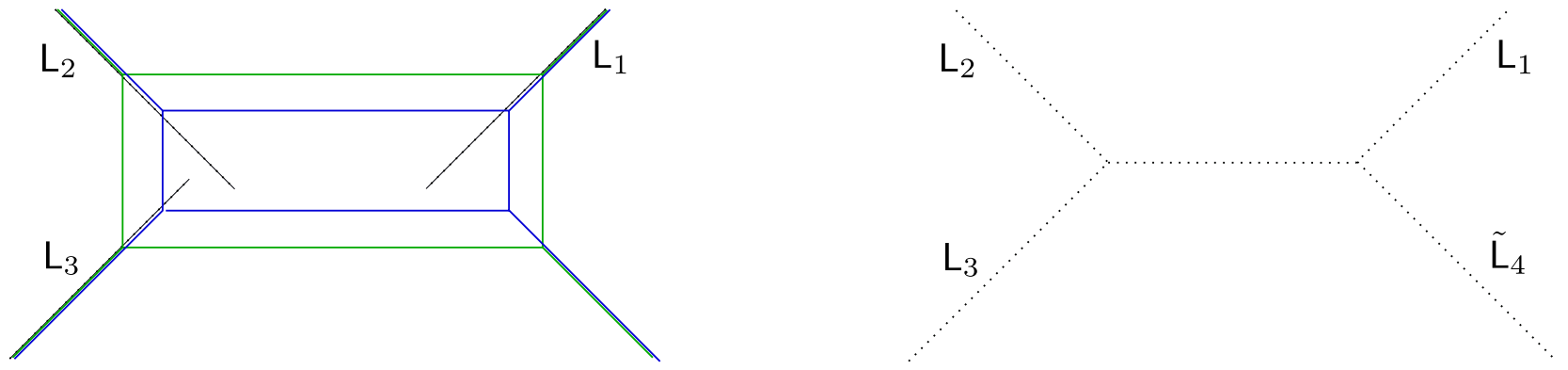

Figure 10: A pencil of tropical curves and the Cayley–Bacharach property

We want to count curves of genus $g = 0$ in this system. There is a finite number of such curves, thus the set $\underline{x}$ of base points is empty. There is a unique rational curve, depicted below together with the corresponding triangulation of the Newton polygon (note that a different choice for the lines $\mathsf{L}_1, \mathsf{L}_2, \mathsf{L}_3$ may lead to another, similar, triangulation).

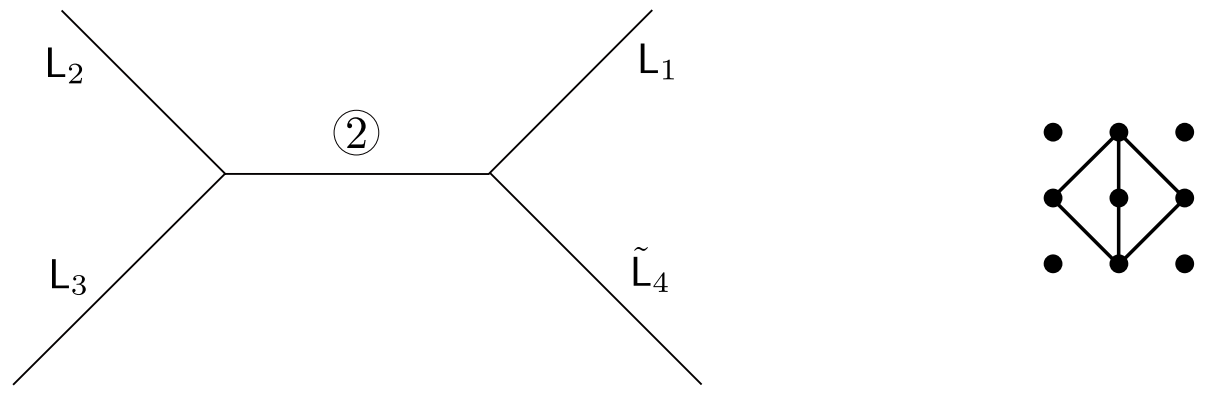

Both triangles in this triangulation yield the same multiplicity,

$$\mu(\Delta_v) = \left| \det \begin{pmatrix} 1 & 1 \\ 1 & -1 \end{pmatrix} \right| = 2.$$

One thus finds $N_0^\Delta\big((1,1,1,0),(0,0,0,1)\big) = 4$, in tune with the number previously found in (2.7).

**(6.11)** We now consider the enumerative problem (2.4): we let $\Delta$ be the Newton polygon indicated below, and set $g = 0$, $\alpha = (1,1,0)$, $\beta = (0,0,1)$. The set $\Omega$ consists of two real values, fixing two lines $\mathsf{L}_1, \mathsf{L}_2$ as indicated on the picture below, and the set $\underline{x}$ of general base points is empty.

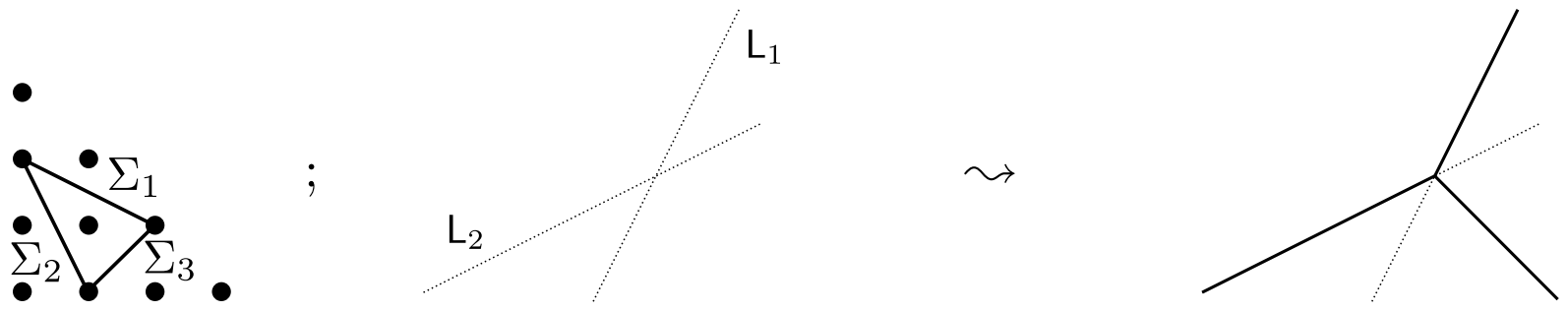

This enumerative problem is very similar to (2.7) and may be solved exactly as in (6.10) above. There is only one tropical map in $\mathcal{E}_0^{\Delta}\big((1,1,0),(0,0,1)\big)$; the corresponding subdivision of $\Delta$ is the trivial one, with only one triangle $\Delta_v$ giving the multiplicity

$$\mu(\Delta_v) = \left|\det\begin{pmatrix} 2 & 1 \\ -1 & -2 \end{pmatrix}\right| = 3.$$

One thus finds $N_0^{\Delta}\big((1,1,0),(0,0,1)\big) = 3$, in tune with the number previously found in (2.4).

**(6.12)** Let us consider one last enumerative problem, which is linked with (2.7). The problem is to enumerate genus 0 curves of type $(2,2)$ on $\mathbf{P}^1 \times \mathbf{P}^1$ tangent to all four points cut out on an anticanonical quadrilateral by a general curve of type $(1,1)$. As this imposes only seven independent linear conditions, this amounts to requiring the tangency at three fixed points and one mobile point.

Correspondingly, we consider the Newton polygon $\Delta$ indicated below, and set $g = 0$, $\alpha = (2,2,2,0)$ and $\beta = (0,0,0,2)$. The set $\Omega$ determines three lines $\mathsf{L}_1, \mathsf{L}_2, \mathsf{L}_3$, and the set $\underline{x}$ is empty.

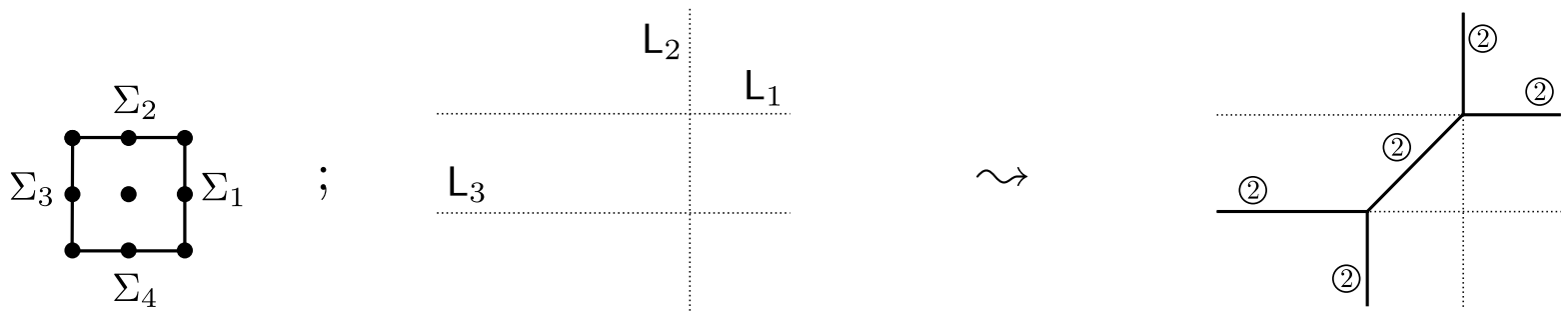

There is only one tropical map $(C,f)$ in $\mathcal{E}_0^{\Delta}(\alpha,\beta)$; note that its shape depends on the relative vertical position of the two parallel lines $\mathsf{L}_1$ and $\mathsf{L}_3$. In any event, there are five edges of weight 2, and the subdivision of $\Delta$ consists of two triangles, each with Euclidean area 2. The Correspondence Theorem thus gives

$$N_0^{\Delta}(\alpha,\beta) = \frac{1}{I^{\alpha}} \times \mu^{\mathbf{C}}(C,f) = \frac{4\times 4}{2^3} = 2.$$

The relation with (2.7) comes from seeing the quartic Del Pezzo surface with four ordinary double points of (2.7) as a double cover of $\mathbf{P}^1 \times \mathbf{P}^1$ branched over an anticanonical quadrilateral. The number 2 found above fits with the number 4 given in (2.7) and (6.10), the factor 2 between 4 and 2 being accounted for by the double cover.

For one more example, slightly more involved than the above ones, you can consult Example (7.18), in which the enumerative problem of (2.5) is considered.

# 7 – Floor diagrams and enumeration

This section is dedicated to the translation of the tropical enumerative problem (6.1) in purely combinatorial terms, under certain technical assumptions. We first define floors and elevators of a plane tropical curve, which will enable its encoding as a floor diagram. Next we give an abstract definition of floor diagrams, and discuss the relevant technical conditions. When this is done we state Brugallé and Mikhalkin's Correspondence Theorem. Many examples are discussed throughout. We use the notation of Section 5.2.

The floor diagrams that I define here are a mixture of those defined in [8] and [2], with a presentation according to my personal taste. As Brugallé taught me, floor diagrams are more a philosophy than a specific object. Besides, the floor diagrams considered in [13] are themselves slightly different from those in [8], see [13, Remark 1.2]. Note that the floor diagrams introduced in [IX] are a particular instance of those presented here, corresponding to the Newton polygons of complete linear systems of plane curves (i.e., the triangles $T_d$ of Example (1.1)).

**(7.1) Definition.** *Let* $\mathsf{d} = (a, b) \in \mathbf{Z}^2$ *be a primitive vector, and let* $(C, f)$ *be a parametrized plane tropical curve, i.e.,* $C = (\bar{\Gamma}, g)$ *is an abstract tropical curve with a tropical map* $f$ *to the plane. An* elevator *of* $(C, f)$ *is an edge* $\Sigma$ *of* $C$ *with oriented direction vector in the plane* $u_\Sigma = \pm\mathsf{d}$ *(see Definition (5.8)). A* floor *of* $(C, f)$ *is a connected component of* $\Gamma \setminus \mathrm{Elev}$*, where* $\Gamma = \bar{\Gamma} \setminus \Gamma^\infty$ *as in Definition (5.7), and* Elev *is the union of all elevators of* $(C, f)$*.*

**(7.2) Example.** Below, I explicit the respective decompositions of a genus 1 and a genus 0 plane tropical cubics into floors and elevators, assuming $\mathsf{d}$ is the vertical direction. On the left, I picture the plane tropical curve; in the middle, I separate the floors and elevators (the former are chains of line segments with no vertical components, the latter are vertical segments or semi-lines); on the right, I provide a schematic picture of the incidence between floors and diagrams, on which the floors are represented by ellipses and the elevators by vertical segments.

On Figure 12, note that we see the curve as a parametrized rational curve, so that the two segments intersecting at the black square are separate in the source curve, as in Example (5.9), case (b). The floors and elevators are defined at the level of the source curve.

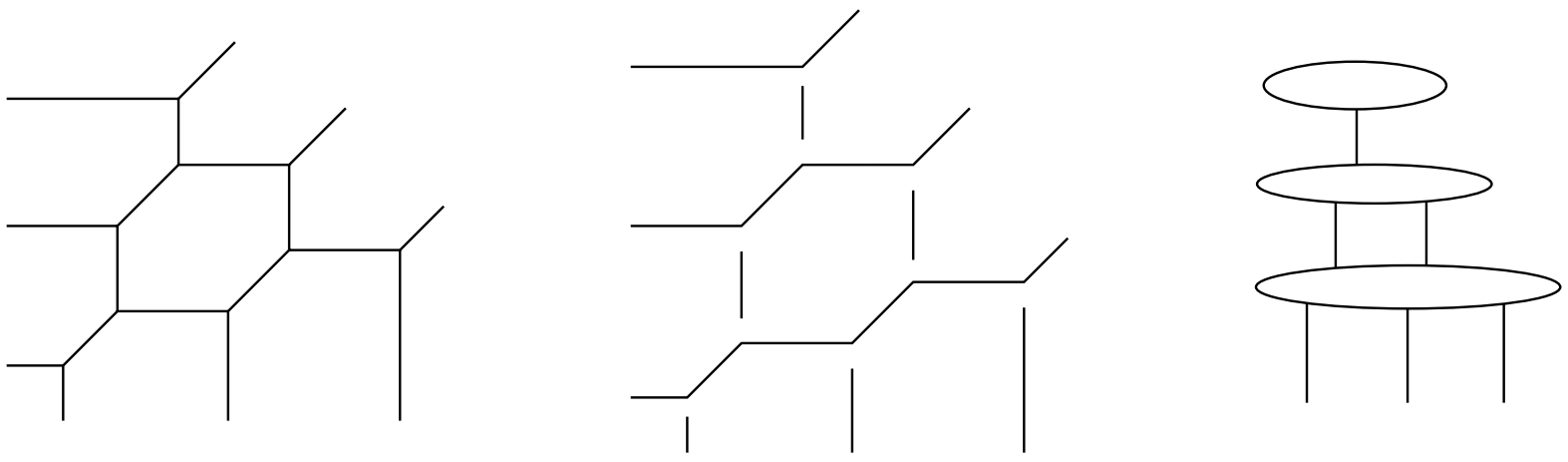

Figure 11: Floor decomposition of a genus 1 tropical cubic

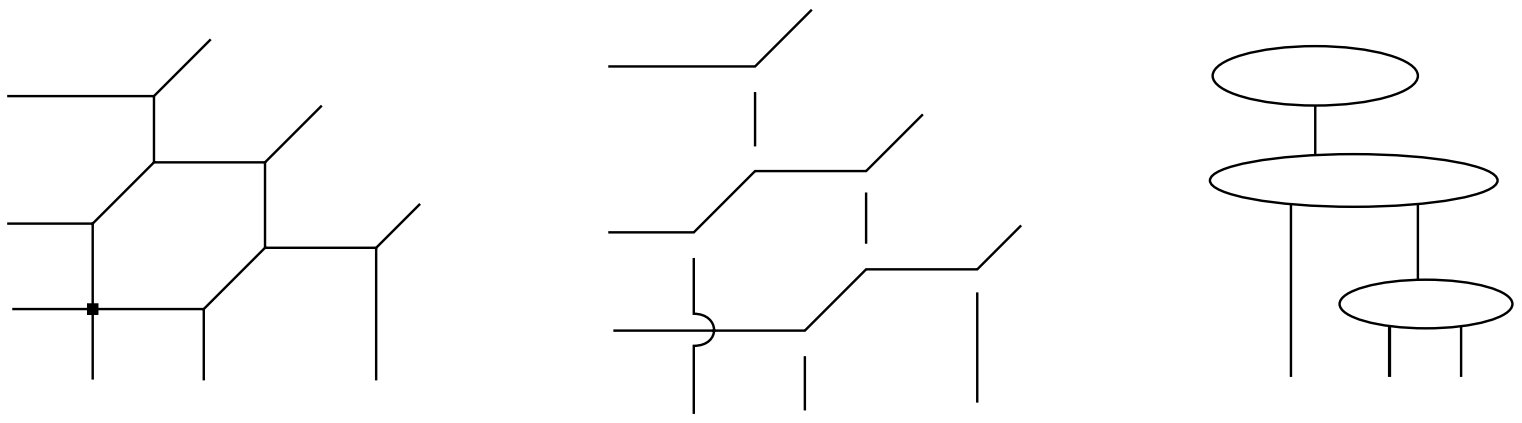

Figure 12: Floor decomposition of a genus 0 tropical cubic

**(7.3) Remark.** In general it is not true that floors are chains of line segments, or that an elevator intersects a given floor only once, see the examples in Figure 13. This will however be the case when we impose the passing through base points $\underline{x}$ sufficiently stretched out along direction $\mathsf{d}$.

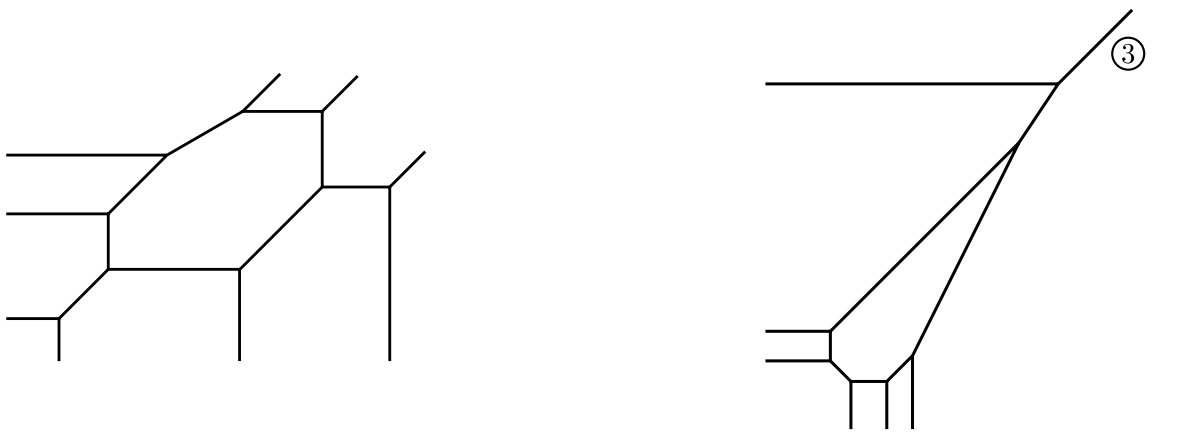

Figure 13: Tropical curves with messed up floor decomposition (see Remark (7.3))

It is also possible that the planar incarnations of two floors intersect, even under our general favourable framework when the points $\underline{x}$ are sufficiently stretched out along direction $\mathsf{d}$. See (7.15.1) for an example.

The feature that will enable the translation of a tropical enumerative problem into an enumeration of floor diagrams is the following: if the the base points $\underline{x}$ are sufficiently streched in the direction $\mathsf{d}$, then for all curve $(C,f) \in \mathcal{E}_g^\Delta(\alpha,\beta)(\underline{x},\Omega)$, each elevator (resp. floor) of $(C,f)$ will pass through exactly one point $x_i$ of $\underline{x}$. In particular, we will see that the cardinality of $\underline{x}$ equals the number of elevators and floors in $(C,f)$.

There are two technical conditions for the translation into an enumeration of floor diagrams to work correctly. The first one is that the Newton polygon $\Delta$ should be "transverse to the direction ${}^{\mathsf{T}}\mathsf{d}$", and the second one is that non-trivial boundary conditions (i.e., conditions of types $\alpha$ and $\beta$) may only be imposed corresponding to specific components of the boundary of $\Delta$, orthogonal to the direction $\mathsf{d}$.

**(7.4) Definition.** *Set* $\mathsf{d} = (a,b) \in \mathbf{Z}^2$ *with* $\gcd(a,b) = 1$. *Let* $\Delta$ *be a Newton polygon. The* left and right boundaries *of* $\Delta$ *with respect to* $\mathsf{d}$ *are, respectively,*

$$\partial_l\Delta = \left\{p \in \Delta :\ \forall t > 0,\ p + t.\mathsf{d}^\perp \notin \Delta\right\} \quad \text{and} \quad \partial_r\Delta = \left\{p \in \Delta :\ \forall t > 0,\ p - t.\mathsf{d}^\perp \notin \Delta\right\}.$$

*The Newton polygon* $\Delta$ *is* $\mathsf{d}^\perp$-transverse *if for all vector* $u \in \mathbf{Z}^2$, *parallel to an edge of* $\partial_l\Delta$ *or* $\partial_r\Delta$ *and with relatively prime coordinates, one has* $\det(u,\mathsf{d}^\perp) = \pm 1$.

The original definition [8, 2.2] corresponds to $\mathsf{d} = (0,1)$, and is called $h$-transversality there, '$h$' standing for 'horizontal'.

The geometric meaning of $\mathsf{d}^\perp$-transversality is that for all tropical curves $\Gamma$ with Newton diagram $\Delta$, any (weight one) line in $\mathbf{R}^2$ directed by $\mathsf{d}$ should intersect each infinite branch $e$ of $\Gamma$ with multiplicity equal to the weight $w(e)$ (except if the branch has the same direction $\mathsf{d}$). This implies that the total number of intersection of such a line with $\Gamma$ is the number of floors of $\Gamma$. The definition of intersection multiplicity to make this work is the following: two lines directed by primitive vectors $u_1, u_2 \in \mathbf{Z}^2$ and with weights $w_1$ and $w_2$ intersect with multiplicity $w_1 \cdot w_2 \cdot \det(u_1,u_2)$.

**(7.5) Examples.** The triangle $T_d$, corresponding to the complete linear system $|dH|$ on $\mathbf{P}^2$ (see (1.1)), is $\mathsf{d}^\perp$-transverse for $\mathsf{d} = (0,1)$.

The trapezium $\mathrm{Tz}^r_{a,b}$, corresponding to the complete linear system $|aH + bF|$ on the Hirzebruch surface $\mathbf{F}_r$ (see (1.6)), is $\mathsf{d}^\perp$-transverse for $\mathsf{d} = (0,1)$, but if $r > 1$ it is not for $\mathsf{d} = (1,0)$.

The two polygons given in (2.5) and (2.7) are $\mathsf{d}^\perp$-transverse for $\mathsf{d} = (0,1)$, but neither for $\mathsf{d} = (1,1)$, nor for $\mathsf{d} = (1,-1)$.

The triangle in (2.4) is not $\mathsf{d}^\perp$-transverse either for $\mathsf{d} = (0,1)$, or $(1,0)$, or $(1,1)$, or $(1,-1)$.

**(7.6)** The second condition to be fulfilled by a tropical enumerative problem in order to be translatable into a floor diagram problem is that one should impose non-trivial boundary conditions only along those edges of $\Delta$ that are orthogonal to $\mathsf{d}$, or equivalently along the one or two edges from $\partial\Delta$ off $\partial_l\Delta \cup \partial_r\Delta$.

In other words, in the notation of (6.1), for all $s = 1,\dots,n$, if $\Sigma_s$ lies on the left or right boundary of $\Delta$ with respect to $\mathsf{d}$, then $\alpha_s = 0$ and $\beta_s = [\ell_{\mathbf{Z}}(\Sigma_s)]$.

We need the following before we can define floor diagrams.

**(7.7)** Let $\Delta$ be a $\mathsf{d}^\perp$-transverse Newton polygon, for some primitive vector $\mathsf{d} \in \mathbf{Z}^2$. Let $\Sigma$ be an edge of $\Delta$ lying on its left or right boundary, and view it as a vector going down with respect to the direction $\mathsf{d}$: this vector writes $\ell_{\mathbf{Z}}(\Sigma).u$, with $u$ a primitive integer vector such that $\det(\mathsf{d}^\perp, u) = +1$. The list of directions associated to $\Sigma$ is the unordered list $[u^\perp, \ldots, u^\perp]$ consisting of the primitive vector $u^\perp$ repeated $\ell_{\mathbf{Z}}(\Sigma)$ times. The *list of left (resp. right) directions* of $\Delta$ is the unordered list $D_l(\Delta)$ (resp. $D_r(\Delta)$) obtained by concatenating the lists of directions of all edges $\Sigma$ lying on the left (resp. right) boundary of $\Delta$.

Note that, for tropical curves $(C, f)$ with Newton polygon $\Delta$, the list of left directions of $\Delta$ correspond to the direction vectors $u_{\infty,\Sigma}$ pointing from infinity inwards of infinite branches of coming from the left with respect to the direction $\mathsf{d}$, whereas the right directions correspond to vectors $u_{\Sigma,\infty}$ pointing outwards to infinity for infinite branches going to the right. In other words, the direction is always from left to right with respect to the direction $\mathsf{d}$.

Lastly, let $d_+(\Delta)$ (resp. $d_-(\Delta)$) be the integral length of the top (resp. bottom) edge of $\Delta$, parallel to the direction $\mathsf{d}^\perp$. If there is no such edge, we set the corresponding value $d_\pm(\Delta)$ to 0.

The $\mathsf{d}^\perp$-transversality condition ensures that $D_l(\Delta)$ and $D_r(\Delta)$ have the same number of elements[8], which we shall call the $\mathsf{d}$-height of the polygon $\Delta$, and denote by $d_{\mathsf{d}}(\Delta)$. Thus

$$\begin{aligned} d_{\mathsf{d}}(\Delta) &= \operatorname{length}(D_l(\Delta)) = \operatorname{Card}(\partial_l(\Delta) \cap \mathbf{Z}^2) - 1 \\ &= \operatorname{length}(D_r(\Delta)) = \operatorname{Card}(\partial_r(\Delta) \cap \mathbf{Z}^2) - 1. \end{aligned}$$

Moreover the following holds,

$$2d_{\mathsf{d}}(\Delta) + d_+(\Delta) + d_-(\Delta) = \operatorname{Card}(\partial\Delta_{\mathbf{Z}}). \tag{7.7.1}$$

**(7.8) Examples.** The following are the examples from (7.5) continued. We take $\mathsf{d} = (0, 1)$ in all cases.
For the triangle $T_d$, one has $d_{\mathsf{d}} = d_- = d$, $d_+ = 0$, and

$$D_l = \big[(1,0), \ldots, (1,0)\big] \quad \text{and} \quad D_r = \big[(1,1), \ldots, (1,1)\big].$$

For the trapezium $\mathrm{Tz}^r_{a,b}$, one has $d_{\mathsf{d}} = a$, $d_- = b + ar$, $d_+ = b$, and

$$D_l = \big[(1,0), \ldots, (1,0)\big] \quad \text{and} \quad D_r = \big[(1,r), \ldots, (1,r)\big].$$

For the polygon in (2.7), one has $d_{\mathsf{d}} = 2$, $d_- = d_+ = 0$, and

$$D_l = D_r = \big[(1,-1), (1,1)\big].$$

For the polygon in (2.5), one has $d_{\mathsf{d}} = 3$, $d_- = d_+ = 0$, and

$$D_l = D_r = \big[(1,-1), (1,1), (1,1)\big].$$

Let us now give an abstract definition of floor diagrams.

[8] by a suitable change of coordinates in $\mathrm{SL}_2(\mathbf{Z})$ we may reduce to the case $\mathsf{d} = (0, 1)$, and then the transversality condition is that any primitive vector $u \in \mathbf{Z}^2$ parallel to an edge on the left or right boundary of $\Delta$ is of the form $(\alpha, \pm 1)$, $\alpha \in \mathbf{Z}$.

**(7.9) Weighted oriented graphs.** Consider a connected, acyclic, oriented graph $\mathcal{D}$, with finitely many vertices and edges. We let $\mathrm{Vert}(\mathcal{D})$ and $\mathrm{Edges}(\mathcal{D})$ be the sets of vertices and edges of $\mathcal{D}$ respectively (or simply Vert and Edges if no confusion is likely), and see Edges as a subset of $\mathrm{Vert}\times\mathrm{Vert}$. We call $\mathcal{D}_{\mathrm{top}}$ the topological incarnation of $\mathcal{D}$ as a 1-dimensional CW-complex.

We assume that there is a distinguished subset of 1-valent vertices of $\mathcal{D}$, which we call the *vertices at infinity*. We let $\mathrm{Vert}^{-\infty}$ (resp. $\mathrm{Vert}^{+\infty}$) be the set of vertices at infinity whose unique adjacent edge is outgoing (resp. ingoing). We let $\mathrm{Vert}^{\mathrm{fin}}$ be the complement in Vert of $\mathrm{Vert}^{-\infty}\coprod\mathrm{Vert}^{+\infty}$, and call those vertices in $\mathrm{Vert}^{\mathrm{fin}}$ *floors*.

We assume that there is no edge joining two vertices at infinity. Thus we have a partition $\mathrm{Edges}=\mathrm{Edges}^{-\infty}\coprod\mathrm{Edges}^{\mathrm{fin}}\coprod\mathrm{Edges}^{+\infty}$, where $\mathrm{Edges}^{\pm\infty}$ is the set of those edges adjacent to a vertex in $\mathrm{Vert}^{\pm\infty}$. Thus $\mathrm{Edges}^{\mathrm{fin}}$ is the set of those edges joining two floors. The edges in $\mathrm{Edges}^{-\infty}\coprod\mathrm{Edges}^{+\infty}$ are called *infinite*, and the others *finite*. In the terminology of Definition (7.1), one calls the edges of $\mathcal{D}$ (either finite or infinite) *elevators*.

We let $\mathsf{D}$ be the abstract set $\mathrm{Vert}^{\mathrm{fin}}\coprod\mathrm{Edges}$, consisting of floors and all edges (finite or infinite) of $\mathcal{D}$. Since $\mathcal{D}$ is acylic, the set $\mathsf{D}$ inherits a partial ordering from the orientation.

We will in addition assume that $\mathcal{D}$ is weighted, which means that it comes equipped with a *weight function* $w:\mathrm{Edges}\to\mathbf{N}^*$. One then defines the *divergence* $\mathrm{div}(v)$ of a vertex $v$ (either a vertex at infinity or a floor) as

$$\mathrm{div}(v)=\textstyle\sum_{(s,v)\in\mathrm{Edges}}w(s,v)-\sum_{(v,t)\in\mathrm{Edges}}w(v,t)$$

("sum of weights of incoming edges minus sum of weights of outgoing edges").

**(7.10) Definition.** *Let $\Delta$ be a $\mathsf{d}^{\perp}$-transverse lattice polygon, and $g$ be a non-negative integer. A (plane)* floor diagram $(\mathcal{D},\theta)$ *of genus $g$ and Newton polygon $\Delta$ with respect to direction $\mathsf{d}$ is the data of an acyclic connected weighted oriented graph $\mathcal{D}$ and a map $\theta:\mathrm{Vert}(\mathcal{D})^{\mathrm{fin}}\to\mathbf{Z}$, such that the following conditions hold:*

*(i) the first Betti number $b_1(\mathcal{D}_{\mathrm{top}})$ equals $g$;*

*(ii) one has $\sum_{v\in\mathrm{Vert}^{-\infty}(\mathcal{D})}\mathrm{div}(v)=-d_-(\Delta)$ and $\sum_{v\in\mathrm{Vert}^{+\infty}(\mathcal{D})}\mathrm{div}(v)=d_+(\Delta)$;*

*(iii) the unordered list $\left[\theta(v),v\in\mathrm{Vert}(\mathcal{D})^{\mathrm{fin}}\right]$ equals $D_l(\Delta)$;*

*(iv) the unordered list $\left[\theta(v)+\mathrm{div}(v).\mathsf{d},v\in\mathrm{Vert}(\mathcal{D})^{\mathrm{fin}}\right]$ equals $D_r(\Delta)$.*

To understand condition (iv), note that it follows from the balancing condition for plane tropical curves that if a floor of some tropical curve has direction at infinity from the left $u$, and is intersected from below by $k$ edges (counted with weights) with direction $\mathsf{d}$ and from above by $l$ edges (counted with weights) with direction $\mathsf{d}$, then its direction at infinity to the right is $u+(k-l).\mathsf{d}$. For conditions (iii) and (iv) to be compatible, one needs that $D_l(\Delta)$ and $D_r(\Delta)$ have the same length, which is granted by the $\mathsf{d}^{\perp}$-transversality condition, see (7.7); then $d_{\mathsf{d}}(\Delta)$ is the number of floors of a floor diagram with Newton polygon $\Delta$.

In condition (ii), $\pm\sum_{v\in\mathrm{Vert}^{\pm\infty}(\mathcal{D})}\mathrm{div}(v)$ is the number of edges in $\mathrm{Edges}^{\pm\infty}$, counted with weights; we may also refer to it as the number of vertices in $\mathrm{Vert}^{\pm\infty}$ counted with weights (even though, strictly speaking, vertices at infinity don't have weights).

The number of vertices of a diagram $\mathcal{D}$ as above is thus, counting the vertices at infinity with weights,

(7.10.1) $$\mathrm{Card}_w(\mathrm{Vert}(\mathcal{D}))=d_{\mathsf{d}}(\Delta)+d_+(\Delta)+d_-(\Delta)=\mathrm{Card}(\partial\Delta_{\mathbf{Z}})-d_{\mathsf{d}}(\Delta),$$

where the last equality is given by Relation (7.7.1) (the subscript '$w$' indicates that the vertices should be counted with weigths). On the other hand, the Euler formula for $\mathcal{D}_{\mathrm{top}}$,

(7.10.2) $$\mathrm{Card}(\mathrm{Edges}(\mathcal{D}))=g-1+\mathrm{Card}(\mathrm{Vert}(\mathcal{D})),$$

implies, since there is exactly one infinite edge for each vertex at infinity,

$$\text{Card}_w(\text{Edges}(\mathcal{D})) = g - 1 + \text{Card}_w(\text{Vert}(\mathcal{D})), \tag{7.10.3}$$

where infinite edges are counted with weights (be careful that, on the contrary, finite edges are not counted with weigths). Relations (7.10.1) and (7.10.3) yield

$$\text{Card}_w(\mathsf{D}) = d_{\mathsf{d}}(\Delta) + \text{Card}_w(\text{Edges}(\mathcal{D})) = \text{Card}(\partial\Delta_{\mathbf{Z}}) + g - 1, \tag{7.10.4}$$

where infinite edges are counted with weights.

We adopt the following graphical convention: floors will be represented by ellipses, and vertices at infinity will be omitted; edges will be represented by vertical line segments, the orientation is from the bottom to the top. The weight of an edge is indicated only if it is greater than 1, as a circled integer. We indicate $\theta(v)$ inside the ellipse $v$, except if $\theta$ takes the same value on all floors in which case we don't indicate it.

**(7.11) Examples.** Let us continue with the examples from (7.5) and (7.8). For each of these I picture below the unique floor diagram corresponding to curves with the maximal genus, with respect to the direction $\mathsf{d} = (0, 1)$.

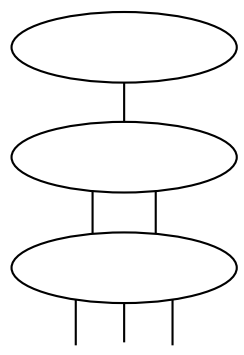

Figure 14: Floor diagram of Newton polygon $T_3$ and genus 1

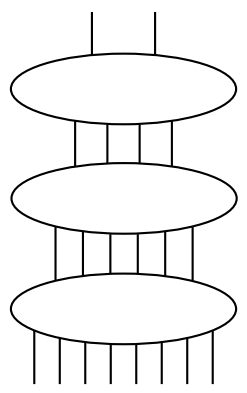

Figure 15: Floor diagram of Newton polygon $\text{Tz}^2_{3,2}$ (see Figure 1) and genus 8

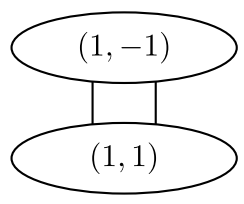

Figure 16: Floor diagram of genus 1 for the Newton polygon of (2.7)

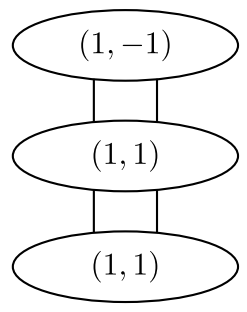

Figure 17: Floor diagram of genus 2 for the Newton polygon of (2.5)

In order to appropriately count floor diagrams, we need to enhance them with markings, which encode the distribution of the base points $\underline{x}$ and of the assigned supports for the boundary conditions (encoded in the set $\Omega$). To understand the following definition, keep in mind that the idea is that if the points $\underline{x}$ are sufficiently stretched in the direction $\mathsf{d}$, then for all curves $(C, f) \in \mathcal{E}^{\Delta}_g(\alpha, \beta)$, the points $\underline{x}$ distribute one on each floor or elevator, except on those elevators that are already taken by a type $\alpha$ boundary condition. Moreover, we shall label the points $\underline{x}$ according to their order along the direction $\mathsf{d}$, which is where the monotonicity condition (ii) comes from.

Due to our requirement (7.6), the boundary conditions reduce to only $\alpha^{\pm}, \beta^{\pm} \in \underline{\mathbf{N}}$, tangency conditions along the top and bottom edges of $\Delta$ respectively, parallel to direction $\mathsf{d}^{\perp}$ (with '+' for the top and '−' for the bottom), at fixed and mobile points respectively.

**(7.12) Definition.** *Let $\Delta$ be a $\mathsf{d}^{\perp}$-transverse lattice polygon, and consider $g \in \mathbf{N}$, and $\alpha^+, \alpha^-, \beta^+, \beta^- \in \underline{\mathbf{N}}$, subject to the conditions*

$$0 \leqslant g \leqslant \mathrm{Card}(\mathring{\Delta}_{\mathbf{Z}}) \quad \textit{and} \quad I\alpha^{\pm} + I\beta^{\pm} = d_{\pm}(\Delta).$$

*Set $s = g - 1 + 2d_{\mathsf{d}}(\Delta) + |\beta^+| + |\beta^-|$. A* marking of type *$(\alpha^+, \alpha^-, \beta^+, \beta^-)$ of a floor diagram $(\mathcal{D}, \theta)$ of genus $g$ and Newton polygon $\Delta$ is a map $m : \{-|\alpha^-| + 1, \dots, 1, \dots s + |\alpha^+|\} \to \mathsf{D}$ such that the following conditions hold:*

*(i) the map $m$ is injective;*

*(ii) if $m(i) > m(j)$ for the order of $\mathsf{D}$, then $i > j$;*

*(iii) for all $k \in \mathbf{N}^*$, if $-|\alpha^-| + 1 + \sum_{l=1}^{k-1} \alpha_l^- \leqslant i \leqslant -|\alpha^-| + \sum_{l=1}^{k} \alpha_l^-$ (resp. $s + 1 + \sum_{l=1}^{k-1} \alpha_l^+ \leqslant i \leqslant s + \sum_{l=1}^{k} \alpha_l^+$), then $m(i) \in \mathrm{Edges}^{-\infty}$ (resp. $m(i) \in \mathrm{Edges}^{+\infty}$), and the edge $m(i)$ has weight $k$;*

*(iv) for all $k \in \mathbf{N}^*$, there are exactly $\beta_k^{\pm}$ elements of $\mathrm{Edges}^{\pm\infty}$ of weight $k$ in $m(\{1, \dots, s\})$.*

A floor diagram $\mathcal{D}$ enhanced with a marking $m$ of type $(\alpha^+, \alpha^-, \beta^+, \beta^-)$ is called a *marked floor diagram* $(\mathcal{D}, m)$ of type $(\alpha^+, \alpha^-, \beta^+, \beta^-)$.

Two marked floor diagrams of the same type are *equivalent* if there exists a homeomorphism of oriented graphs $\phi : \mathcal{D} \cong \mathcal{D}'$ such that $w = w' \circ \phi$, $\theta = \theta' \circ \phi$, and $m = \phi \circ m'$.

Now, the fundamental result is the following. The enumerative application is given in Corollary (7.16) below.

**(7.13) Theorem** (Brugallé–Mikhalkin [8, Prop. 5.9])**.** *Maintain the situation of Definition (7.12) above, and consider the corresponding enumerative problem as in (6.1). If the points $\underline{x}$ are sufficiently strectched along direction $\mathsf{d}$, then there is a bijection between $\mathcal{E}_g^{\Delta}(\alpha, \beta)(\underline{x}, \Omega)$ and the set of equivalence classes of marked floor diagrams of type $(\alpha^+, \alpha^-, \beta^+, \beta^-)$.*

I shall not prove this result and refer instead to [8, §5], but I will try to convey some ideas why it is true. First of all, it follows from Definition (7.1) how to associate a floor diagram to a tropical curve, with the same genus and Newton polygon. Moreover, each tangency condition of order $k$ along the boundary imposed to the tropical curve translates into a condition that the floor diagram has an infinite edge of weight $k$.

These features of the curve being fixed, we need to let it move so as to meet the points $\underline{x}$ and $\Omega$. One has to convince oneself that the degrees of freedom in order to do so are to let floors move along the direction $\mathsf{d}$ (which correspondingly stretches the elevators), and let elevators move along the direction $\mathsf{d}^{\perp}$. By the way, note that indeed the number of floors and elevators (which is also the cardinality of $\mathsf{D}$) equals the number of points in $\underline{x}$ and $\Omega$, i.e., $s + |\alpha^+| + |\alpha^-|$ in the above notation, by the version of formula (7.10.4) where we no longer count infinite edges with weights: the number of vertices of $\mathcal{D}$ is

$$\mathrm{Card}(\mathrm{Vert}(\mathcal{D})) = d_{\mathsf{d}}(\Delta) + |\alpha^+| + |\beta^+| + |\alpha^-| + |\beta^-|$$

(compare with (7.10.1)), hence by the Euler formula (7.10.2)

$$\begin{aligned} \mathrm{Card}(\mathsf{D}) &= \mathrm{Card}(\mathrm{Edges}(\mathcal{D})) + \mathrm{Card}(\mathrm{Vert}(\mathcal{D})) - \big(|\alpha^+| + |\beta^+| + |\alpha^-| + |\beta^-|\big) \\ &= g - 1 + 2d_{\mathsf{d}}(\Delta) + |\alpha^+| + |\beta^+| + |\alpha^-| + |\beta^-| = s + |\alpha^+| + |\alpha^-|, \end{aligned}$$

as we said.

How the points $\underline{x}$ and $\Omega$ are taken up by the curve is encoded in the marking of the corresponding diagram. The points $\Omega$ must be taken up by pairwise distinct infinite elevators, wich correspond to infinite edges of the floor diagram. Next, when the points $\underline{x}$ are sufficiently distant one from another in direction d, it happens that they are taken up one by one by the floors and remaining elevators. The condition that the marking $m$ be increasing, i.e., (ii) in Definition (7.12), is satisfied when one orders the points $\underline{x}$ following direction d; note that the increasing condition is meaningful only in the range $\{1, \ldots, s\}$ (which corresponds to points $\underline{x}$), as it is otherwise implied by the requirement that $m(i)$ is in $\mathrm{Edges}^{-\infty}$ if $i < 1$ and in $\mathrm{Edges}^{+\infty}$ if $i > s$.

Lastly one has to convince oneself that, given the incidence relations prescribed by a marked floor diagram, there is a unique corresponding tropical curve. The idea is this: consider a floor in the diagram; in the plane it corresponds to a chain of line segments, intersecting various elevators coming from below or from above (with respect to direction d) in a way prescribed by the diagram; each of these elevators has to pass through one base point prescribed by the marking, so that the line supporting each elevator in the plane is known. The slope of the leftmost (with respect to d) line segment in the chain is prescribed by the function $\theta$ attached to the diagram, and this in turn determines the slopes of the other segments via the balancing condition, see also the proof of Corollary (7.16) below. At this point, we know the slopes of all segments in the floor, and we know that the changes of slopes occur along the lines supporting the elevators; therefore the incarnation of the floor in the plane is determined up to translation along the direction d. Therefore, there is a unique translate passing through the base point prescribed by the marking. It is arguably best to experiment this construction by one's self, playing around with examples; I provide some material to do so in the following example.

**(7.14) Example.** We display all marked floor diagrams of rational plane cubics with trivial boundary conditions in Figures 18, 19, and 20 below.

For the first diagram (Figure 18), we picture the corresponding tropical curve for two different distributions of the base points (always stretched out in the vertical direction). For both of them, I try to explain how to construct the curves in practice by showing the constructions lines that I have used to find the curves, following the procedure explained in the previous paragraph.

The five next marked diagrams (Figure 19) correspond to five different markings of the same diagram. I give the corresponding plane tropical curves for a given distribution of the base points (the same distribution for all five curves). I use black circles to picture the base points, and black squares to picture the "nodes" of the plane tropical curve (one on each curve in this case).

The three remaining marked diagrams are pictured on Figure 20, and I leave it as an exercise to construct the corresponding tropical curves for a given sufficiently stretched distribution of base points.

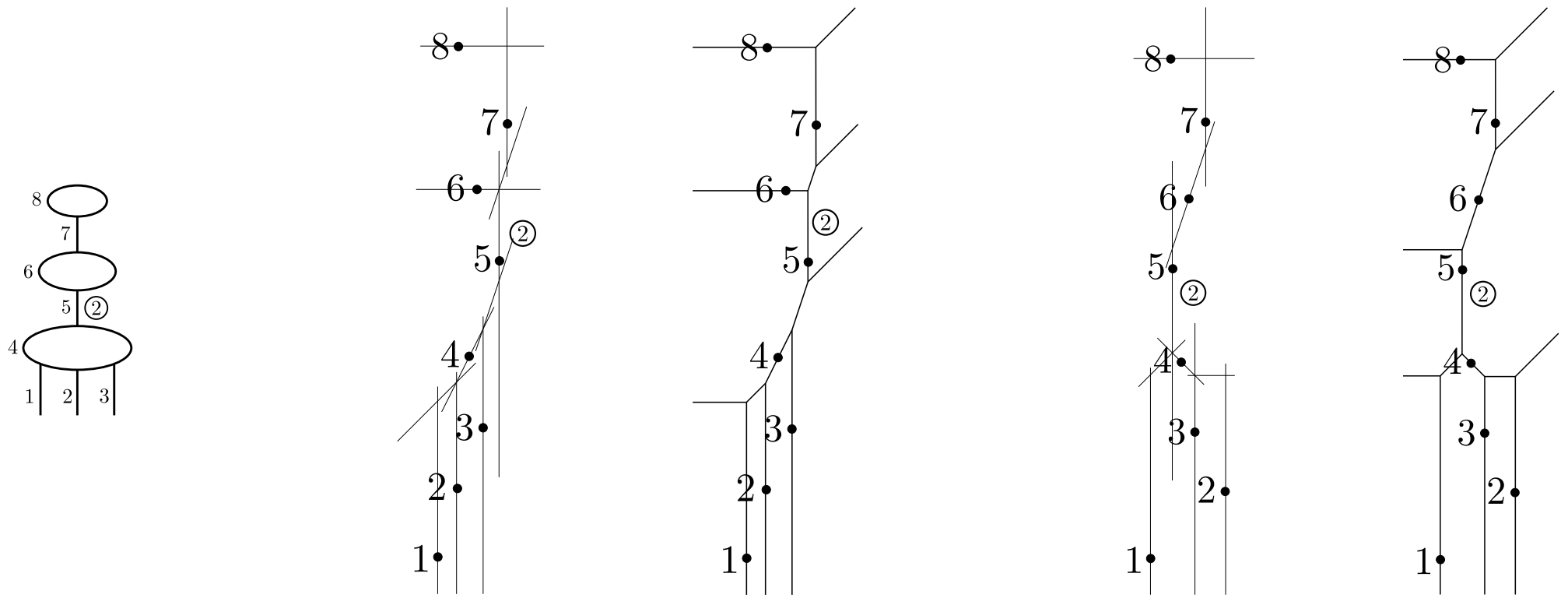

Figure 18: Tropical cubics with the same marked floor diagram for two different distributions of the base points (for each, on the left: construction lines)

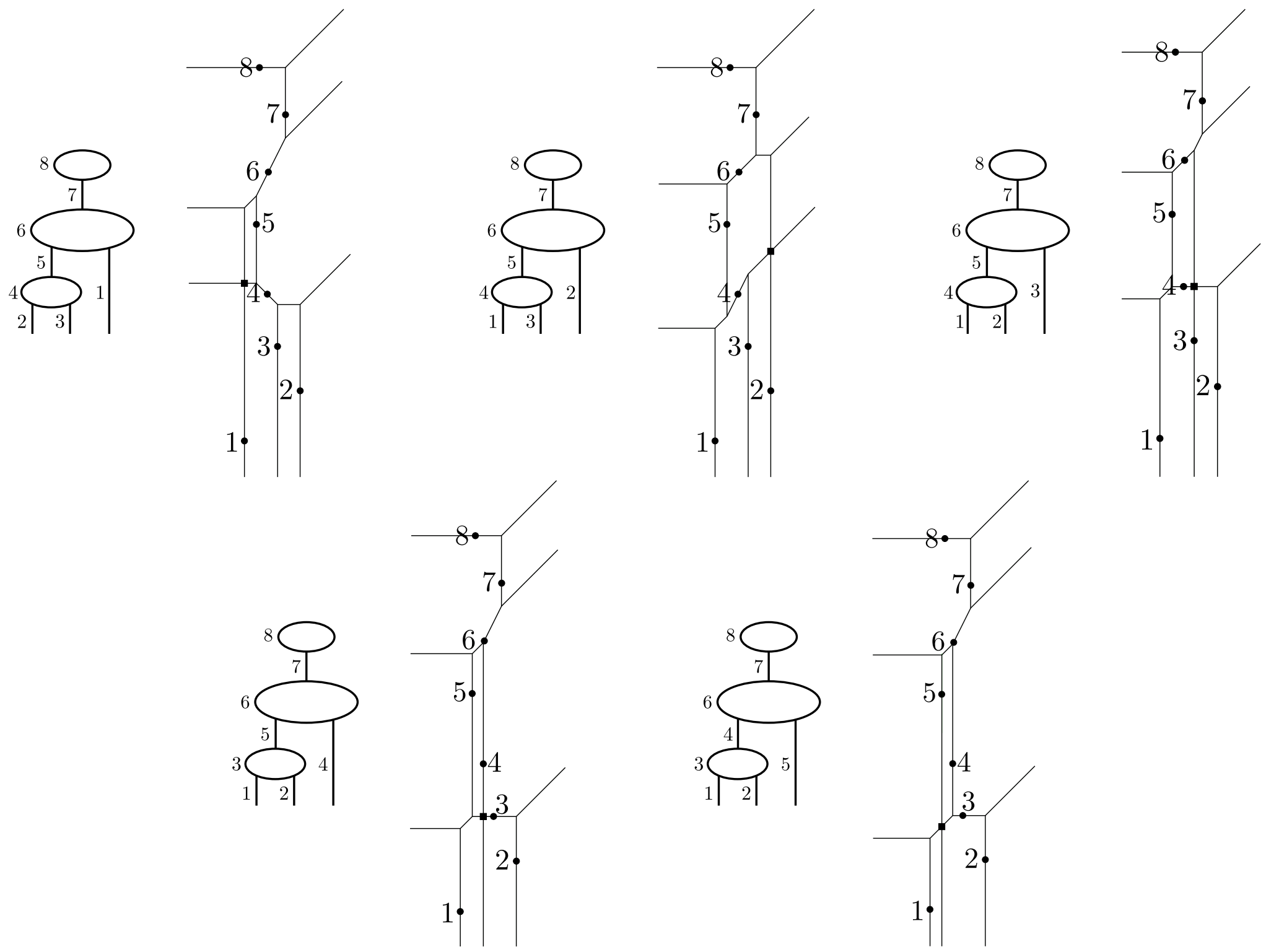

Figure 19: Tropical cubics with the same floor diagram, all possible markings

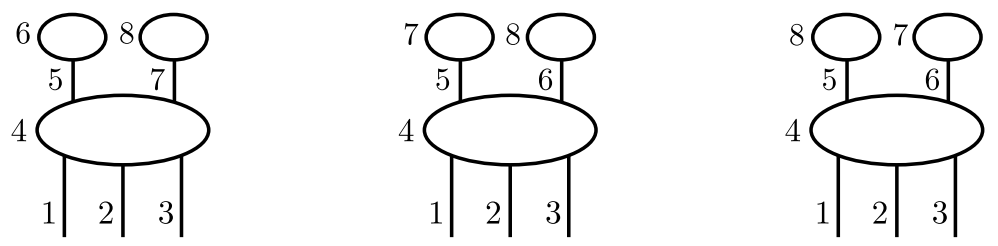

Figure 20: Remaining marked floor diagrams for rational plane cubics

**(7.15) Warnings.** Here are a few caveats about the correspondence between tropical curves and floor diagrams, which I feel are helpful to properly understand what is going on.

*(7.15.1)* There is no reason why the left infinite directions of floors (those encoded by $\theta$) would be ordered following the order in which they appear on the Newton polygon: this is the reason why the lists $D_l(\Delta)$ and $D_r(\Delta)$ are unordered.

An example is given in the figure below (you should try and draw the tropical curve for yourself). We take $\mathsf{d} = (0, 1)$; the Newton polygon pictured on the left is indeed $\mathsf{d}^{\perp}$-transverse. Note that the floor diagram has genus 3, while the arithmetic genus of the Newton polygon is 5. The singularities of a curve with the indicated floor diagram are to be found as intersection points between the planar incarnation of two distinct floors.

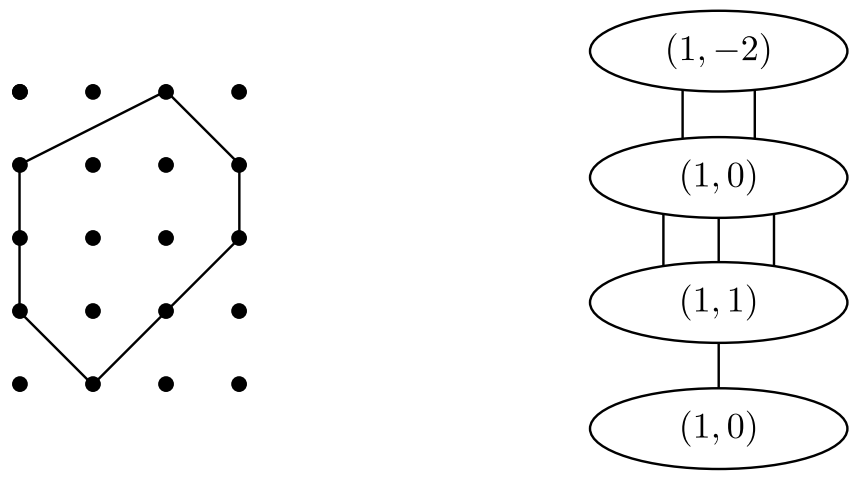

*(7.15.2)* In general the datum of the floor diagram does not impose a particular subdivision of the Newton polygon. This can be observed on the examples of Figure 18 above.

Finally, let us give the application of Theorem (7.13) to the enumerative problem (2.2).

**(7.16) Corollary.** *Let $\Delta$ be a $\mathsf{d}^\perp$-transverse lattice polygon, and consider $g \in \mathbf{N}$, and $\alpha^+, \alpha^-, \beta^+, \beta^- \in \underline{\mathbf{N}}$, subject to the conditions*

$$0 \leqslant g \leqslant \mathrm{Card}(\mathring{\Delta}_{\mathbf{Z}}) \quad \textit{and} \quad I\alpha^\pm + I\beta^\pm = d_\pm(\Delta).$$

*The number $N_g^\Delta(\alpha^+, \alpha^-, \beta^+, \beta^-)$ equals the number of equivalence classes of marked floor diagrams of genus $g$, Newton polygon $\Delta$, and type $(\alpha^+, \alpha^-, \beta^+, \beta^-)$, where each marked diagram $(\mathcal{D}, m)$ is to be counted with the multiplicity*

$$\mu^{\mathbf{C}}(\mathcal{D}, m) = I^\beta \prod_{e \in \mathrm{Edges}^{\mathrm{fin}}} w(e)^2.$$

*Proof.* By Mikhalkin's correspondence theorem (6.8), the number $N_g^\Delta(\alpha^+, \alpha^-, \beta^+, \beta^-)$ equals $I^{-\alpha}$ times the number of parametrized plane tropical curves $(C, f) \in \mathcal{E}_g^\Delta(\alpha^+, \alpha^-, \beta^+, \beta^-)(\underline{x}, \Omega)$, for a general choice of points $\underline{x}$ and $\Omega$, each tropical curve being counted with the multipicity $\mu^{\mathbf{C}}(C, f)$ defined in (6.6). In turn, this number counts marked floor diagrams if we take the points $\underline{x}$ as in Theorem (7.13). It thus only remains to show that the multiplicities of the tropical curves in $\mathcal{E}_g^\Delta(\alpha^+, \alpha^-, \beta^+, \beta^-)(\underline{x}, \Omega)$ relate to that of the corresponding marked floor diagrams in the appropriate way.

Let us consider a curve $(C, f)$ as above. In the corresponding subdivision of the Newton polygon $\Delta$, each sub-polygon is a triangle corresponding to a vertex at which a floor intersects an elevator and changes slope. Let us thus consider a floor of $(C, f)$, and name $u_0, \ldots, u_n$ the direction vectors of the various line segments composing this floor, oriented and numbered from left to right with respect to direction $\mathsf{d}$ as indicated on the figure below.

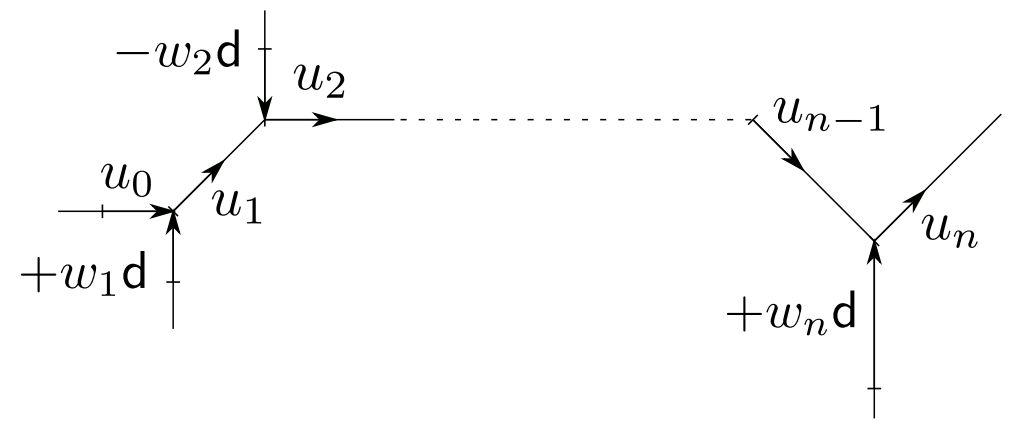

Figure 21: Slopes of line segments along a floor

Label the vertices on this floor from 1 to $n$, and for each let $w_i$ be the weight of the corresponding elevator, and $\varepsilon_i = +1$ (resp. $-1$) if it arrives from the bottom (resp. from the top). The balancing condition says that, for all $i = 0, \ldots, n$,

$$u_i = u_{i-1} + \varepsilon_i w_i \mathsf{d}. \tag{7.16.1}$$

By $\mathsf{d}^\perp$-transversality, one has $\det(u_0, \mathsf{d}) = \pm 1$. Then it follows from the induction relation (7.16.1) that $\det(u_i, \mathsf{d}) = \pm 1$ for all $i = 0, \ldots, n$, hence the vertex $\#i$ contributes to the multiplicity $\mu^{\mathbf{C}}(C, f)$ by $w_i$, cf. (6.6.1).

Thus, the multiplicity $\mu^{\mathbf{C}}(C, f)$ may be computed as a product of contributions from the elevators, the contribution of an elevator of weight $w$ being $w$ (resp. $w^2$) if it is attached to one (resp. two) floor. Therefore

$$\mu^{\mathbf{C}}(C, f) = \prod_{e \in \mathrm{Edges}^{\pm\infty}} w(e) \cdot \prod_{e \in \mathrm{Edges}^{\mathrm{fin}}} w(e)^2,$$

where the products on the right-hand-side are indexed on the edges of the floor diagram associated to $(C, f)$. This coincides with the multiplicity $\mu^{\mathbf{C}}(\mathcal{D}, m)$ after we multiply by the factor $I^{-\alpha}$ indicated by Mikhalkin's correspondence theorem. □

Applications of Corollary (7.16) to the enumeration of plane cubics with various tangency conditions along a fixed line have already been given in [IX, (1.7)]. Let us now give a few other examples, in continuation of Section 2.3.

**(7.17) Examples.** We continue with some of the examples from (7.5), (7.8), and (7.11).

*(7.17.1) Class of the quartic surface of (2.7).* This amounts to counting marked floor diagrams of type $\Delta$ and genus 0 (with trivial boundary conditions), where $\Delta$ is the square displayed in (2.7), with respect to the direction $\mathsf{d} = (0, 1)$. There is only one such diagram up to equivalence, which is the following.

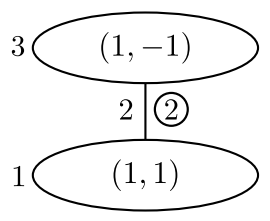

Figure 22: Floor diagram of genus 0 for the Newton polygon of (2.7)

Thus the number of irreducible rational curves in a pencil of hyperplane sections of the surface $(xy - t^2 = zw - t^2 = 0) \subseteq \mathbf{P}^4$ is 4, by Corollary (7.16).

*(7.17.2) Invariants of the octic surface of (2.5).* We want to count irreducible genus 1 (resp. genus 0) curves in a pencil (resp. a web) of hyperplane sections of the octic surface in $\mathbf{P}^7$ corresponding to the rectangle $\Delta$ displayed in (2.5). The corresponding marked floor diagrams, with respect to direction $\mathsf{d} = (0, 1)$, are the following.

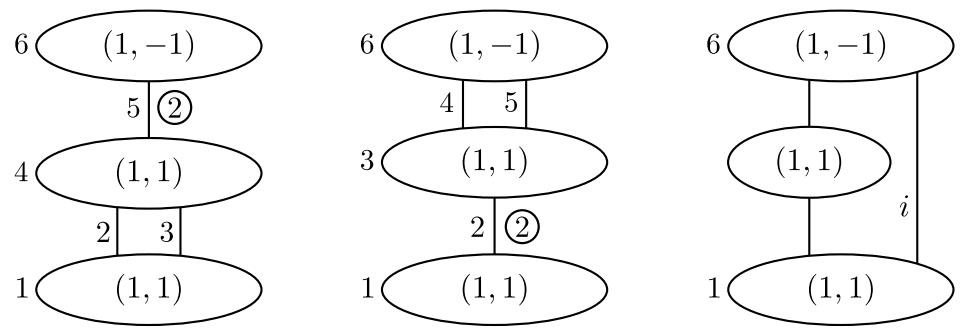

Figure 23: Genus 1 floor diagrams

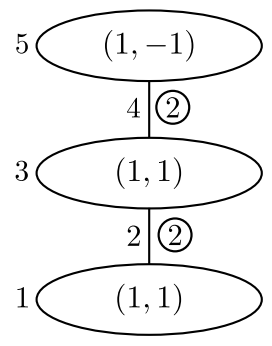

Figure 24: Genus 0 floor diagram

In the genus 1 case, the rightmost diagram has four possible markings, as $i$ may take any value in $[\![1, 4]\!]$, which then determines the rest of the marking (i.e., the marking of the one floor and two elevators with nothing indicated on the figure). We thus find $N_1^\Delta = 12$, and $N_0^\Delta = 16$.

Note that the number in genus 1 does not coincide with the one given in (2.5), which is not a contradiction because this is not the same enumerative problem: in the computation above, no base points have been imposed along the toric anticanonical cycle, whereas in (2.5) all the points along the toric anticanonical cycle are fixed. The latter computation cannot be performed purely in terms of floor diagrams, because we want too many boundary conditions with respect to what is possible, according to the condition given in (7.6). We solve the tropical problem corresponding to the enumerative problem of (2.5) in the next example.

**(7.18) Example.** So, let us consider the enumerative problem of (2.5): we let $\Delta$ be the Newton polygon indicated below, and fix $\alpha = (2, 1, 2, 0)$, $\beta = (0, 0, 0, 1)$. Thus the set $\Omega$ of boundary conditions consists of five real values, fixing five lines $\mathsf{L}_1, \mathsf{L}_1', \mathsf{L}_2, \mathsf{L}_3, \mathsf{L}_3'$ as indicated on the picture below.

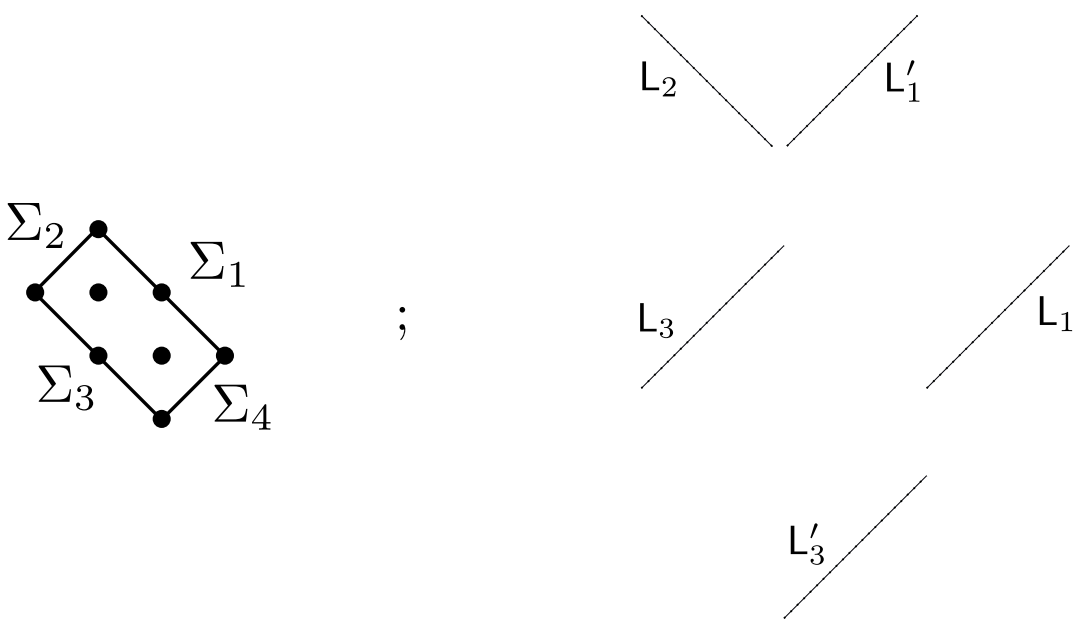

The tropical curves defined by this data move in a 2-dimensional family, in the way pictured below (three curves are pictured, in green, blue, and purple, respectively): the deformations of the curves amount to change the sizes of the two rectangles by letting their corners slide along the lines $\mathsf{L}_i$ and $\mathsf{L}'_i$. As the curves move in this family, one notices that the semi-line corresponding to the edge $\Sigma_4$ of the Newton polygon is always supported on the same line, which we shall call $\tilde{\mathsf{L}}_4$: this is yet another tropical incarnation of the toric Cayley–Bacharach property stated in Lemma (2.1).

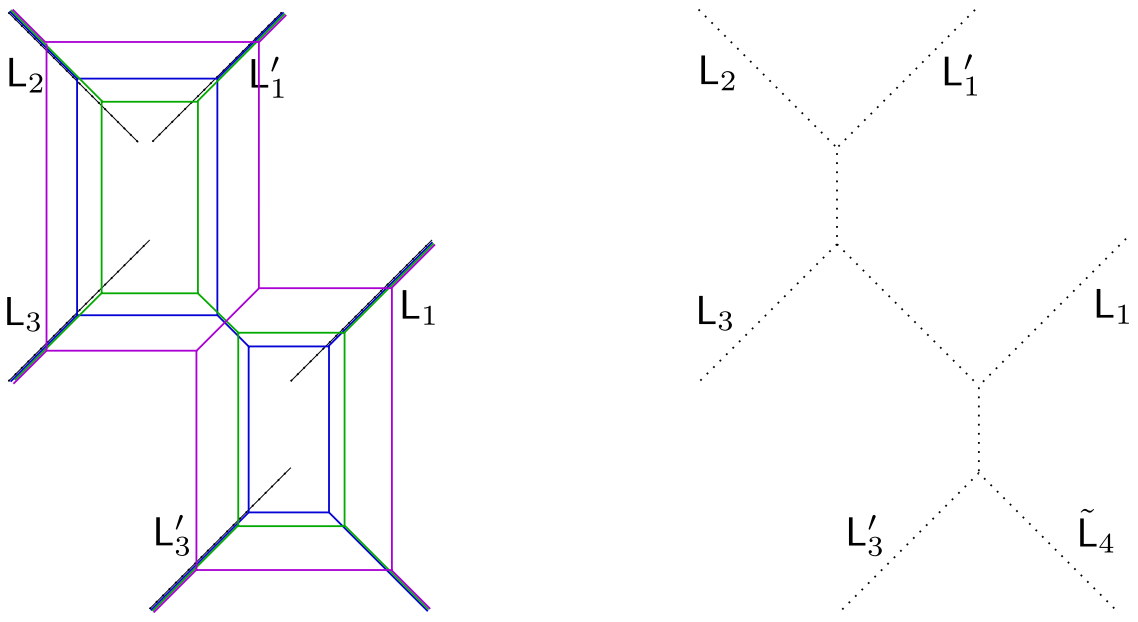

Figure 25: A pencil of tropical curves and the Cayley–Bacharach property

We want to count curves of genus $g = 1$ and $g = 0$ in this system. Let us focus on the case of genus $g = 1$. We need to fix a set $\underline{x}$ of base points consisting of one point only, which is represented by a black dot on Figure 26 below. With the configuration depicted on this figure, we find the four physical curves represented below: the first two count with multiplicity 1 each, and the last two with multiplicity 4 each. This confirms the number $N_1^{\Delta} = 10$ given in (2.5) and found in [10].

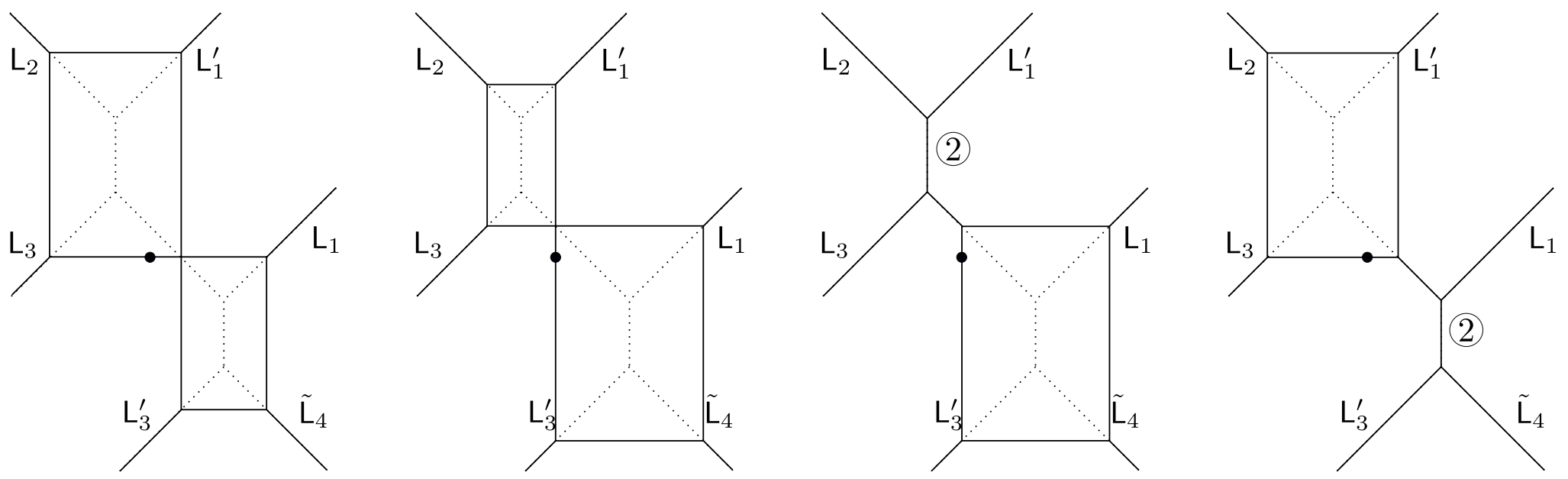

Figure 26: Genus 1 curves in the pencil defined by one base point

Note that these four curves all have a perfectly fine floor decomposition, with corresponding (unmarked) floor diagram one of the diagrams in Figure 23 above. It is however not true that the base points and boundary conditions are taken up exactly one by one by the floors and elevators.

To complete the picture, note that there is a unique genus $g = 0$ curve in this tropical linear system, which you can see by filling out the spaces between the dots in the picture on the right of Figure 25. It counts with multiplicity 16, and has the (unmarked) floor diagram given in Figure 24.

## 8 – Tropical representation of algebraic curves

In this section we elaborate on the idea that tropical curves are great for the graphic representation of algebraic curves in toric surfaces, and illustrate this with meaningful examples. A natural attempt to represent complex algebraic plane curves is to draw the set of their real points, in $\mathbf{R}^2$ say, but as anyone who has tried knows, this does not reflect the reality of things very well. For instance, the set of real points may be disconnected, it may contain isolated points, and it does not respect the Bezout Theorem. These problems do not occur with tropical curves.

**(8.1) Theorem** (Bezout Theorem, see [25, Section 4])**.** *Let $C, D \subseteq \mathbf{R}^2$ be two plane tropical curves with Newton polygon $\Delta$. If $C$ and $D$ intersect at finitely many points, then the number of their intersection points, counted with multiplicities, is $c_1(\mathcal{O}_{X_\Delta}(1))^2 = 2\mathscr{A}(\Delta)$, in the notation of (1.9).*

The intersection multiplicity of two line segments with primitive direction vectors $u$ and $v$ and multiplicities $m$ and $n$ is defined as $m \cdot n \cdot \big|\det(u, v)\big|$. Note that the theorem says not only that all intersection points are visible in the tropical shadow, but also that none of them has wandered at infinity to some component of the anticanonical cycle.

The next property is merely an observation, as it follows directly from the definition of the Newton polygon of a plane tropical curve. At a more conceptual level, plane tropical curves are by design projected shadows of curves in the torus $(\mathbf{C}^*)^2$; it is therefore expected that they behave well with respect to the completion of $(\mathbf{C}^*)^2$ in a toric surface by the addition of a toric boundary.

**(8.2) Remark.** Let $C \subseteq \mathbf{R}^2$ be a plane tropical curve with Newton polygon $\Delta$. The number, counted with multiplicities, of infinite branches of $C$ in the direction corresponding to a segment $\delta$ of the boundary of $\Delta$ equals the intersection number $c_1(\mathcal{O}_{X_\Delta}(1)) \cdot C_\delta$, where $C_\delta$ is the component of the anticanonical cycle $\partial_\Delta$ corresponding to $\delta$.

Now, the principles of "my" tropically inspired way of graphically representing algebraic curves on (surfaces birational to) toric surfaces are the following.

**(8.3) Graphic representation.** Let $C$ be a curve in a toric surface $X_\Delta$. We represent the surface $X_\Delta$ by drawing its anticanonical cycle $\partial_\Delta$ as a polygon (which looks much like the Newton polygon $\Delta$); then the surface $X_\Delta$ is the inside of this polygon. To represent the curve $C$, we draw a (finite length portion of a) tropical curve $C_{\mathrm{trop}}$ with the same Newton polygon as $C$, inside the polygon representing $X_\Delta$: the infinite branches of $C_{\mathrm{trop}}$ are represented by finite segments hitting the appropriate segment of $\partial_\Delta$, and the finite segments of $C_{\mathrm{trop}}$ do not hit $\partial_\Delta$. The concept of tropical toric varieties, see [5, Section 7.3], gives a mathematical meaning to this graphic process.

Note that the Newton polygon of $C$ may differ from $\Delta$: for instance, we may want to represent a fibre of the ruling of a Hirzebruch surface $\mathbf{F}_e$, or a curve $C$ in $\mathbf{F}_e$ such that the complete linear system $|C|$ represents $\mathbf{F}_e$ as a cone over a rational normal curve. In this case

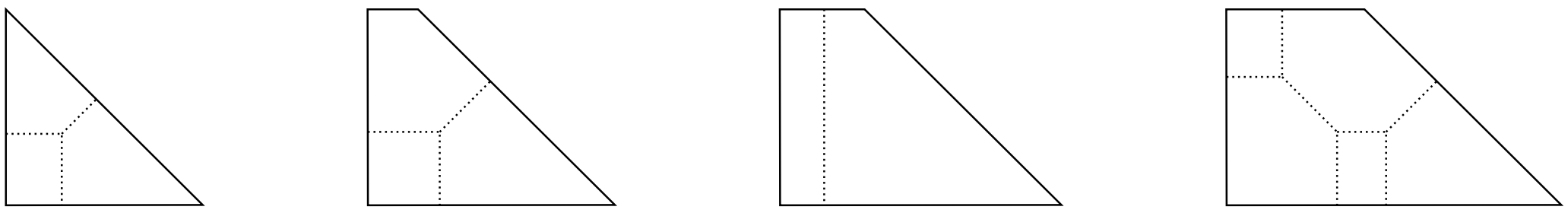

Figure 27: Curves linearly equivalent to $H$ in $\mathbf{P}^2$, $H$ in $\mathbf{F}_1$, $F$ in $\mathbf{F}_1$, and $H+F$ in $\mathbf{F}_1$ (in the notation introduced in (0.1))

the curve $C$ is not ample in $X_\Delta$, and some segments of $\reflectbox{K}_\Delta$ will not be hit by any infinite branch of $C_{\text{trop}}$. We give first examples in Figure 27.

One more nice feature of this way of picturing curves is that it lets graphically see the genus, as it is the first Betti number of the real object that we are actually seeing. Thus, in Figure 28, we see that the curve on the left (resp., on the right) has genus one (resp., two).

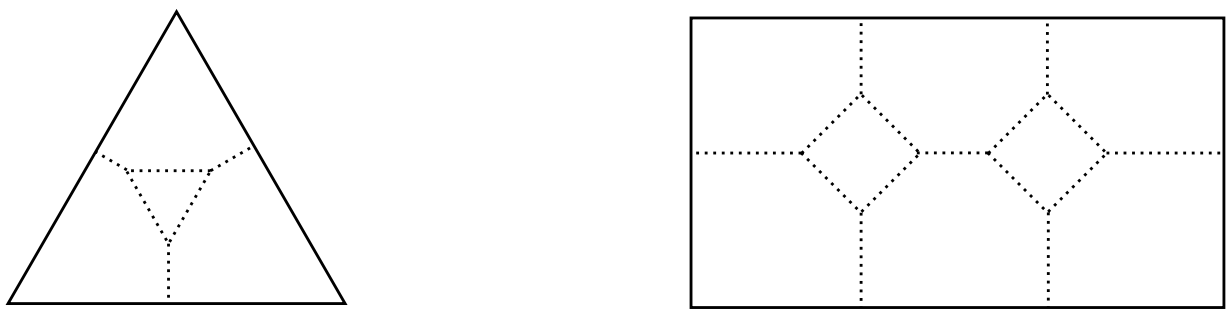

Figure 28: Hyperplane sections of the cubic of (2.4) and the quartic of (2.5)

Note that what we are considering is really representations of algebraic curves inspired by tropical geometry, and not actual tropical curves. In particular, we may forget about the integral structure underlying $\mathbf{R}^2$ (i.e., , $\mathbf{Z}^2$), the balancing condition, the subdivision of the Newton polygon of the curve, and so on. This allows for an essentialization of our drawings: we may for instance simply represent a Hirzebruch surface by a rectangle instead of a trapezium (or by whatever quadrilateral we fancy), and simplify the representations of curves as indicated on Figure 29.

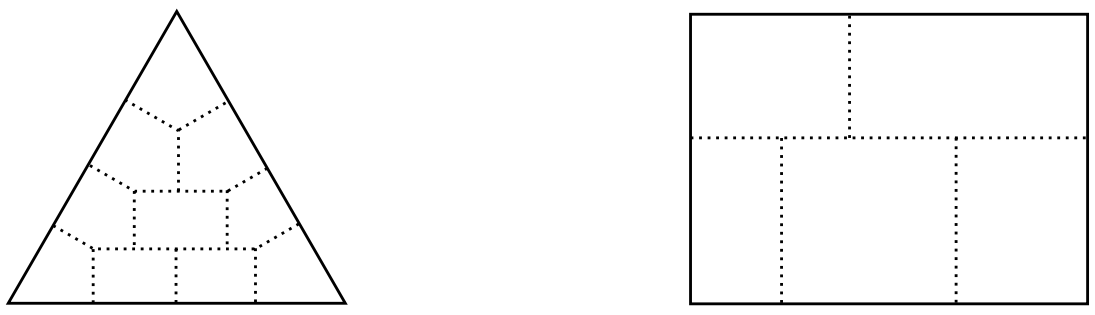

Figure 29: Curves linearly equivalent to $3H$ in $\mathbf{P}^2$, and $H+F$ in $\mathbf{F}_1$

In the same vein it frequently happens, as we have seen in [I], that one needs to represent curves that actually lie on some blow-up of a toric surface at points lying on the toric boundary. In that case we adapt our representation in a rather flexible way without much trouble.

Our interest resides especially in representing curves lying on reducible surfaces, the irreducible components of which are toric surfaces (or blow-ups thereof). In that case we may simply tile our polygons so as to reflect the dual graph of the various irreducible components, somehow as in (5.6). A nice example is given in [I, Figure $**$]. In that example, the dual graph of our reducible surface is a tiling of the sphere $\mathbf{S}^2$, and we have tweaked some of the sides of our polygons to obtain a better looking drawing; incidentally, we have also tweaked some line segments of the representation of our curve. Note that the surface $T^1$ is actually outside of its polygon, as we have performed a kind of stereographic projection from a point inside this surface $T^1$, in order to represent the sphere $\mathbf{S}^2$ on the plane. You may observe that the representation of the curve is topologically the same as the union of the edges of a tetrahedron, which is none other than the dual of the tetrahedron (in the sense of Section 2.3) we had started with.

**(8.4) Application to Brugallé degenerations** (see [IX])**.** I will simply give an example which, I believe, displays all I want to show. Let us follow a degeneration of a plane cubic along the degeneration of $\mathbf{P}^2$ to a chain of Hirzebruch $\mathbf{F}_1$ surfaces $S_0, \ldots, S_9$, with a (useless) projective plane $S_{10}$ at one end of the chain. The plane cubic degenerates (in our example) to a curve which restricts to curves linearly equivalent to $H$ on $S_9$, $H+F$ on $S_7$, $H+2F$ on $S_4$, and multiples of $F$ on all other $\mathbf{F}_1$ surfaces (and 0 on the plane), in the notation introduced in (0.1). We may see our plane cubic as a hyperplane section of the Veronese surface $v_3(\mathbf{P}^2) \subseteq \mathbf{P}^9$, and its degeneration as a hyperplane section of a "reducible toric surface" as in (5.6), made of a projective plane, attached to a rational normal scroll of type $(1,2)$ along a line, and the latter attached to a scroll of type $(2,3)$ along a conic.

In Figure 30a, we have drawn "proper tropical curves" inside the various polygons (those on which the curve has class $F$ are skipped). Now the wonder is that when we put all the pieces together, what we obtain is a plane tropical cubic! In Figure 30b, we have drawn "essentialized tropical curves"; if you look at their patching together with the right kind of eyes, you should see only floors and elevators.

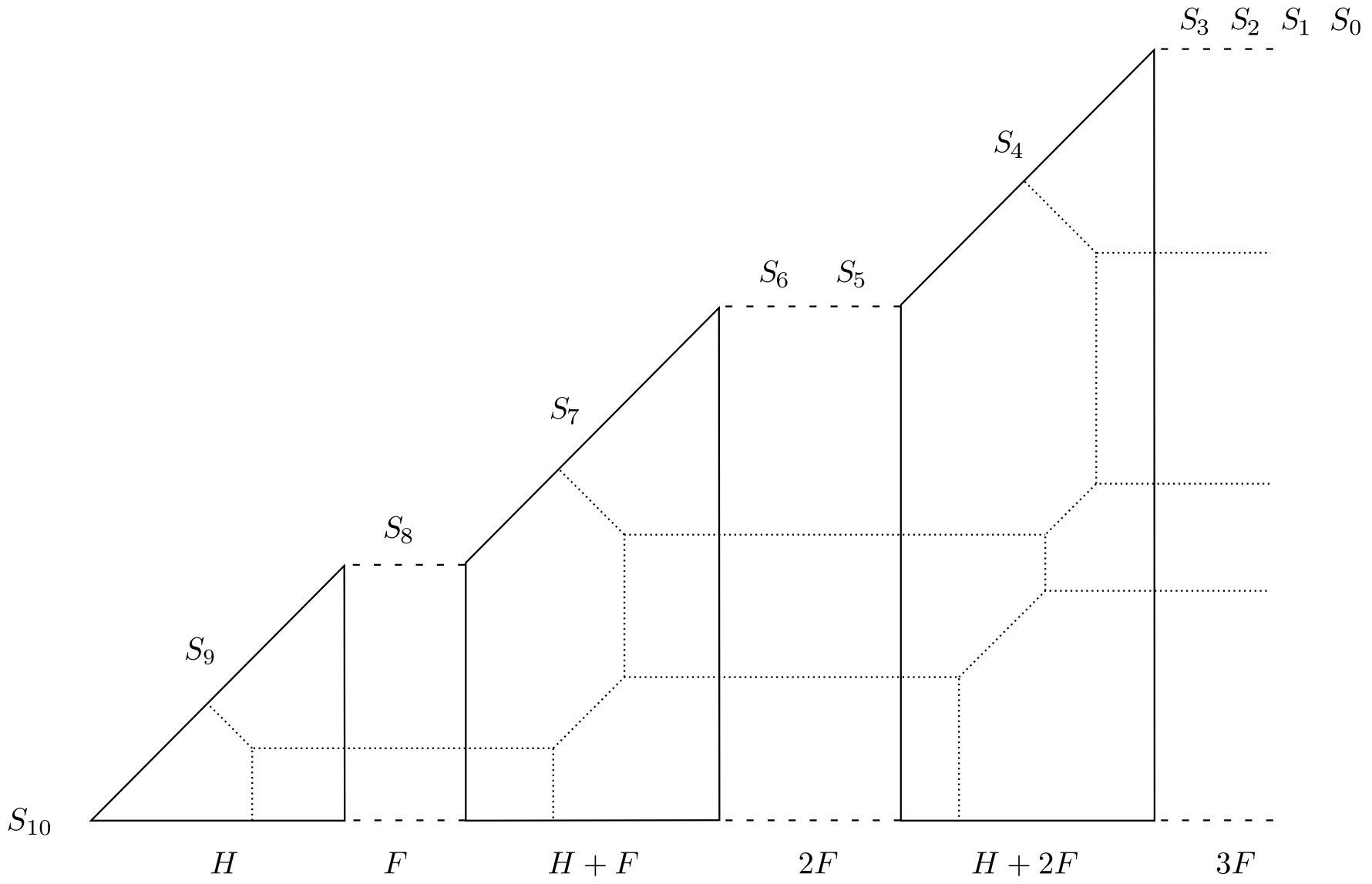

Figure 30a: Smooth plane cubic in a degeneration of $v_3(\mathbf{P}^2)$ to a union of scrolls (tropical version)

## References

### Chapters in this Volume[9]

[9]see `https://www.math.univ-toulouse.fr/~tdedieu/#enumTV`

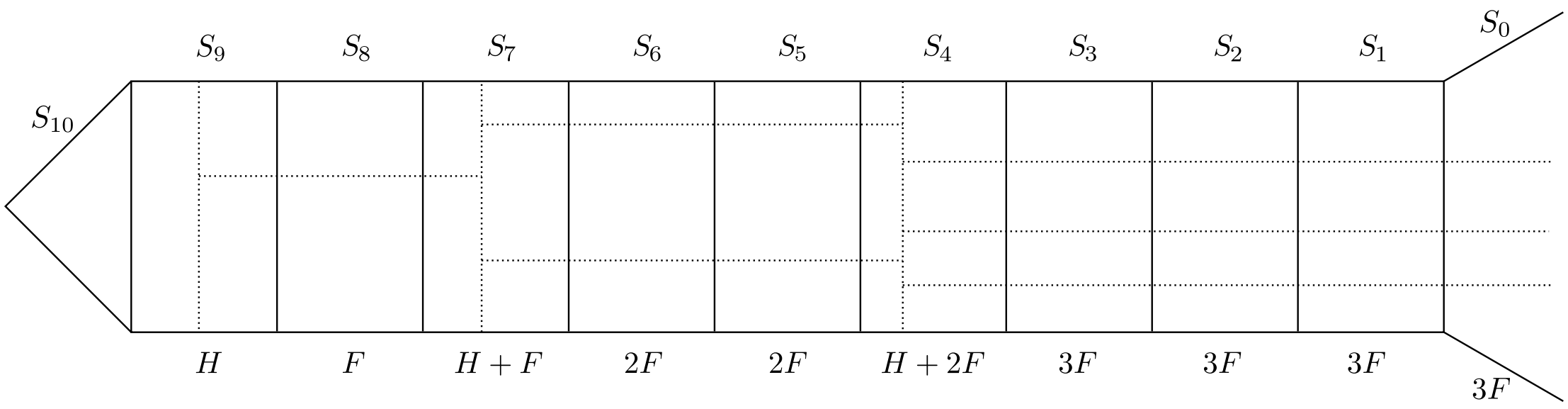

Figure 30b: Smooth plane cubic in a degeneration of $v_3(\mathbf{P}^2)$ to a union of scrolls (arts déco version)

**External References**